\documentclass[reqno]{amsart}
\usepackage{mathptmx,amsmath,mathbbol,amssymb,mathbbol,amsthm}

\usepackage{epic}
\usepackage{pstricks}
\usepackage{amsmath,amssymb,colordvi}
\usepackage{url}

 \usepackage[dvips,bottom=1.4in,right=1in,top=1in, left=1in]{geometry}

\newtheorem{theorem}{Theorem}[section]

\newtheorem{lemma}[theorem]{Lemma}
\newtheorem{corollary}[theorem]{Corollary}

\theoremstyle{definition}

\newtheorem{example}[theorem]{Example}
\newtheorem{remark}[theorem]{Remark}

\def\ud{\mathrm{d}}

\def\sign{\operatorname{sign}}

\def\RR{\mathbb{R}}
\def\R{\mathbb{R}}
\def\NN{\mathbb{N}}
\def\N{\mathbb{N}}

\def\st{\ :\ }
\def\Om{\Omega}

\def\var{\varepsilon}

\def\bOm{\bar{\Om}}

\def\pOm{\partial \Omega}
\def \sgn{\operatorname{sgn}}

\def\Cap{\operatorname{Cap}}

\def\cA{{\mathcal A}}

\def\E{{\mathcal E}}

\DeclareMathOperator{\Dom}{dom}

\DeclareMathOperator{\argmin}{arg\,min}
\newcommand{\dom}[1]{\Dom\,{#1}}

\newcommand\ringring[1]{%
  {
   \mathop{\kern0pt #1}\limits^{
     \vbox to-1.85ex{
       \kern-2ex 
       \hbox to 0pt{\hss\normalfont\kern.1em \r{}\kern-.45em \r{}\hss}%
       \vss 
     }
   }
  }
}

\numberwithin{equation}{section}

\newcounter{aufzi}
\renewcommand\theaufzi{(\alph{aufzi})}

\newcounter{aufzii}
\renewcommand\theaufzii{(\roman{aufzii})}

\newcounter{aufziii}
\renewcommand\theaufziii{(\arabic{aufziii})}

\newcounter{aufziv}
\renewcommand\theaufziv{\Roman{aufziv}}

\begin{document}

\title[Nonlinear, nonlocal  operators]{A unified view of nonlinear, nonlocal operators and qualitative properties of associated elliptic and parabolic problems}

\author{Ralph Chill}

\address{R.~Chill, Institut f\"ur Analysis, Fakult\"at Mathematik, Technische Universit\"at Dresden, 01062 Dresden (Germany)}
\email{ralph.chill@tu-dresden.de}

\author{Mahamadi Warma}
\address{M.~Warma,  Department of Mathematical Sciences and The Center for Mathematics and Artificial Intelligence, George Mason University, Fairfax, VA 22030, USA}
\email{mwarma@gmu.edu}



\thanks{The work of MW is partially supported by the US Army Research Office (ARO) under Award NO: W911NF-20-1-0115. This work started when MW was visiting the Institute of Analysis at the Technische Universit\"at Dresden (Germany). He would like to thank the Institute of Analysis for their support.}

\begin{abstract}
We put together a general framework to deal with elliptic and parabolic equations associated with (nonlinear) nonlocal (fractional order) operators. Many well-known nonlocal operators enter into our framework, and in addition one may introduce many other, new nonlocal operators that have not yet been considered in the literature. We use the abstract theory of (nonlinear) semigroups generated by subgradients of proper, lower semicontinuous and convex functionals on Hilbert spaces to build a rigorous and applicable framework that works for many classical elliptic operators but also nonlocal or sometimes fractional operators. After recalling the notion of a nonlinear semigroup generated by subgradients and $j$-subgradients of the associated energy functions, we introduce a general class of (nonlinear) nonlocal elliptic type operators
and define rigorously subgradients and $j$-subgradients of such functionals that generate (nonlinear) submarkovian semigroups and hence, the abstract Cauchy problem associated with these subgradients and/or $j$-subgradients is wellposed. The existence and the qualitative properties of solutions to these Cauchy problems and the corresponding semigroups are investigated. More precisely, we show some comparison and maximum principles, submarkovian, domination, ultracontractivity properties, and some H\"older type estimates for these semigroups of operators. These results are usually useful in several branches of pure and applied partial differential equations. We finish the paper by giving several examples of nonlocal operators in Euclidean spaces, graphs, metric random walk spaces, fractional Brownian motions, and L\'evy flights, that fit in our general framework. 
\end{abstract}

\maketitle

\tableofcontents

\section{Introduction}

We study the nonlocal evolution problem
\begin{equation} \label{eq.nonlocal}
 \begin{cases}
 \displaystyle \partial_t \hat{u} (t,x) + P.V.\int_{\hat{\Omega}} \varphi (x,y,\hat{u}(t,x)-\hat{u}(t,y)) \; \ud \mu(y) = f(t,x) & \quad (t,x) \in ]0,\infty [\times\Omega , \\
\displaystyle  P.V.\int_{\hat{\Omega}} \varphi (x,y,\hat{u}(t,x)-\hat{u}(t,y)) \; \ud \mu(y) = 0 & \quad (t,x) \in ]0,\infty [ \times (\hat{\Omega} \setminus \Omega ) , \\
 \displaystyle \text{possibly boundary conditions} &\quad  \\ 
 \displaystyle \hat{u}(0,x) = u_0 (x) & \quad x\in\Omega .
 \end{cases}
\end{equation}
Here, $(\hat{\Omega} , \mathcal{A},\mu)$ is a measure space, and $\Omega \subseteq\hat{\Omega}$ is a measurable subset. The kernel function $\varphi : \hat{\Omega} \times \hat{\Omega} \times \R \to \R$, the right hand side $f : [0,\infty [ \times \Omega \to \R$ and the initial value $u_0 :\Omega \to\R$ are given, the function $\hat{u} : [0,\infty [\times \hat{\Omega} \to\R$ is the unknown, and P.V. stands for Cauchy principal value. The evolution problem \eqref{eq.nonlocal} actually is a system of a general nonlinear evolution equation in $\Omega$ (first line) coupled with a stationary equation in $\hat{\Omega}\setminus\Omega$ (second line), which we call the exterior condition, and a boundary condition in $\partial\hat{\Omega}$ or possibly some other abstract boundary which depends on the kernel $\varphi$. The exterior condition is due to the nonlocal nature of the operator involved in the first line of \eqref{eq.nonlocal} that we shall explain below. The boundary condition will be made explicit in Sections \ref{sec.neumann}, \ref{sec.robin} and \ref{section-ER}, at least for some special cases.

A nonlocal partial differential equation (PDE) often is an integro-differential equation involving (singular) integral terms or pseudo-differential terms. 
Nonlocal models differ from the classical PDE models. In the classical PDEs interactions between two domains occur only due to contact. In the case of nonlocal PDEs interactions can occur at a distance. These interactions also explain the exterior conditions given in the second line of \eqref{eq.nonlocal}.

In recent years, there has been a growing interest in nonlocal models because of their relevance in a wide spectrum of practical applications. Indeed, there is a plethora of situations in which a nonlocal equation gives a significantly better description than a local PDE of the problem one wants to analyze. A widely studied class of nonlocal models involves fractional order operators, which have nowadays emerged as a modeling alternative in various branches of science. They usually describe anomalous diffusion. Typical examples in which nonlocal or fractional equations appear are models in turbulence (Bakunin \cite{bakunin2008}), population dynamics (De Roos \& Persson \cite{de2002}), image processing (Gilboa \& Osher \cite{gilboa2008}), laser design (Longhi \cite{longhi2015}), and porous media flow (Vazquez \cite{vazquez2012}). In addition, a number of stochastic models associated with fractional operators have been introduced in the literature to explain anomalous diffusion. Among them we quote the fractional Brownian motion, the continuous time random walk, the L\'evy flights, the Schneider gray Brownian motion, and more generally, random walk models based on evolution equations of single and distributed fractional order in  space (see e.g. Dubkov, Spagnolo \& Uchaikin \cite{dubkov2008}, Gorenflo, Mainardi \& Vivoli \cite{gorenflo2007}, Mandelbrot \& Van Ness \cite{mandelbrot1968}, Schneider \cite{schneider1990}). In general, a fractional diffusion operator corresponds to a diverging jump length variance in the random walk. Finally, we can refer to Antil \& Rautenberg \cite{antil2019} and Weiss, van Bloemen \& Antil \cite{weiss2018} for the relevance of fractional operators in geophysics and imaging science.

In the wide spectrum of nonlocal and fractional models, of particular interest are those involving the fractional Laplacian. From a mathematical perspective, there is today a well established and rich literature on this operator and its employment in PDE models. Among many other contributions, we recall here the works of Biccari \& Hernandez-Santamari\'a \cite{biccari2018-1}, Biccari, Warma \& Zuazua \cite{biccari2017,biccari2018-2}, Caffarelli \& Silvestre \cite{CaSi07,caffarelli2009}, Leonori et al. \cite{leonori2015}, Ros-Oton \& Serra \cite{ROSe14}, Servadei \& Valdinoci \cite{servadei2012,servadei2013}, Warma \cite{Wa15}. Also, it is well known that the fractional Laplacian is the generator of s-stable L\'evy processes, and it is often used in stochastic models with applications, for instance, in mathematical finance (Levendorskii \cite{levendorskii2004}, Pham \cite{pham1997}).

The growing interest in nonlocal and fractional models has opened a very challenging field in applied mathematical research, since many of the existing techniques in the analysis of PDEs were not adapted to treat nonlocal effects. Unlike classical PDEs, where the underlying dynamics is governed by local interactions and differential operators, nonlocal equations often require tailored analytical approaches. In fact, many of the standard tools in the theory of PDEs either fail to apply or must be significantly modified to account for the global influence embedded in nonlocal operators. 

If $\Omega\subseteq\hat{\Omega}\subseteq\R^N$ are open, $\theta\in ]0,1[$, and 
\begin{equation}\label{phi1}
\varphi (x,y,s) = \frac{s}{|x-y|^{N+2\theta}}, 
\end{equation}
then the problem \eqref{eq.nonlocal} is a linear evolution problem with a nonlocal operator. If also $\Omega=\R^N$ the associated operator is well-known in the literature as the fractional Laplace operator (see Example \ref{flrn} for the precise description) and if $\Omega=\hat\Omega$ is a proper subset of $\R^N$, then it is known as the regional fractional Laplace operator (see e.g. Example \ref{ex.9.3} and Gal \& Warma \cite{GaWa2021}, Guan \cite{Gua}, Guan \& Ma \cite{GuMa1,GuMa2}, Warma \cite{Wa15,War-NODEA,War18} for the precise description).
If $\Omega\subsetneq\hat{\Omega}$,  then there appears the stationary equation in \eqref{eq.nonlocal}, the so-called {\em Neumann exterior condition} (see Section \ref{sec.neumann}).  Below in Section \ref{sec.robin}, we shall also consider Robin and Dirichlet exterior or complementary conditions, but we shall also study situations when $\hat{\Omega}$ is not a subset of $\R^N$, for example a graph or a random walk space. In this direction and if $\hat{\Omega} =\R^N$, a great progress has been made regarding the existence and regularity of solutions of the linear evolution equation \eqref{eq.nonlocal} (see e.g. Caffarelli \& Silvestre \cite{CaSi07}, Chen, Kim \& Song \cite{ChKiSo10}, Chen \& Kumagai \cite{ChKu08}, Claus \& Warma \cite{ClWa20}, Gal \& Warma \cite{Ga-Wa-CPDE,KSW} and the references therein)
and the associated stationary (elliptic) problem (see e.g. Abels \& Grubb \cite{abels23}, Caffarelli \& Stinga \cite{CaSt}, Di Nezza, Palatucci \& Valdinoci \cite{DNPaVa12}, Dipierro, Ros-Oton \& Valdinoci \cite{DiROVa17}, Gal \& Warma \cite{Ga-Wa-CPDE}, Grubb \cite{G-JFA,grubb15}, Guan \cite{Gua}, Guan \& Ma \cite{GuMa1,GuMa2} and their references).

If $\Omega\subseteq\hat{\Omega}\subseteq\R^N$ are open, $p\in ]1,\infty[$, $\theta\in ]0,1[$, and 
\begin{equation}\label{phi2}
\varphi (x,y,s) = \frac{|s|^{p-2}s}{|x-y|^{N+p\theta}},
\end{equation}
then the problem \eqref{eq.nonlocal} is a quasi-linear evolution problem with a quasi-linear nonlocal operator. If $\Omega=\R^N$ the associated operator is known in the literature as the fractional $p$-Laplace operator (see e.g. Section \ref{sec.neumann}) and if $\Omega=\hat\Omega$ is a proper subset of $\R^N$, then it is known as the regional fractional $p$-Laplace operator (see Section \ref{sec.lnnod} and Gal \& Warma \cite{Ga-Wa-CPDE}, Warma \cite{War-NODEA,War18} and the references therein). The case of general $\Omega\subseteq\hat{\Omega}=\R^N$ has been considered in Foghem \cite{Fo25,Fo20,Fo25b}; see also Foghem \& Kassmann \cite{FoKa24} for the linear case.

The success obtained by researchers regarding the preceding special cases is due to the fact that the function spaces needed to study these problems are fractional order Sobolev spaces whose well-known properties allow one to obtain some satisfactory results. But several well-known results in the local case are still unknown in the nonlocal case, specifically the numerical analysis and precise elliptic and parabolic regularity of solutions of such PDEs, and more importantly the controllability and observability properties of such PDEs. We refer to the recent papers by Biccari, Warma \& Zuazua \cite{BWZ-HB,BWZ-BC} and their references for a precise discussion on this topic. 

In the present work our aim is to present a unified, general framework for the system \eqref{eq.nonlocal} that includes the above mentioned examples but also allows one to think of new examples, which include different sets of points $(x,y)\in\hat{\Omega}\times\hat{\Omega}$  between which a nonlocal interaction is possible (see the sets ${\mathbf S}_\Phi$ defined in Section \ref{sec.neumann}). Also, the kernels $\varphi$ need not depend on the distance between $x$ and $y$, and in the dependence on the third variable may be a more general monotone function. And last but not least, the choice of the sets $\Omega$ and $\hat{\Omega}$ is relatively free; $\hat{\Omega}$ need not be equal to $\R^N$, $\Omega$ and $\hat{\Omega}$ need not be subsets of $\R^N$. Of course, there appears the problem to identify for example the boundary conditions on $\partial\hat\Omega$ or some other subset of the closure of $\hat{\Omega}$. Boundary conditions depend on an appropriate integration by parts formula. This is crucial in the study of elliptic and parabolic PDEs. We can mention that in this direction even in the case where $\hat\Omega=\R^N$ and $\varphi$ is given by \eqref{phi1} or \eqref{phi2},  nowadays the integration by parts formulas for fractional PDEs are satisfactory for smooth functions (see, for example, Dipierro, Ros-Oton \& Valdinoci \cite{DiROVa17}, Guan \cite{Gua}, Warma \cite{War-NODEA} and their references). These satisfactory integration by parts formulas are based on Fourier transforms (in the linear case) and/or on spaces of smooth functions (see, for example, Di Nezza, Palatucci \& Valdinoci \cite{DNPaVa12}, Dipierro, Ros-Oton \& Valdinoci \cite{DiROVa17}, Grubb \cite{G-JFA} and their references). Unfortunately, this is not always useful in the weak formulation of fractional PDEs on domains since most of the time we lack a density argument. In addition, it is known that the exterior conditions considered in the present paper and elsewhere are always satisfied almost everywhere in $\hat\Omega\setminus\Omega$ (see e.g. Claus \& Warma \cite{ClWa20} and their references), but the huge problem is the equation in $\Omega$. Even if such equation is satisfied in the sense of distributions in $\Omega$,  in contrast with the classical local PDEs,  even if the right hand side belongs to $L_{loc}^1(\Omega)$, we cannot immediately conclude that we have a strong solution. The classical theory of subgradients on Hilbert spaces, dating back at least to the 1960s/70s proves to be very successful here. We shall give in Section \ref{sec.neumann} a rigorous weak formulation of the associated stationary problem that works for our general setting. These results allow us to obtain several qualitative properties of solutions of the system \eqref{eq.nonlocal} and the associated stationary system.

In order to study well-posedness and qualitative behavior of \eqref{eq.nonlocal} we write the system \eqref{eq.nonlocal} as an abstract gradient system on the Hilbert space $L^2(\Omega)$. Well-posedness then follows from the classical theory developed for example in the monographs by Brezis \cite{Br73}, Barbu \cite{Bar10}, Miyadera \cite{Miya-1992} and others. There is a little twist in the correct formulation of the abstract gradient system, as the system \eqref{eq.nonlocal} actually is a system on the possibly larger set $\hat{\Omega}$, while we get well-posedness only on $L^2 (\Omega )$. This obstacle is overcome by the recent theory of $j$-subgradients.

The rest of the paper is organized as follows.  In Section \ref{sec.preliminaries} we recall some known results on subgradients and $j$-subgradients of proper, lower semicontinuous convex functionals and the properties of the semigroups generated by these negative subgradients. Our work really starts in Section \ref{sec.neumann}. First, we introduce the general assumptions on the function $\varphi$ and the measure space $(\hat{\Omega} , \mathcal{A},\mu)$. Under these assumptions we introduce Musielak-Orlicz-Sobolev type spaces and show that the system \eqref{eq.nonlocal} can be rewritten as an abstract Cauchy problem 
where the associated operator is the $j$-subgradient of a suitable proper, lower semicontinuous and convex energy function. We show that the associated Cauchy problem, hence the system \eqref{eq.nonlocal} is well-posed (see Theorem \ref{thm.e.lsc} and Corollary \ref{cor.neumann}). In the second part of this section we state and prove the above mentioned 
weak formulation of the stationary problem
(Theorem \ref{thm.identification.2}) and give the relation between strong and weak solutions (Lemma \ref{lem-3.4}). We also give examples of concrete nonlocal operators on open subsets of $\RR^N$ and on graphs. In Section \ref{NDF} we show that the semigroup generated by the negative $j$-subgradient we constructed is order preserving and $L^\infty$-contractive on $L^2 (\Omega )$, that is, the associated energy function is a nonlinear Dirichlet form (Theorem \ref{theo-subm}). These results give a comparison and maximum principles. With the help of the classical notion of regular Dirichlet forms we introduce a general notion of capacity associated to our Musielak-Orlicz-Sobolev type spaces and define the associated polar sets and quasi-continuous functions in Section \ref{sec5-cap}. Section \ref{sec.robin} is devoted to other boundary and exterior conditions. More precisely, first we give some admissible potentials that help us to define consistently the Robin boundary and exterior conditions. Second, we give some generation results of submarkovian semigroups of the associated negative $j$-subgradients, where the main results are Theorems \ref{Sub} and \ref{thm.identification.b.2}. Section \ref{sec-dom} is devoted to the domination of (nonlinear) semigroups. We introduce the notions of dominated and totally dominated (nonlinear) semigroups that are characterized in terms of the associated functionals. In particular we show that under suitable conditions on a given function $B$ and a Borel measure $\nu$ on the closure of $\hat{\Omega}$, the Dirichlet semigroup $S^{D,\nu}$ is totally dominated by the Robin semigroup $S^{B,\nu}$ which is also totally dominated by the Neumann semigroup $S^N$; see Theorem \ref{theo-7.1} and Remark \ref{rem-to-dom}. In Section \ref{sec-UlDN} we show that the the semigroups generated by the Dirichlet and Neumann boundary and exterior conditions are ultracontractive in the sense that they map $L^q(\Omega)$ ($q\in [1,\infty]$) into $L^\infty(\Omega)$. More precisely we give an $L^q-L^\infty$-H\"older continuity estimates of the semigroups generated by the Dirichlet and Neumann boundary and exterior conditions; see Theorems \ref{theo-ultra} and \ref{theo-ultra-NEC}. Finally in Section \ref{section-ER} we give further concrete examples that enter in our general framework. To be more precise, we give examples of fractional powers of general elliptic operators with various boundary conditions, we introduce other nonlocal operators on domains in $\R^N$ and their variants on graphs, and we briefly touch metric random walk spaces. 

\section{Semigroups generated by negative subgradients} \label{sec.preliminaries}

In order to solve the nonlocal problem \eqref{eq.nonlocal} for initial values in $L^2 (\Omega )$, we use the classical theory of gradient systems in Hilbert spaces, which we briefly recall. 

Let $H$ be a real Hilbert space with scalar product $\langle\cdot,\cdot\rangle_H$ and norm  $\|\cdot\|_H$, let $\E : H \to ]-\infty , \infty ]$ be a proper, lower semicontinuous, convex function with {\em effective domain} $\dom{\E} = \{u\in H:\; \E(u) <\infty \}$, and let 
\begin{align*}
\partial\E := \{ (u,f) \in H \times H \st u\in {\rm dom}\,\E \text{ and }  
  \E (u+v) - \E (u) \geq \langle f , v\rangle_H \; \forall\; v\in H\}
\end{align*}
be its {\em subgradient} with {\em domain} $\dom{\partial\E} := \{ u\in H\st \exists f\in H \text{ s.t. } (u,f)\in\partial\E\}$. By Brezis \cite[Exemple 2.3.4]{Br73}, the subgradient is $m$-accretive ($=$ maximal monotone), and by \cite[Proposition 2.11]{Br73}, $\overline{\dom{\partial\E}}^H = \overline{\dom{\E}}^H$. Therefore, by  \cite[Th\'eor\`eme 3.4]{Br73}, for every initial value $u_0 \in D := \overline{\dom{\E}}^H$ and for every right hand side $f\in L^1_{loc} ([0,\infty[;H)$ the Cauchy problem
\begin{equation}\label{CP}
\dot u + \partial\E (u) \ni f \text{ on } [0,\infty [ , \qquad u(0) = u_0 ,
\end{equation}
admits a unique mild solution $u\in C([0,\infty[;H)$. In particular, the negative subgradient $-\partial\E$ generates a strongly continuous contraction semigroup $S = (S(t))_{t\geq 0}$ on $D$ (see  e.g. \cite[Th\'eor\`eme 3.1, Remarque 3.2]{Br73}). This means that there exist (nonlinear) contractions $S(t) \in C( D , D )$ ($t\geq 0$) such that 
\[
S(0)={\rm Id}_{D} \text{ and } S(t+s) = S(t)\circ S(s)\;  (t, s\geq 0) ,
\]
 and for every initial value $u_0\in D$ the function $u := S(\cdot )u_0 \in C([0,\infty [ ;H)$ is the unique mild solution of the Cauchy problem \eqref{CP} with right hand side $f=0$. By \cite[Th\'eor\`eme 3.6]{Br73}, the semigroup $S$ is regularizing in the sense that for every initial value $u_0\in D$ the solution is a strong solution on $]0,\infty [$, that is, $u\in C([0,\infty [ ; H) \cap W^{1}_{2,loc} (]0,\infty [ ; H)$, $u(t) \in \dom{\partial\E}$ for almost every $t\geq 0$, and the differential inclusion in \eqref{CP} is satisfied for almost every $t\geq 0$. See also Barbu \cite{Bar10} or Miyadera \cite{Miya-1992} for this classical theory. \\

As an extension of the theory of subgradients, we recall from Chill, Hauer \& Kennedy \cite{ChHaKe16} the notion of $j$-subgradients. Let $V$ be a locally convex vector space, $H$ a Hilbert space, let $\E : V \to [0,\infty ]$ be a proper, lower semicontinuous, convex function, and let $j: V\to H$ be a bounded, linear operator. We define the {\em $j$-subgradient} of $\E$ by
\begin{equation}\label{eq-jsubgradient}
  \partial_j\E := \left\{ (u,f) \in H \times H \st  \exists \hat{u}\in \dom{\E}\subseteq V \text{ s.t. } j(\hat{u}) = u \text{ and }
   \E (\hat{u} + \hat{v} ) - \E (\hat{u} ) \geq \langle f , j(\hat{v}) \rangle_H  \;\forall\; \hat{v}\in V \right\} .
\end{equation}
By the definition, the $j$-subgradient of the energy (defined on the space $V$) is an operator on the Hilbert space $H$. We say that the energy $\E$ is {\em $j$-elliptic} if there exists $\omega\in\R$ such that, for every $c\in\R$, the sublevel set 
\[
\left\{ \hat{u}\in  V \st \E (\hat{u}) + \frac{\omega}{2} \, \| j(\hat{u})\|_{H}^2 \leq c\right\}
\]
is weakly compact. By \cite[Theorem 2.6]{ChHaKe16}, if $\E$ is proper, lower semicontinuous, convex and $j$-elliptic, then the $j$-subgradient of the energy $\E$ is $m$-accretive ($=$ maximal monotone), but even more is true. By \cite[Corollary 2.7]{ChHaKe16}, there exists a proper, lower semicontinuous, convex energy $\E^{H} :H \to [0,\infty ]$ such that $\partial_j \E = \partial\E^{H}$; here, the latter operator is the classical subgradient defined above. In this sense, $j$-subgradients are just subgradients, but the ``generalization'' is a more flexible approach to subgradients. By the above mentioned classical theory, $-\partial_j \E$ 
generates a nonlinear strongly continuous contraction semigroup $S$ on $\overline{\dom{\E^H}}^H$. The energy $\E^{H}$ is unique up to an additive constant. By \cite[Theorem 2.9]{ChHaKe16}, and up to an additive constant,
\begin{equation} \label{eq.energy.h}
\E^{H} (u) = \inf_{j(\hat{u}) = u} \E (\hat{u} ) .
\end{equation}
It follows from the representation of $\E^H$ that $\dom{\E^H} = j (\dom{\E})$, and therefore the semigroup $S$ generated by $-\partial_j\E$ is defined on $\overline{j(\dom{\E})}^H$.

For given $(u,f)\in\partial_j\E$, every element $\hat{u}\in V$ such that $j(\hat{u}) = u$ and 
\[
\E (\hat{u} + \hat{v}) - \E (\hat{u}) \geq \langle f,j(\hat{v})\rangle_H \quad \text{ for every } \hat{v}\in V 
\]
is called {\em an elliptic extension} of $u\in H$. By definition of the $j$-subgradient, for every $u\in\dom{\partial_j\E}$ an elliptic extension always exists, and is given by
\[
\hat{u} = \argmin_{j(\hat{w}) = u} (\E (\hat{w}) -\langle f,j(\hat{w})\rangle_H ) = \argmin_{j(\hat{w}) = u} \E (\hat{w})  .
\]

\section{Realization of nonlocal elliptic operators as subgradients on $L^2 (\Omega )$} \label{sec.neumann}

Let $(\hat{\Omega} ,\cA , \mu )$ be a $\sigma$-finite 
measure space, and let $\Omega\subseteq \hat{\Omega}$ be a measurable subset. 

Fix a function $\Phi : \hat{\Omega} \times \hat{\Omega} \times \R \to [0,\infty ]$ such that
\begin{equation}  \label{cond.phi}
\begin{cases}
\Phi (\cdot , \cdot , s) \text{ is measurable on } \hat{\Omega}\times\hat{\Omega} \text{ for every } s\in\R, \\
 \Phi (x,y,0 ) = 0 \text{ for every } x, y\in \hat{\Omega } , \\
 \Phi (x,y,\cdot ) \text{ is lower semicontinuous and convex for every } x, y\in\hat{\Omega} , \\
 \Phi (x,y,\cdot ) \text{ is even for every } x, y\in\hat{\Omega} , \text{ and} \\
 \Phi (x,y,\cdot ) = \Phi (y,x,\cdot ) \text{ for every }  x, y\in\hat{\Omega} .
\end{cases}
\end{equation}
Given such a function, we define the {\em Fenchel conjugate} $\Phi^* : \hat{\Omega}\times\hat{\Omega}\times\R \to [0,\infty ]$ of $\Phi$ by
\[
 \Phi^* (x,y,t) := \sup \{ st -\Phi (x,y,s) \st s\in\R \}, \quad ((x,y,t)\in \hat{\Omega}\times\hat{\Omega}\times\R ) .
\]
We say that $\Phi$ satisfies the {\em $\Delta_2$-condition} if there exist a constant $K_{\Delta_2}>0$ and an integrable, positive function $h_{\Delta_2}:\hat{\Omega}\times\hat{\Omega} \to [0,\infty [$ such that, for every $(x,y,s)\in \hat{\Omega}\times\hat{\Omega}\times\R$,
\[
 \Phi (x,y,2s) \leq K_{\Delta_2} \, \Phi (x,y,s) + h_{\Delta_2} (x,y) .
\]

In the following, we usually assume that both the function $\Phi$ and its Fenchel conjugate $\Phi^*$ satisfy the $\Delta_2$-condition. 
One can easily compute that if the function $\Phi (x,y,\cdot )$ is constant $0$, then the Fenchel conjugate $\Phi^* (x,y,t) = \infty$ for $t\not= 0$, and vice versa: if $\Phi (x,y,s) = \infty$ for $s\not= 0$, then the Fenchel conjugate $\Phi^*(x,y,\cdot )$ is constant $0$. These two extremal functions are explicitly allowed; they have no influence on the $\Delta_2$-condition in the sense that they satisfy the inequality above for any constant $K_{\Delta_2} >0$ and $h_{\Delta_2} (x,y) = 0$. Under the assumption that $\Phi$ and $\Phi^*$ both satisfy the $\Delta_2$-condition, these extremal functions (at least for $(x,y)\in\hat{\Omega}\times\hat{\Omega}$ in a set of positive measure) are the only cases in which $\Phi (x,y,\cdot )$ or $\Phi^* (x,y,\cdot )$ can take the values $\infty$ somewhere. We define the sets
\begin{align*}
 {\mathbf N}_\Phi & := \{ (x,y)\in\hat{\Omega}\times\hat{\Omega} \st \Phi (x,y,s ) = 0 \text{ for every } s\in\R\} , \\
 {\mathbf C}_\Phi & := \{ (x,y)\in\hat{\Omega}\times\hat{\Omega} \st \Phi (x,y,s) = \infty \text{ for every } s\in\R\setminus\{ 0\} \}, \text{ and} \\
 {\mathbf S}_\Phi & := (\hat{\Omega}\times\hat{\Omega})\setminus {\mathbf N}_\Phi .
\end{align*}

By the first four lines in condition \eqref{cond.phi}, the function $\Phi$ is a continuous and convex pseudomodular (see e.g. Musielak \cite[p.1]{Mus83}). The corresponding {\em Musielak-Orlicz space} 
\[
 L^\Phi (\hat{\Omega}\times\hat{\Omega} ) :=\left \{ w \in L^0 (\hat{\Omega}\times\hat{\Omega}) \st \exists \lambda >0 : \int_{\hat{\Omega}} \int_{\hat{\Omega}} \Phi \left(x,y,\frac{w(x,y)}{\lambda}\right) \; \ud\mu (y) \; \ud\mu (x)<\infty \right\} 
\]
is a complete seminormed space for the natural Minkowski seminorm (see e.g. \cite[Theorem 7.7, p.35]{Mus83}),
\[
\|w\|_{L^\Phi} := \inf\left \{ \lambda>0 \st \int_{\hat{\Omega}} \int_{\hat{\Omega}} \Phi \left(x,y,\frac{w(x,y)}{\lambda} \right) \;  \ud\mu (y) \; \ud\mu (x)  \leq 1\right \}.
\]
 Every function in $L^\Phi (\hat{\Omega}\times\hat{\Omega} )$ vanishes almost everywhere on ${\mathbf C}_\Phi$, and the seminorm $\| w\|_{L^\Phi}$ equals $0$ if and only if $w=0$ almost everywhere in ${\mathbf S}_\Phi$. In other words, 
 \[
 L^\Phi ({\mathbf S}_\Phi) = \left \{ w \in L^0 ({\mathbf S}_\Phi) \st \exists \lambda >0 : \iint_{{\mathbf S}_\Phi} \Phi \left(x,y,\frac{w(x,y)}{\lambda}\right) \; \ud\mu (y) \; \ud\mu (x)<\infty \right\} 
 \] 
 is a Banach space (see e.g. Musielak \cite{Mus83} and Rao \& Ren \cite{RaRe91}).

We then define the {\em Musielak-Orlicz-Sobolev type spaces}
\begin{align} 
\nonumber W^{\Phi,0} (\hat{\Omega}) & := \left\{ \hat{u}\in L^0 (\hat{\Omega} ) \st  \exists \lambda >0 : \int_{\hat{\Omega}} \int_{\hat{\Omega}} \Phi (x,y,(\hat{u}(x) -\hat{u}(y))/\lambda ) \; \ud\mu (y) \; \ud\mu (x) <\infty \right\} , \text{ and} \\
\label{eq.space.w} W^{\Phi , 2} (\hat{\Omega} , \Omega ) & := \left\{ \hat{u}\in W^{\Phi,0} (\hat{\Omega}) \st \hat{u}|_\Omega \in L^2 (\Omega )  \right\} ,
\end{align}
and equip these spaces with the seminorms
\begin{align*}
\| \hat{u}\|_{W^{\Phi ,0}(\hat{\Omega})} & := \inf \left\{ \lambda >0 \st \int_{\hat{\Omega}} \int_{\hat{\Omega}} \Phi (x,y,(\hat{u}(x) -\hat{u}(y))/\lambda ) \; \ud\mu (y) \; \ud\mu (x) \leq 1\right \} , \text{ and} \\
\| \hat{u}\|_{W^{\Phi ,2}(\hat{\Omega},\Omega)} & := \| \hat{u}|_\Omega \|_{L^2 (\Omega )} + \| \hat{u}\|_{W^{\Phi ,0}(\hat{\Omega})}.
\end{align*}
We call $W^{\Phi,0}(\hat{\Omega})$ the {\em finite energy space}. These spaces appear naturally in the study of nonlocal operators, see also Felsinger, Kassmann \& Voigt \cite{FeKaVo15} and Servadei \& Valdinoci \cite{SeVa14} in the linear case. By definition of the spaces $W^{\Phi,0} (\hat{\Omega})$ and $W^{\Phi,2} (\hat{\Omega} ,\Omega )$, the {\em difference operator} 
\begin{equation} \label{eq.D}
\begin{split}
    D : W^{\Phi ,0} (\hat{\Omega} ) \; \to L^\Phi (\hat{\Omega}\times\hat{\Omega} ) , \qquad
    \hat{u}  \mapsto D\hat{u}:\;( (x,y) \mapsto \hat{u} (x) - \hat{u} (y)) ,
\end{split}
\end{equation} 
and the {\em restriction operator} 
\begin{equation} \label{eq.j}
j: W^{\Phi , 2} (\hat{\Omega} , \Omega )  \to L^2 (\Omega ) , \;
\hat{u}  \mapsto j(\hat{u}) :=\; \hat{u}_{|\Omega} ,
\end{equation}
are well-defined, linear and continuous. The operator $j$ is not necessarily injective; think of the case when the measure of $\hat{\Omega}\setminus\Omega$ is nonzero. 

\begin{lemma} \label{lem.banach}
Let $\Phi :\hat{\Omega} \times\hat{\Omega}\times\R \to [0,\infty]$ be a function satisfying the measurability and convexity condition \eqref{cond.phi}, and the $\Delta_2$-condition. Suppose that the set ${\mathbf S}_\Phi$ satisfies the {\em thickness condition}, that is,
\begin{equation} \label{cond.connected}
    \begin{cases}
        & \displaystyle\text{there exist a sequence } (A^{(0)}_n)_{n\in\N} \text{ of measurable subsets of } \Omega \\
        &\displaystyle \text{and a set } \hat{N}\subseteq\hat{\Omega} \text{ of measure zero such that, when we define and choose recursively for } j\geq 0 \\
        &\displaystyle I^{(j)} := \{ n\in\N \st \mu (A^{(j)}_n ) >0 \} , \\
        & \displaystyle(B^{(j+1)}_n)_{n\in\N} \text{ any sequence of measurable sets such that }  (A^{(j)}_n \times B^{(j+1)}_n) \subseteq {\mathbf S}_\Phi \text{ for every } n\in\N, \\
        & \displaystyle(A^{(j+1)}_n)_{n\in\N} \text{ any sequence of measurable sets such that } \bigcup_{n} A^{(j+1)}_n = \bigcup_{n} B^{(j+1)}_n , \\
        & \displaystyle\text{then } \hat{\Omega} \setminus \hat{N} = \bigcup_{j\in\N} \bigcup_{n\in I^{(j)}} A^{(j)}_n .
    \end{cases}
\end{equation}
Then the space $W^{\Phi ,2} (\hat{\Omega} , \Omega )$ is a Banach space. Moreover, every convergent sequence $(\hat{u}_k)$ in $W^{\Phi ,2} (\hat{\Omega} ,\Omega )$ has a subsequence which converges almost everywhere in $\hat{\Omega}$. 

If, in addition, the Fenchel conjugate $\Phi^*$ satisfies the $\Delta_2$-condition, then the space $W^{\Phi ,2} (\hat{\Omega} , \Omega )$ is reflexive.
\end{lemma}

\begin{proof}
As remarked above, the space $W^{\Phi ,2} (\hat{\Omega} , \Omega )$ is indeed a seminormed space. 

Let $(\hat{u}_k)$ be a Cauchy sequence in $W^{\Phi ,2} (\hat{\Omega} , \Omega )$. Then $(j(u_k)) = (\hat{u}_k|_{\Omega})$ is a Cauchy sequence in $L^2 (\Omega )$ and $(D\hat{u}_k)$ is a Cauchy sequence in $L^\Phi ({\mathbf S}_\Phi )$, where $j$ is the restriction operator and $D$ is the difference operator (see \eqref{eq.j} and \eqref{eq.D}). There exist a subsequence (which we denote again by $(\hat{u}_k)$) and sets $N^{(0)}\subseteq\Omega$ and $\hat{M}\subseteq\hat{\Omega}\times\hat{\Omega}$ of measure zero such that $\lim_{k\to\infty} \hat{u}_k(x)$ exists for every $x\in\Omega\setminus N^{(0)}$, and $\lim_{k\to\infty} (\hat{u}_k(x) - \hat{u}_k(y))$ exists for every $(x,y)\in {\mathbf S}_\Phi \setminus \hat{M}$ (norm convergence and modular convergence are equivalent when $\Phi$ satisfies the $\Delta_2$-condition, see Musielak \cite[5.2, p. 18]{Mus83}, and modular convergence implies, after subsequence, the convergence almost everywhere). As $A^{(0)}_n$ is a subset of $\Omega$, $\lim_{k\to\infty} \hat{u}_k(x)$ exists for every $x\in A^{(0)}_n \setminus N^{(0)}$ ($n\in\N$). Define $I^{(0)}$ as in \eqref{cond.connected}.

Assume that, for some $j\in\N$, there exists a set $N^{(j)}\subseteq\hat{\Omega}$ of measure zero such that, for every $n\in I^{(j)}$, $\lim_{k\to\infty} \hat{u}_k(x)$ exists for every $x\in A^{(j)}_n \setminus N^{(j)}$ (as this is for example the case when $j=0$). Let $N^{(j+1)}_n := \{ y\in B^{(j+1)}_n \st (A^{(j)}_n \setminus N^{(j)}) \times \{ y\} \subseteq \hat{M} \}$ ($n\in I^{(j)}$). Then $(A^{(j)}_n \setminus N^{(j)}) \times N^{(j+1)}_n \subseteq \hat{M}$. As $A^{(j)}_n$ has positive measure and as $N^{(j)}$ has measure zero, the difference $A^{(j)}_n \setminus N^{(j)}$ has positive measure. As $\hat{M}$ has measure zero, the set $N^{(j+1)}_n$ necessarily has measure zero. Moreover, as the difference $A^{(j)}_n \setminus N^{(j)}$ has positive measure, it is nonempty. By definition of the set $N^{(j+1)}_n$, for every $y\in B^{(j+1)}_n \setminus N^{(j+1)}_n$ there exists $x\in A^{(j)}_n\setminus N^{(j)}$ such that $(x,y)\not\in \hat{M}$. Hence, for every $y\in B^{(j+1)}_n \setminus N^{(j+1)}_n$ and for every such $x\in A^{(j)}_n\setminus N^{(j)}$, $\lim_{k\to\infty} \hat{u}_k (x)$ exists and $\lim_{k\to\infty} (\hat{u}_k (x) - \hat{u}_k(y))$ exists. This implies that $\lim_{k\to\infty} \hat{u}_k(y)$ exists for every $y\in B^{(j+1)}_n \setminus N^{(j+1)}_n$. When we let $N^{(j+1)} = \bigcup_{n\in I^{(j)}} N^{(j+1)}_n$, then $N^{(j+1)}\subseteq\hat{\Omega}$ is a set of measure zero and for every $y\in B^{(j+1)}_n \setminus N^{(j+1)}$, $\lim_{k\to\infty} \hat{u}_k(y)$ exists. If one chooses arbitrarily the measurable sets $A^{(j+1)}_n$ in such a way that $\bigcup_{n} A^{(j+1)}_n = \bigcup_{n} B^{(j+1)}_n$, then for every $y\in A^{(j+1)}_n \setminus N^{(j+1)}$ the limit $\lim_{k\to\infty} \hat{u}_k(y)$ exists. 

This induction argument shows that $\lim_{k\to\infty} \hat{u}_k(x)$ exists for every $x\in\bigcup_{j\in\N} \bigcup_{n\in I^{(j)}} (A^{(j)}_n \setminus N^{(j)})$, that is, by assumption on this union, for almost every $x\in\hat{\Omega}$. When we denote by $\hat{u}$ this pointwise limit, then Fatou's lemma implies
\begin{align*}
  \int_{\hat{\Omega}} \int_{\hat{\Omega}} \Phi (x,y,\hat{u}(x)-\hat{u}(y))\;\ud\mu (y) \; \ud\mu(x) & \leq \liminf_{k\to\infty} \int_{\hat{\Omega}} \int_{\hat{\Omega}} \Phi (x,y,\hat{u}_k(x)-\hat{u}_k(y))\;\ud\mu (y) \; \ud\mu(x) < \infty .
\end{align*}
Together with the convergence of $(\hat{u}_k|_\Omega)$ in $L^2 (\Omega )$ this implies $\hat{u}\in W^{\Phi,2} (\hat{\Omega} ,\Omega )$. But Fatou's lemma also implies $\lim_{k\to\infty} \hat{u}_k = \hat{u}$ in $W^{\Phi ,2} (\hat{\Omega} , \Omega )$. Hence, the space $W^{\Phi ,2} (\hat{\Omega} , \Omega )$ is complete.

When $\Phi$ and $\Phi^*$ satisfy the $\Delta_2$-condition,  then $L^\Phi ({\mathbf S}_\Phi )$ is a reflexive Banach space 
by \cite[Corollary 14, p.21]{RaRe91}, or Hudzik \cite[Theorem 1.4]{Hu84}. As $W^{\Phi , 2} (\hat{\Omega} , \Omega )$ can be seen as a closed subspace of $L^2(\Omega ) \times L^\Phi ({\mathbf S}_\Phi )$ in a natural way, the space $W^{\Phi , 2} (\hat{\Omega} , \Omega )$ is reflexive, too.
\end{proof}

Throughout the article, we say that $\Phi$ satisfies the {\em standard conditions} if 
\begin{equation} \label{cond.standard}
    \begin{cases}
         \Phi \text{ satisfies the measurability and convexity condition \eqref{cond.phi}}, \\
         \Phi \text{ and } \Phi^* \text{ satisfy the $\Delta_2$-condition, and} \\
         {\mathbf S}_\Phi \text{ satisfies the thickness condition \eqref{cond.connected}.}
    \end{cases}
\end{equation}
Given a function $\Phi$ satisfying the standard conditions, we consider the energy $\E : W^{\Phi , 2} (\hat{\Omega} , \Omega ) \to [0,\infty ]$ given by
\begin{equation} \label{energy.e}
\begin{cases}
 \E (\hat{u}) & := \displaystyle \frac12 \, \int_{\hat{\Omega}} \int_{\hat{\Omega}} \Phi (x,y,\hat{u}(x) - \hat{u}(y)) \; \ud\mu (y) \; \ud\mu (x) \\ 
  & = \displaystyle\frac12\, \iint_{{\mathbf S}_\Phi} \Phi (x,y,\hat{u}(x) - \hat{u}(y)) \; \ud\mu (y) \; \ud\mu (x) .
 \end{cases}
\end{equation} 
The integral over $\hat{\Omega}\times\hat{\Omega}$ is equal to the integral over ${\mathbf S}_\Phi$, because $\Phi (x,y,\cdot ) = 0$ for every $(x,y)\in{\mathbf N}_\Phi$. In other words, the difference $\hat{u} (x) - \hat{u} (y)$ does not contribute anything to the energy if $(x,y)\in {\mathbf N}_\Phi$. Using terminology from graph theory, one could say that the vortices $x$ and $y$ are not {\em connected by an edge} if $(x,y)\in {\mathbf N}_\Phi$, although this analogy only serves for an intuition: on the one hand $\hat{\Omega}$ is in general not a graph and on the other hand the two singletons $\{ x\}$ and $\{ y\}$ may have measure zero. However, if there exists a measurable set $B$ of positive measure such that $B\times\hat{\Omega} \subseteq {\mathbf N}_\Phi$, then the set $B$ is not connected to any point in $\hat{\Omega}$ by "an edge". The values of functions $\hat{u}\in W^{\Phi ,2} (\hat{\Omega} ,\Omega)$ on the set $B$ do not contribute anything to the energy $\E$. It is therefore reasonable to assume that there is no measurable set $B\subseteq\hat{\Omega}$ of positive measure such that $B\times\hat{\Omega} \subseteq {\mathbf N}_\Phi$. If there was such a set, then one may replace $\hat{\Omega}$ and $\Omega$ by the smaller sets $\hat{\Omega} \setminus B$ and $\Omega\setminus B$, respectively. 

\begin{theorem} \label{thm.e.lsc}
Assume that $\Phi$ satisfies the standard condition \eqref{cond.standard}. Then the energy $\mathcal E : W^{\Phi ,2} (\hat{\Omega} , \Omega )\to [0,\infty ]$ defined in \eqref{energy.e} is proper, continuous, convex and $j$-elliptic, where $j$ is the restriction operator from $W^{\Phi , 2} (\hat{\Omega} , \Omega )$ into $L^2 (\Omega )$. The $j$-subgradient $\partial_j\E$ is equal to the subgradient of a proper, lower semicontinuous and convex function $\E^{L^2} :L^2 (\Omega )\to [0,\infty ]$, and it is therefore $m$-accretive. For every initial value $u_0\in\overline{j(W^{\Phi ,2} (\hat{\Omega} ,\Omega ))}^{L^2(\Omega)}$ and every right hand side $f\in L^2 (0,T;L^2 (\Omega ))$ the abstract gradient system 
\begin{equation} \label{cp.j}
    \dot{u} + \partial_j\E (u) \ni f \text{ on } [0,\infty [ , \qquad u(0) = u_0 ,
\end{equation}
admits a unique strong solution $u\in C([0,\infty [ ; L^2 (\Omega )) \cap W^{1}_{2,loc} (]0,\infty [ ; L^2 (\Omega ))$.
In particular, the negative $j$-subgradient of $\E$ generates a strongly continuous semigroup of (nonlinear) contractions $S^N$ on $\overline{j(W^{\Phi ,2} (\hat{\Omega} ,\Omega ))}^{L^2 (\Omega )}$. For every $(u,f)\in\partial_j\E$ there exists an elliptic extension $\hat{u}\in W^{\Phi,2} (\hat{\Omega},\Omega )$ of $u$, namely 
\[
\hat{u} = \argmin_{j(\hat{w})=u} (\E (\hat{w}) - \int_\Omega f\hat{w}) = \argmin_{j(\hat{w})=u} \E (\hat{w}) ,
\]
and the elliptic extension is unique if $\Phi (x,y,\cdot )$ is strictly convex for every $(x,y)\in {\mathbf S}_\Phi$. 
\end{theorem}

\begin{proof}
The energy $\E$ is proper, lower semicontinuous and convex on $W^{\Phi , 2} (\hat{\Omega} , \Omega )$. Indeed, by the $\Delta_2$-condition and by the definition of the space $W^{\Phi ,2} (\hat{\Omega} , \Omega )$, the function $\E$ takes finite values everywhere on $W^{\Phi , 2} (\hat{\Omega} , \Omega )$, and so $\E$ is proper. By Fatou's lemma and by the convexity of $\Phi$ in the last variable (which implies continuity of the function $\Phi (x,y,\cdot )$ for every $(x,y)\in (\hat{\Omega}\times\hat{\Omega})\setminus {\mathbf C}_\Phi$), the function $\E$ is lower semicontinuous and convex (strictly convex, when $\Phi (x,y,\cdot )$ is strictly convex for every $(x,y)\in{\mathbf S}_\Phi$). As the effective domain $\dom{\E}$ is equal to the entire space $W^{\Phi , 2} (\hat{\Omega} , \Omega )$, and as $\E$ is lower semicontinuous, convex, the energy is thus continuous (see Borwein \& Vanderwerff \cite[Proposition 4.1.5]{BoVa10}). 

By reflexivity of the space $W^{\Phi , 2} (\hat{\Omega} , \Omega )$ and by continuity and convexity of the energy $\E$, the weak compactness of the sublevel sets $\{ \hat{u} \in W^{\Phi , 2} (\hat{\Omega} , \Omega ) \st \E (\hat{u}) + \| \hat{u}|_\Omega \|_{L^2 (\Omega )} \leq c\}$ follows simply from their boundedness. The boundedness of the sublevel sets is, however, clear from the definition of the space $W^{\Phi , 2} (\hat{\Omega} , \Omega )$. In other words, $\E$ is $j$-elliptic. The existence of a proper, lower semicontinuous and convex function $\E^{L^2} : L^2 (\Omega ) \to [0,\infty ]$ such that $\partial_j\E = \partial\E^{L^2}$, as well as the existence, possibly uniqueness, and characterization of the elliptic extensions follow from the general theory; see Section \ref{sec.preliminaries}. The same is true for the fact that the negative $j$-subgradient generates a semigroup of contractions and that the abstract gradient system is wellposed.
\end{proof}

Let us identify the $j$-subgradient of $\E$.
We define the {\em associate space} of $L^\Phi (\hat{\Omega}\times\hat{\Omega})$ by
\begin{align*}
    L^\Phi (\hat{\Omega}\times\hat{\Omega})_a & := \{ w\in L^0 (\hat{\Omega}\times\hat{\Omega} ) \st wz\in L^1 (\hat{\Omega}\times\hat{\Omega}) \text{ for every } z\in L^\Phi (\hat{\Omega}\times\hat{\Omega} ) \} .
\end{align*}

\begin{lemma} \label{lem.identification}
Assume that $\Phi$ satisfies the standard condition \eqref{cond.standard}. In addition, assume that $\Phi$ is continuously partially differentiable in the third variable in the sense that there exists a function $\varphi : \hat{\Omega} \times\hat{\Omega}\times\R\to\R$ such that 
\begin{equation}  \label{cond.phi.derivative}
\begin{cases}
\varphi (\cdot , \cdot , s) \text{ is measurable and symmetric on } \hat{\Omega}\times\hat{\Omega} \text{ for every } s\in\R, \\
 \varphi (x,y,0 ) = 0 \text{ for every } x, y\in \hat{\Omega } , \\
 \varphi (x,y,\cdot ) \text{ is continuous, increasing and odd for every } x, y\in\hat{\Omega} , \text{ and} \\
\displaystyle \Phi (x,y,s) = \int_0^s \varphi (x,y,\tau ) \; \ud\tau \text{ for every }  x, y\in\hat{\Omega} , s\in\R .
\end{cases}
\end{equation}
Then the energy $\hat{\E} : W^{\Phi,0} (\hat{\Omega}) \to \R$ given by
\[
\hat{\E} (\hat{u}) = \frac12 \int_{\hat{\Omega}} \int_{\hat{\Omega}} \Phi (x,y,\hat{u}(x) - \hat{u}(y))\; \ud\mu (y) \; \ud\mu (x) \quad (\hat{u}\in W^{\Phi,0} (\hat{\Omega} ))
\]
is G\^ateaux differentiable on $W^{\Phi ,0} (\hat{\Omega})$, for every $\hat{u}$, $\hat{v}\in W^{\Phi ,0} (\hat{\Omega} )$ the function
\begin{align*}
\hat{\Omega}\times\hat{\Omega} & \to \R , \quad (x,y) \mapsto \varphi (x,y,\hat{u} (x) - \hat{u} (y)) (\hat{v} (x) - \hat{v} (y))
\end{align*}
is integrable and 
\begin{equation} \label{eq.e.gateaux}
\hat{\E}'(\hat{u}) \hat{v} = \frac12 \int_{\hat{\Omega}} \int_{\hat{\Omega}} \varphi (x,y,\hat{u} (x) - \hat{u} (y)) (\hat{v} (x) - \hat{v} (y)) \; \ud\mu (y) \; \ud\mu (x) .
\end{equation}
In particular, the energy $\E$ from \eqref{energy.e}, which is the restriction of the energy $\hat{\E}$ to the space $W^{\Phi,2} (\hat{\Omega},\Omega )$, is G\^ateaux differentiable and the representation \eqref{eq.e.gateaux} holds with $\hat{\E}$ replaced by $\E$ and for every $\hat{u}$, $\hat{v}\in W^{\Phi,2} (\hat{\Omega},\Omega )$.
\end{lemma}

\begin{proof}
The fact that the energy $\hat{\E}$ is G\^ateaux differentiable on $W^{\Phi ,0} (\hat{\Omega})$ follows from the $\Delta_2$-condition and from the assumption that $\Phi$ is partially continuously differentiable in the third variable (see condition \eqref{cond.phi.derivative}). Let us give the argument.

By the $\Delta_2$-condition, there exist a constant $K_{\Delta_2}\geq 0$ and an integrable function $h_{\Delta_2}$ such that, for every $(x,y)\in\hat{\Omega}\times\hat{\Omega}$ and every $s\geq 0$,
\begin{align*}
    (K_{\Delta_2}-1) \Phi (x,y,s) + h_{\Delta_2} (x,y) & \geq \Phi (x,y,2s) - \Phi (x,y,s) \\
    & = \int_s^{2s} \varphi (x,y,\tau ) \; \ud\tau \\
    & \geq s \, \varphi (x,y,s) ;
\end{align*}
in the last inequality we used the fact that $\varphi$ is increasing in the last variable. This estimate between $\Phi (x,y,s)$ and $s\, \varphi (x,y,s)$ implies that for every $w\in L^\Phi (\hat{\Omega}\times\hat{\Omega})$ one has
\begin{align*}
\lefteqn{\int_{\hat{\Omega}} \int_{\hat{\Omega}} | \varphi (x,y,w(x,y)) | \, |w(x,y)|\; \ud \mu (y) \; \ud\mu (x) }\\
& = \int_{\hat{\Omega}} \int_{\hat{\Omega}}  \varphi (x,y,|w(x,y)|) \, |w(x,y)|\; \ud \mu (y) \; \ud\mu (x) \\
& \leq (K_{\Delta_2} -1) \int_{\hat{\Omega}} \int_{\hat{\Omega}} \Phi (x,y,w(x,y)) \; \ud \mu (y) \; \ud\mu (x) + \int_{\hat{\Omega}} \int_{\hat{\Omega}} h_{\Delta_2} (x,y) \; \ud \mu (y) \; \ud\mu (x) \\
& <\infty ,
\end{align*}
that is, for every $w\in L^\Phi (\hat{\Omega}\times\hat{\Omega})$, the function
\begin{align*}
\hat{\Omega}\times\hat{\Omega} & \to \R , \quad
(x,y) \mapsto \varphi (x,y,w(x,y)) w(x,y) 
\end{align*}
is integrable. As a consequence, for every $w$, $z\in L^\Phi (\hat{\Omega}\times\hat{\Omega} )$,
\begin{align}
    \nonumber \lefteqn{\int_{\hat{\Omega}} \int_{\hat{\Omega}}  | \varphi (x,y,w(x,y)) | \, |z(x,y)|\; \ud \mu (y) \; \ud\mu (x)} \\
    \nonumber & = \int_{\hat{\Omega}} \int_{\hat{\Omega}} \varphi (x,y,|w(x,y)|) \, |z(x,y)|\; \ud \mu (y) \; \ud\mu (x) \\
    \nonumber & \leq \int_{\hat{\Omega}} \int_{\hat{\Omega}} | \varphi (x,y,|w(x,y)| + |z(x,y)|) | \, (|w(x,y)| + |z(x,y)|)\; \ud \mu (y) \; \ud\mu (x) \\
    \label{eq.associate} & < \infty ,
\end{align}
that is, for every $w$, $z\in L^\Phi (\hat{\Omega}\times\hat{\Omega} )$, the function
\begin{align*}
\hat{\Omega}\times\hat{\Omega} & \to \R , \quad
(x,y) \mapsto \varphi (x,y,w(x,y)) z(x,y) 
\end{align*}
is integrable. In other words, for every $w\in L^\Phi (\hat{\Omega}\times\hat{\Omega} )$ the function $\varphi (\cdot , \cdot , w(\cdot ,\cdot ))$ belongs to the associate space $L^\Phi (\hat{\Omega}\times\hat{\Omega})_a$. 

Now, let $\hat{u}$, $\hat{v}\in W^{\Phi ,0} (\hat{\Omega})$ and $\lambda >0$. By the mean value theorem, for almost every $(x,y)\in \hat{\Omega}\times\hat{\Omega}$ there exists $\xi_\lambda (x,y)\in \R$ such that $|\xi_\lambda (x,y)| \leq \lambda \, |\hat{v} (x) - \hat{v} (y)|$ and 
\begin{align*}
    \Phi (x,y, \hat{u} (x) - \hat{u} (y) + \lambda (\hat{v} (x) - \hat{v}(y))) - \Phi (x,y,\hat{u} (x) - \hat{u}(y)) & = \lambda \, \varphi (x,y,\hat{u} (x) - \hat{u} (y) + \xi_\lambda (x,y)) (\hat{v} (x) - \hat{v}(y)) .
\end{align*}
Integrating this equality, dividing by $\lambda$, and taking the limit as $\lambda\to 0+$ yields the G\^ateaux differentiability of $\E$ and 
\begin{align*}
\hat{\E}' (\hat{u}) \hat{v} & = \lim_{\lambda\to 0+} \frac{\hat{\E} (\hat{u} + \lambda \hat{v}) - \hat{\E} (\hat{u})}{\lambda} \\
& = \lim_{\lambda\to 0+} \frac12 \int_{\hat{\Omega}} \int_{\hat{\Omega}} \varphi (x,y,\hat{u} (x) - \hat{u}(y) + \xi_\lambda (x,y)) (\hat{v} (x) - \hat{v}(y)) \; \ud\mu (y) \; \ud \mu (x) \\
& = \frac12 \int_{\hat{\Omega}} \int_{\hat{\Omega}} \varphi (x,y,\hat{u} (x) - \hat{u} (y)) (\hat{v} (x) - \hat{v} (y)) \; \ud\mu (y) \; \ud\mu (x) .
\end{align*}
Here, we have used the continuity of $\varphi$ in the last variable, and Lebesgue's dominated convergence theorem with the dominating function 
\[
(x,y) \mapsto \varphi (x,y, |\hat{u}(x)-\hat{u}(y)| + |\hat{v}(x) -\hat{v}(y)|) \, |\hat{v}(x) - \hat{v}(y)| ,
\]
which is integrable as noted above.
\end{proof}

Note carefully that Lemma \ref{lem.identification} does not imply that for every $\hat{u}$, $\hat{v}\in W^{\Phi,0} (\hat{\Omega})$ (or even $\hat{u}$, $\hat{v}\in W^{\Phi,2} (\hat{\Omega},\Omega )$) the function
\begin{align*}
\hat{\Omega}\times\hat{\Omega} & \to \R , \quad 
(x,y) \mapsto \varphi (x,y,\hat{u} (x) - \hat{u}(y)) \hat{v} (x)
\end{align*}
is integrable. If it was integrable, then by Fubini's theorem, for almost every $x\in\{|\hat{v}| >0\}$, the function
\begin{align*}
    \hat{\Omega} & \to \R , \quad
    y \mapsto \varphi (x,y,\hat{u}(x) - \hat{u}(y)) 
\end{align*}
was integrable, but it is known from standard examples that in general this function is not integrable. Nevertheless, one can give a meaning to its integral as a Cauchy principal value. Let us explain what we understand by this concept in the case of general measure spaces. \\

First, fix a subspace ${\mathcal D}\subseteq W^{\Phi ,0} (\hat{\Omega} )$ of so-called test functions. This subspace need not be closed. It should not be too large in order to give some flexibility in the definition of the Cauchy principal value. On the other hand, it should not be too small in order to separate the points of the {\em associate space}
\begin{align*}
{\mathcal D}_a & := \{ \hat{u}\in L^0 (\hat{\Omega}) \st \hat{u}\hat{v}\in L^1 (\hat{\Omega}) \text{ for every } \hat{v}\in {\mathcal D} \} .
\end{align*}
Here we say that ${\mathcal D}$ {\em separates the points of the associate space ${\mathcal D}_a$} if for every $\hat{f}\in {\mathcal D}_a$
\[
\int_{\hat{\Omega}} \hat{f} \hat{v} = 0 \text{ for every } \hat{v}\in {\mathcal D} \quad \Rightarrow \quad \hat{f} = 0 .
\]

Let $\hat{f}\in{\mathcal D}_a$ and ${\mathbf S}_\Phi (x):=\{y\in\hat{\Omega}:\; (x,y)\in {\mathbf S}_\Phi\}$.
We say that a function $\hat{u}\in W^{\Phi,0} (\hat{\Omega} )$ is a {\em strong solution} of 
\begin{equation} \label{eq.cpv}
 P.V.\int_{\hat{\Omega}} \varphi (x,y,\hat{u}(x)-\hat{u}(y))\;\ud\mu (y) = P.V.\int_{{\mathbf S}_\Phi (x)} \varphi (x,y,\hat{u}(x)-\hat{u}(y))\;\ud\mu (y) \stackrel{!}{=} \hat{f} (x)    \quad (x\in\hat{\Omega}) ,
\end{equation}
or that the {\em Cauchy principal value} of the ''integral'' $\int_{\hat{\Omega}} \varphi (x,y,\hat{u}(x) -\hat{u}(y)) \, \ud\mu (y)$ {\em exists in a strong sense} and equals $\hat{f}$, if
\begin{equation} \label{cond.cauchy.pv}
    \begin{cases}
    & \displaystyle\text{there exists an increasing sequence } ({\mathbf S}_n) \text{ of measurable and symmetric subsets of } \hat{\Omega}\times\hat{\Omega} \\
    & \displaystyle\text{ such that } \bigcup_n {\mathbf S}_n = {\mathbf S}_\Phi \text{ and, when one defines } {\mathbf S}_n(x) := \{ y\in\hat{\Omega} \st (x,y)\in {\mathbf S}_n\} ,  \\
    & \displaystyle \text{then for every } \hat{v}\in{\mathcal D} \\
    & \displaystyle\iint_{{\mathbf S}_n} | \varphi (x,y,\hat{u} (x) - \hat{u} (y)) \hat{v} (x) | \; \ud\mu (y) \; \ud\mu (x) <\infty \text{ for every } n\in\N , \text{ and} \\
    & \displaystyle\lim_{n\to\infty} \int_{\hat{\Omega}} \int_{{\mathbf S}_n (x)} \varphi (x,y,\hat{u}(x) - \hat{u}(y)) \; \ud\mu (y) \hat{v} (x) \; \ud\mu (x) = \int_{\hat{\Omega}} \hat{f} (x) \hat{v} (x) \; \ud\mu (x)  . \\
    \end{cases}
\end{equation}
We say that $\hat{u}\in W^{\Phi,0} (\hat{\Omega})$ is a {\em weak solution} of \eqref{eq.cpv}, or that the Cauchy principal value of the ''integral'' $\int_{\hat{\Omega}} \varphi (x,y,\hat{u}(x) -\hat{u}(y)) \, \ud\mu (y)$ {\em exists in a weak sense} and equals $\hat{f}$, if for every $\hat{v}\in{\mathcal D}$ 
\[
\frac12 \int_{\hat{\Omega}} \int_{\hat{\Omega}} \varphi (x,y,\hat{u}(x)-\hat{u} (y)) (\hat{v}(x) - \hat{v}(y)) \;\ud\mu(y) \;\ud\mu(x) = \int_{\hat{\Omega}} \hat{f} (x) \hat{v}(x) \; \ud\mu(x) .
\]
Note carefully that the above definitions of the Cauchy principal value in a strong or weak sense depend on the choice of the space ${\mathcal D}$. 

\begin{lemma}\label{lem-3.4}
    Every strong solution of \eqref{eq.cpv} is a weak solution of \eqref{eq.cpv}.
\end{lemma}

\begin{proof}
Let $\hat{u}\in W^{\Phi,0} (\hat{\Omega})$ be a strong solution of \eqref{eq.cpv}. Let $({\mathbf S}_n)$ be the sequence of measurable and symmetric subsets of $\hat{\Omega}\times\hat{\Omega}$ from the definition of strong solution / the Cauchy principal value (see \eqref{cond.cauchy.pv}). Then the function $(x,y)\mapsto \varphi (x,y,\hat{u}(x)-\hat{u}(y)) \hat{v}(x)$ is integrable on ${\mathbf S}_n$ for every $n\in\N$ and every $\hat{v}\in{\mathcal D}$ (see \eqref{cond.cauchy.pv}). By the integrability of the function $(x,y)\mapsto \varphi (x,y,\hat{u}(x)-\hat{u}(y)) (\hat{v}(x) -\hat{v}(y))$ over $\hat{\Omega}\times\hat{\Omega}$ (Lemma \ref{lem.identification}), for every $n\in\N$ and every $\hat{v}\in{\mathcal D}$ the function $(x,y)\mapsto \varphi (x,y,\hat{u}(x)-\hat{u}(y)) \hat{v}(y)$ is integrable on ${\mathbf S}_n$, too. As the set ${\mathbf S}_n$ is symmetric and as the function $\varphi$ is symmetric with respect to the first two variables and odd in the third variable,
\begin{align*}
    \lefteqn{- \iint_{{\mathbf S}_n} \varphi (x,y,\hat{u}(x) - \hat{u}(y)) \hat{v}(y) \; \ud\mu (y) \; \ud\mu (x)} \\
    & = - \iint_{{\mathbf S}_n} \varphi (y,x,\hat{u}(x) - \hat{u}(y)) \hat{v} (y) \; \ud\mu (y) \; \ud\mu (x) \\
    & = \iint_{{\mathbf S}_n} \varphi (y,x,\hat{u}(y) - \hat{u}(x)) \hat{v} (y) \; \ud\mu (y) \; \ud\mu (x) \\
    & = \iint_{{\mathbf S}_n} \varphi (x,y,\hat{u}(y) - \hat{u}(x)) \hat{v} (x) \; \ud\mu (y) \; \ud\mu (x) .
\end{align*}
This implies that, for every $n\in\N$ and every $\hat{v}\in{\mathcal D}$,
\[
\frac12 \iint_{{\mathbf S}_n} \varphi (x,y,\hat{u} (x) - \hat{u} (y)) (\hat{v} (x) - \hat{v} (y)) \; \ud\mu (y) \; \ud\mu (x) = \iint_{{\mathbf S}_n} \varphi (x,y,\hat{u}(y) - \hat{u}(x)) \; \ud\mu (y) \hat{v} (x) \; \ud\mu (x) .
\]
Note that $\bigcup_n {\mathbf S}_n = {\mathbf S}$, and $\varphi (x,y,\cdot ) = 0$ for $(x,y)\not\in {\mathbf S}$. By Lebesgue's dominated convergence theorem, by the preceding equality, and by the definition of Cauchy principal value / strong solution (see \eqref{cond.cauchy.pv}), for every $\hat{v}\in{\mathcal D}$, 
\begin{align*} 
\frac12 \int_{\hat{\Omega}} \int_{\hat{\Omega}} \varphi (x,y,\hat{u} (x) - \hat{u} (y)) (\hat{v} (x) - \hat{v} (y)) \; \ud\mu (y) \; \ud\mu (x) & = \lim_{n\to\infty} \frac12 \iint_{{\mathbf S}_n} \varphi (x,y,\hat{u} (x) - \hat{u} (y)) (\hat{v} (x) - \hat{v} (y)) \; \ud\mu (y) \; \ud\mu (x) \\
& = \lim_{n\to\infty} \iint_{{\mathbf S}_n} \varphi (x,y,\hat{u}(y) - \hat{u}(x)) \; \ud\mu (y) \hat{v} (x) \; \ud\mu (x) \\
& = \int_{\hat{\Omega}} \hat{f}(x) \hat{v} (x) \; \ud\mu (x) .
\end{align*}
Hence, $\hat{u}$ is a weak solution of \eqref{eq.cpv}.
\end{proof}

\begin{theorem} \label{thm.identification.2}
Assume that $\Phi$ satisfies the standard condition \eqref{cond.standard}, and that $\Phi$ is continuously partially differentiable in the third variable (condition \eqref{cond.phi.derivative}). Let $\E$ be the energy defined in \eqref{energy.e}, and let $\partial_j\E$ be its $j$-subgradient on $L^2 (\Omega )$, and let $u$, $f\in L^2 (\Omega )$. Assume that ${\mathcal D}\subseteq W^{\Phi ,2} (\hat{\Omega},\Omega )$ and that ${\mathcal D}$ separates the points of the associate space ${\mathcal D}_a$. 

Then $(u,f)\in \partial_j\E$ if and only if there exists $\hat{u} \in W^{\Phi ,2} (\hat{\Omega},\Omega )$ such that $\hat{u}|_{\Omega} = u$ and $\hat{u}$ is a weak solution of 
\begin{equation} \label{eq.nonlocal.weak.identification}
\begin{cases}
     \displaystyle P.V.\int_{\hat{\Omega}} \varphi (x,y,\hat{u} (x) - \hat{u} (y)) \; \ud\mu (y) = f(x) \quad &\text{ in } \Omega , \\
     \displaystyle P.V.\int_{\hat{\Omega}} \varphi (x,y,\hat{u} (x) - \hat{u} (y)) \; \ud\mu (y) = 0 \quad &\text{ in } \hat{\Omega}\setminus\Omega , 
\end{cases}
\end{equation}
and, moreover,  for every $\hat{v}\in W^{\Phi ,2} (\hat{\Omega},\Omega)$,
\begin{equation} \label{eq.nonlocal.weak.identification.bc}
\begin{split}
& \displaystyle \frac12 \int_{\hat{\Omega}} \int_{\hat{\Omega}} \varphi (x,y,\hat{u}(x)-\hat{u}(y)) (\hat{v} (x) -\hat{v} (y)) \; \ud\mu(y) \; \ud\mu(x)  \\
    & \phantom{\ud x} \displaystyle = \int_{\hat{\Omega}} P.V.\int_{\hat{\Omega}} \varphi (x,y,\hat{u} (x) - \hat{u} (y)) \; \ud\mu (y) \hat{v} (x) \; \ud\mu (x) .
\end{split}
\end{equation}
Every such $\hat{u}$ is an elliptic extension of $u$ and satisfies
\[
\hat{u} = \argmin_{j(\hat{w}) = u} ( \E (\hat{w}) - \int_\Omega f \hat{w} ) = \argmin_{j(\hat{w}) = u} \E (\hat{w}).
\]
The elliptic extension is unique when $\Phi (x,y,\cdot )$ is strictly convex for every $(x,y)\in{\mathbf S}_\Phi$.
\end{theorem}

\begin{proof}
Assume that $(u,f)\in\partial_j \E$. By definition of the $j$-subgradient, there exists an elliptic extension $\hat{u}\in W^{\Phi ,2} (\hat{\Omega},\Omega )$ such that $\hat{u}|_\Omega = u$ and, for every $\hat{v}\in W^{\Phi ,2} (\hat{\Omega},\Omega )$, $\E (\hat{u} + \hat{v} ) - \E (\hat{u}) \geq \int_\Omega f (x) \hat{v} (x) \; \ud\mu (x)$. In this inequality, one may replace $\hat{v}$ by $\lambda\hat{v}$ for $\lambda >0$, divide the inequality by $\lambda$ and then let $\lambda$ tend to $0$. From the G\^ateaux differentiability of $\E$ (Lemma \ref{lem.identification}) 
follows that for every $\hat{v}\in W^{\Phi ,2} (\hat{\Omega},\Omega )$
\begin{equation} \label{eq.gateaux}
\E' (\hat{u}) \hat{v}  
\geq \int_\Omega f (x) \hat{v} (x) \; \ud\mu (x) .
\end{equation}
Replacing $\hat{v}$ by $-\hat{v}$ in this inequality, one actually sees that the inequality sign can be replaced by an equality sign. Together with the representation of $\E'(\hat{u})\hat{v}$ from Lemma \ref{lem.identification}, this yields that $\hat{u}$ is a weak solution of \eqref{eq.nonlocal.weak.identification}, but as the preceding inequality not only holds for all $\hat{v}\in{\mathcal D}$, one actually obtains \eqref{eq.nonlocal.weak.identification.bc}.

Conversely, assume that $u$, $f\in L^2 (\Omega )$ are such that there exists $\hat{u} \in W^{\Phi ,2} (\hat{\Omega},\Omega )$ such that $\hat{u}|_{\Omega} = u$, $\hat{u}$ is a weak solution of \eqref{eq.nonlocal.weak.identification}, and \eqref{eq.nonlocal.weak.identification.bc} holds. Then the convexity of $\E$ (which implies that the function $\lambda \mapsto \frac{\E (\hat{u} + \lambda \hat{v}) - \E (\hat{u})}{\lambda}$ is increasing on $]0,\infty [$), Lemma \ref{lem.identification} and the equalities \eqref{eq.nonlocal.weak.identification.bc} and \eqref{eq.nonlocal.weak.identification} imply that, for every $\hat{v}\in W^{\Phi,2} (\hat{\Omega},\Omega )$,
\begin{align*}
\E (\hat{u} + \hat{v}) - \E (\hat{u}) & \geq \lim_{\lambda\to 0+} \frac{\E (\hat{u} + \lambda \hat{v}) - \E (\hat{u})}{\lambda} \\
& = \E'(\hat{u})\hat{v} \\
& = \frac12 \int_{\hat{\Omega}} \int_{\hat{\Omega}} \varphi (x,y,\hat{u}(x)-\hat{u}(y)) (\hat{v} (x) -\hat{v} (y)) \; \ud\mu(y) \; \ud\mu(x) \\
& = \int_{\hat{\Omega}} \left( P.V.\int_{\hat{\Omega}} \varphi (x,y,\hat{u} (x) - \hat{u} (y)) \; \ud\mu (y) \right) \hat{v} (x) \; \ud\mu (x) \\
& = \int_\Omega f(x) \hat{v}(x) \; \ud\mu (x) .
\end{align*}
This implies $(u,f)\in\partial_j\E$. 
\end{proof}

The two equations in \eqref{eq.nonlocal.weak.identification} are to be understood as one Cauchy principal value on $\hat{\Omega}$. The data $f\in L^2 (\Omega )$ appears in the first equation. The second equation in \eqref{eq.nonlocal.weak.identification} represents a complementary condition on $\hat{\Omega}\setminus\Omega$. We call this complementary condition {\em Neumann exterior condition} for the elliptic extension $\hat{u}$ of $u$.  

If we replaced the Cauchy principal value on the right hand side of equation \eqref{eq.nonlocal.weak.identification.bc} by the function $f$ (extended by $0$ outside $\Omega$), then equation \eqref{eq.nonlocal.weak.identification.bc} would be stronger than equation \eqref{eq.nonlocal.weak.identification}, simply because the space ${\mathcal D}$ used in the definition of the Cauchy principal value is a subspace of $W^{\Phi,2} (\hat{\Omega},\Omega )$. However, once the Cauchy principal value is determined, then 
equation \eqref{eq.nonlocal.weak.identification.bc} corresponds to a ''boundary condition'' in many examples. The ''boundary'' depends on the choice of $\Omega$, $\hat{\Omega}$, $S_\Phi \subseteq\hat{\Omega}\times\hat{\Omega}$, and the space ${\mathcal D}$. The term ''boundary condition'' becomes for example clear when there is an integration by parts formula, that is, equation \eqref{eq.nonlocal.weak.identification.bc} with an additional integral on the ''boundary''. More concretely, by the weak definition of the Cauchy principal value, the equality in \eqref{eq.nonlocal.weak.identification.bc} automatically holds for every $\hat{v}\in{\mathcal D}$. Assume that the equality also holds on the closure of the space ${\mathcal D}$ in $W^{\Phi,2} (\hat{\Omega},\Omega )$, which one may denote by $\mathring{W}^{\Phi,2} (\hat{\Omega},\Omega )$. By Lemma \ref{lem.identification}, the left hand side of \eqref{eq.nonlocal.weak.identification.bc} is a representation of $\E'(\hat{u})\hat{v}$, so the difference $\E' (\hat{u}) - P.V.\int_{\hat{\Omega}} \varphi (\cdot ,y,\hat{u} (\cdot ) - \hat{u} (y)) \; \ud\mu (y)$ vanishes on $\mathring{W}^{\Phi,2} (\hat{\Omega},\Omega )$. Hence, the difference can be seen as a bounded, linear functional on the quotient space $W^{\Phi ,2} (\hat{\Omega},\Omega ) / \mathring{W}^{\Phi ,2} (\hat{\Omega},\Omega )$. Typically this corresponds to a ''boundary condition''.

\begin{remark}
It is an open problem whether for every $(u,f)\in\partial_j\E$ there exists an elliptic extension $\hat{u}\in W^{\Phi,2} (\hat{\Omega},\Omega )$ for which the Cauchy principal value exists in a strong sense, that is, $\hat{u}$ is a strong solution of \eqref{eq.nonlocal.weak.identification}. This is even not clear in the situation when $\Omega$ and $\hat{\Omega}$ are open subsets of $\R^N$, even not in the case of quadratic energies. This open problem has been raised by Claus \& Warma \cite{ClWa20}. 

When $\Omega$ and $\hat{\Omega}$ are locally finite graphs, however, this problem has a simple positive answer because in that case the integral $\int_{{\mathbf S} (x)} \varphi (x,y,\hat{u}(x)-\hat{u}(y))\;\ud\mu(y)$ actually is a finite sum (see Example \ref{ex.3.2} below).
\end{remark}

\begin{corollary}[\bf Evolution equation with Neumann exterior conditions and Neumann boundary conditions] \label{cor.neumann}
Take the assumptions of Theorem \ref{thm.identification.2}. For every initial value $u_0\in \overline{j(W^{\Phi ,2} (\hat{\Omega},\Omega ))}^{L^2(\Omega)}$ and every right hand side $f\in L^2_{loc} ([0,\infty [;L^2 (\Omega ))$ the problem
\begin{equation} \label{eq.nonlocal.corollary}
 \begin{cases}
 \displaystyle \partial_t \hat{u} (t,x) + P.V.\int_{\hat{\Omega}} \varphi (x,y,\hat{u}(t,x)-\hat{u}(t,y)) \; \ud\mu(y) = f(t,x) &  \text{in } ]0,\infty [\times\Omega , \\
\displaystyle  P.V.\int_{\hat{\Omega}} \varphi (x,y,\hat{u}(t,x)-\hat{u}(t,y)) \; \ud\mu(y) = 0 & \text{in } ]0,\infty [ \times (\hat{\Omega} \setminus \Omega ) , \\
\displaystyle \frac12 \int_{\hat{\Omega}} \int_{\hat{\Omega}} \varphi (x,y,\hat{u} (t,x) - \hat{u} (t,y)) (\hat{v} (x) - \hat{v} (y)) \; \ud\mu (y) \; \ud\mu (x) = & \\ 
\phantom{\ud x} \displaystyle \int_{\hat{\Omega}} P.V.\int_{\hat{\Omega}} \varphi (x,y,\hat{u}(t,x)-\hat{u}(t,y)) \; \ud\mu(y) \hat{v} (x) \; \ud\mu (x) &  \text{in } ]0,\infty [ ,\;\forall \hat{v}\in W^{\Phi,2} (\hat{\Omega},\Omega ), \\
 \displaystyle \hat{u}(0,x) = u_0 (x) & \mbox{in }\Omega ,
 \end{cases}
\end{equation}
admits a unique solution $u\in C([0,\infty [;L^2 (\Omega )) \cap W^{1}_{2,loc} (]0,\infty [ ; L^2 (\Omega ))$ in the sense that for almost every $t\geq 0$ the function $u(t)\in L^2 (\Omega )$ admits an elliptic extension $\hat{u} (t)\in W^{\Phi ,2} (\hat{\Omega} ,\Omega )$, that is $\hat{u} (t)|_\Omega = u(t)$, and the first two equalities of \eqref{eq.nonlocal.corollary} hold in a weak sense. 
\end{corollary}

\begin{proof}
 Let $\E$ be the energy defined in \eqref{energy.e}. By Theorem \ref{thm.e.lsc}, this energy is proper, continuous, convex and $j$-elliptic, where $j$ is the restriction operator defined in \eqref{eq.j}. Moreover, the abstract gradient system \eqref{cp.j} with the $j$-subgradient of $\E$ is wellposed in $L^2(\Omega )$ in the sense that for every initial value $u_0\in\overline{j(W^{\Phi ,2} (\hat{\Omega},\Omega))}^{L^2(\Omega)}$ and every right hand side $f\in L^2_{loc} ([0,\infty [;L^2 (\Omega ))$ it admits a unique solution $u$ satisfying $u\in C([0,\infty [ ; L^2 (\Omega ))\cap W^{1}_{2,loc} (]0,\infty [;L^2 (\Omega ))$. For almost every $t> 0$, $(u(t),f(t) - \dot{u} (t))\in\partial_j\E$. The rest of the statement follows from the identification of the $j$-subgradient of $\E$ in Theorem \ref{thm.identification.2}.
\end{proof}

\begin{example}[\bf Nonlocal operators on open subsets in $\R^N$] \label{ex.3.1}
Let $\Omega\subseteq\hat{\Omega}\subseteq\R^N$ be two open, nonempty sets, equipped with the Lebesgue measure. 
Let $p\in ]1,\infty [$, and let $k:\hat{\Omega}\times\hat{\Omega}\to [0,\infty]$ be a measurable and symmetric kernel such that, for some constants $\theta\in ]0,1[$, $C_1$, $C_2 >0$ and for all $x$, $y\in\hat{\Omega}$,
\begin{equation} \label{cond.k}
C_1 \, \frac{1}{|x-y|^{N+\theta p}} \leq k(x,y) \leq C_2 \, \frac{1}{|x-y|^{N+\theta p}} .
\end{equation}
Finally, let ${\mathbf S} \subseteq \hat{\Omega}\times\hat{\Omega}$ be a symmetric, open set, define ${\mathbf S} (x) := \{ y\in\hat{\Omega} \st (x,y) \in {\mathbf S} \}$ for every $x\in\hat{\Omega}$, and 
\[
\partial_{\mathbf S}\hat{\Omega} := \{ x\in\overline{\hat{\Omega}} \st x\not\in {\mathbf S} (x)\} .
\]
The latter set is called the {\em ${\mathbf S}$-boundary} in the closure of $\hat{\Omega}$. We assume that ${\mathbf S}$ satisfies the thickness condition \eqref{cond.connected}. This condition is for example satisfied when (a) ${\mathbf S} = \hat{\Omega}\times\hat{\Omega}$ or (b) ${\mathbf S}=(\hat{\Omega}\times\Omega)\cup (\Omega\times\hat{\Omega})$. Indeed, there is nothing to show in the first case. In the second case, choose $A^{(0)}_0 = \Omega$ and the other $A^{(0)}_n$ ($n\geq 1$) arbitrary measurable subsets of $\Omega$. As $k(x,y) >0$ for every $x$, $y\in\hat{\Omega}$, one can choose $B^{(1)}_0 = A^{(1)}_0 = \hat{\Omega}$, and therefore $\bigcup_{j\in\N} \bigcup_{n\in I^{(j)}} A^{(j)}_n = \hat{\Omega}$. Clearly, one may think of other examples of sets ${\mathbf S}$ which satisfy the thickness condition.\\

For given $f\in L^2 (\Omega )$, consider the stationary problem of finding $u\in L^2 (\Omega )$ such that
\begin{equation} \label{eq.nonlocal.weak.example}
\begin{cases}
    \multicolumn{2}{l}{\text{there exists an extension of } u \text{ to a function } \hat{u}\in W^{\theta,{\mathbf S}}_{p,2} (\hat{\Omega},\Omega ), \text{ that is, } \hat{u}|_{\Omega} = u , \text{ such that}} \\
    \displaystyle P.V.\int_{{\mathbf S} (x)} k(x,y) |\hat{u} (x) - \hat{u} (y)|^{p-2} (\hat{u} (x) - \hat{u} (y)) \; \ud\mu (y) = f(x) & \text{for } x\in\Omega , \\
   \displaystyle P.V.\int_{{\mathbf S} (x)} k(x,y) |\hat{u} (x) - \hat{u} (y)|^{p-2} (\hat{u} (x) - \hat{u} (y)) \; \ud\mu (y) = 0 & \text{for } x\in\hat{\Omega}\setminus\Omega , \\
   \displaystyle \mathcal N_{\partial_{\mathbf S}\hat{\Omega}}^{2-2\theta,p}\hat{u} (x) = 0 & \text{for } x\in\partial_{\mathbf S}\hat{\Omega}. 
\end{cases}
\end{equation}
Consider also for an initial value $u_0 \in L^2 (\Omega )$ and a given right hand side $f\in L^2_{loc} ([0,\infty [ ; L^2 (\Omega ))$ the associated evolution problem
\begin{equation} \label{eq.nonlocal.weak.example.cp}
\begin{cases}
   \displaystyle  \partial_t \hat{u} (t,x) + P.V.\int_{{\mathbf S} (x)} k(x,y) |\hat{u} (t,x) - \hat{u} (t,y)|^{p-2} (\hat{u} (t,x) - \hat{u} (t,y)) \; \ud\mu (y) = f(t,x) \quad &\text{for } (t,x)\in ]0,\infty [ \times \Omega , \\
   \displaystyle P.V.\int_{{\mathbf S} (x)} k(x,y) |\hat{u} (t,x) - \hat{u} (t,y)|^{p-2} (\hat{u} (t,x) - \hat{u} (t,y)) \; \ud\mu (y) = 0 \quad &\text{for } (t,x)\in ]0,\infty [ \times (\hat{\Omega}\setminus\Omega) , \\
     \displaystyle \mathcal N_{\partial_{\mathbf S}\hat{\Omega}}^{2-2\theta,p}\hat{u} (t,x) = 0 & \text{for } (t,x)\in ]0,\infty [\times \partial_{\mathbf S}\hat{\Omega} , \\
   \displaystyle \hat{u} (0,x) = u_0 (x)&\text{for } x\in\Omega .
\end{cases}
\end{equation}
The boundary operator $\mathcal N_{\partial_{\mathbf S}\hat{\Omega}}^{2-2\theta,p}$ in the pevious systems will be clarified in \eqref{eq.operator.n} below.

Here,
\[
W^{\theta,{\mathbf S}}_{p,2} (\hat{\Omega},\Omega ) := \left\{ \hat{u}\in L^0 (\hat{\Omega} ) \st \hat{u}|_\Omega \in L^2 (\Omega ) \text{ and } \iint_{{\mathbf S}}  \frac{| \hat{u} (x) - \hat{u} (y) |^p}{|x-y|^{N+\theta p}}  \; \ud y \; \ud x <\infty \right\} ,
\]
is a fractional Sobolev type space equipped with the norm
\[
\| \hat{u} \|_{W^{\theta,{\mathbf S}}_{p,2} (\hat{\Omega},\Omega )} := \| \hat{u} \|_{L^2 (\Omega )} + \left( \iint_{{\mathbf S}}  \frac{| \hat{u} (x) - \hat{u} (y) |^p}{|x-y|^{N+\theta p}} \; \ud y \; \ud x \right)^{\frac{1}{p}} .
\]

Both problems fit into the setting of this article, when one defines the kernel $\Phi : \hat{\Omega}\times\hat{\Omega}\times\R \to [0,\infty [$ by
\[
\Phi (x,y,s) = \frac{1}{p}\, k(x,y) \, |s|^p \, 1_{{\mathbf S}} (x,y) ,
\] 
where $1_{\mathbf S}$ is the characteristic function of the set ${\mathbf S}$. This kernel satisfies the standard condition \eqref{cond.standard}, and it is continuously partially differentiable in the third variable. The space $W^{\Phi ,2} (\hat{\Omega},\Omega )$ is equal to the fractional Sobolev type space $W^{\theta,{\mathbf S}}_{p,2} (\hat{\Omega},\Omega )$, and it is a Banach space by Lemma \ref{lem.banach}. \\

The associated energy $\E : W^{\theta,{\mathbf S}}_{p,2} (\hat{\Omega},\Omega ) \to [0,\infty ]$ given by
\begin{align*}
\E (\hat{u} ) & = \frac{1}{2p}\, \iint_{{\mathbf S}} k(x,y) \, | \hat{u} (x) - \hat{u} (y) |^p \; \ud y \; \ud x \\
 & = \frac{1}{2p}\, \int_{\hat{\Omega}} \int_{{\mathbf S}(x)} k(x,y) \, | \hat{u} (x) - \hat{u} (y) |^p \; \ud y \; \ud x 
\end{align*} 
is G\^ateaux differentiable and, for every $\hat{u}$, $\hat{v}\in W^{\theta,{\mathbf S}}_{p,2} (\hat{\Omega},\Omega )$,
 \[
 \E'(\hat{u} ) \hat{v} = \frac12 \, \iint_{{\mathbf S}} k(x,y) |\hat{u} (x) - \hat{u} (y)|^{p-2} (\hat{u} (x) - \hat{u} (y)) (\hat{v} (x) - \hat{v} (y)) \; \ud y \; \ud x ,
 \]
the function under the integral being integrable (Lemma \ref{lem.identification}). \\

The Cauchy principal value is to be understood in a weak sense and with respect to the space of test functions
\[
{\mathcal D} := \{ w|_{\hat{\Omega}} \st w\in C^\infty_c (\R^N ) \text{ and } {\rm supp}(w) \cap \partial_{\mathbf S}\hat{\Omega}=\emptyset \} .
\]
It depends on the ${\mathbf S}$-boundary, as does the associate space 
\[
{\mathcal D}_a = L^1_{loc} (\overline{\hat{\Omega}} \setminus \partial_{\mathbf S}\hat{\Omega} ) .
\]
The sets
\[
{\mathbf S}_\varepsilon := \{ (x,y)\in{\mathbf S} \st |x-y|>\varepsilon \} \qquad (\varepsilon >0) 
\]
are chosen independently of the functions $\hat{u}$ (instead of the real parameter $\varepsilon >0$ one might also choose a null sequence $(\varepsilon_n)$). As for every $\varepsilon >0$ the kernel $k$ is bounded and integrable on ${\mathbf S}_\varepsilon$, and in particular $k\in L^p ( {\mathbf S}_\varepsilon)$, as $|D\hat{u}|^{p-1} \in L^{p'} (\hat{\Omega}\times\hat{\Omega})$ for every $\hat{u}\in W^{\theta,{\mathbf S}}_{p,2} (\hat{\Omega},\Omega )$, and as every $\hat{v}\in{\mathcal D}$ is bounded, then for every $\hat{u}\in W^{\theta,{\mathbf S}}_{p,2} (\hat{\Omega},\Omega )$ and every $\hat{v}\in{\mathcal D}$
\[
\iint_{{\mathbf S}_\varepsilon} k(x,y) |\hat{u} (x) - \hat{u} (y)|^{p-1} |\hat{v} (x)| \; \ud y \; \ud x <\infty .
\]
Hence, for given $\hat{f}\in L^1_{loc} (\overline{\hat{\Omega}} \setminus \partial_{\mathbf S}\hat{\Omega} )$ the function $\hat{u}\in W^{\theta,{\mathbf S}}_{p,2} (\hat{\Omega},\Omega )$ is a strong solution of
\[
P.V.\int_{{\mathbf S}(x)} k(x,y) |\hat{u} (x) - \hat{u} (y)|^{p-2} (\hat{u} (x) - \hat{u} (y)) \; \ud y = \hat{f} (x) \qquad (x\in\hat{\Omega})
\]
if, for some null sequence $(\varepsilon_n)$ of positive real numbers, for every $\hat{v}\in{\mathcal D}$,
\[
\lim_{n\to\infty} \int_{\hat{\Omega}} \int_{{\mathbf S}_{\varepsilon_n} (x)} k(x,y) |\hat{u} (x) - \hat{u} (y)|^{p-2} (\hat{u} (x) - \hat{u} (y)) \; \ud y \hat{v} (x) \; \ud x = \int_{\hat{\Omega}} \hat{f} (x) \hat{v} (x) \; \ud x .
\]
This is for example the case when 
\[
\lim_{\varepsilon\to 0+} \int_{{\mathbf S}_\varepsilon (x)} k(x,y) |\hat{u} (x) - \hat{u} (y)|^{p-2} (\hat{u} (x) - \hat{u} (y)) \; \ud y = \hat{f} (x) \text{ in } L^1_{loc} (\overline{\hat{\Omega}} \setminus \partial_{\mathbf S}\hat{\Omega} ) .
\]

\begin{lemma}
 Let $p\in ]1,\infty[$.
 If $0<\theta < \frac{p-1}{p}$ and $\hat{u}\in C^{0,\alpha} (\hat{\Omega}) \cap L^\infty (\hat{\Omega})$ for some $\frac{\theta p}{p-1} <\alpha\leq 1$, or if $\frac{p-1}{p}<\theta<1$ and $\hat{u}\in C^{1,\alpha} (\hat{\Omega})\cap L^\infty (\hat{\Omega})$ for some $\theta p - (p-1) <\alpha\leq 1$, then
 \[
 \lim_{\varepsilon\to 0+} \int_{{\mathbf S}_\varepsilon (x)} k(x,y) |\hat{u}(x) - \hat{u}(y)|^{p-2} (\hat{u}(x)-\hat{u}(y)) \;\ud y
 \]
 exists uniformly for $x$ from compact subsets of $\overline{\hat{\Omega}} \setminus \partial_{\mathbf S}\hat{\Omega}$.
\end{lemma}

\begin{proof}
{\bf First case: $0<\theta < \frac{p-1}{p}$ and $\hat{u}\in C^{0,\alpha} (\hat{\Omega}) \cap L^\infty (\hat{\Omega})$ for some $\frac{\theta p}{p-1} <\alpha\leq 1$.} In this case, there exists a constant $C\geq 0$ such that $|\hat{u}(x)-\hat{u}(y)|\leq C\,|x-y|^\alpha$ for every $x$, $y\in\hat{\Omega}$ and $\| \hat{u}\|_{L^\infty(\hat{\Omega})}\leq C$. As a consequence,
\[
|\hat{u}(x)-\hat{u}(y)|\leq C\,\min\{ |x-y|^\alpha , 2\} .
\]
Hence,
\begin{align*}
 \int_{{\mathbf S}(x)} k(x,y) |\hat{u}(x)-\hat{u}(y)|^{p-1} \;\ud y & \leq C_2 C^{p-1} \, \int_{\R^N} \frac{\min \{ | x-y|^{\alpha (p-1)}, 2^{p-1}\}}{| x-y|^{N+\theta p}} \; \ud y <\tilde{C} ,
\end{align*}
for some finite constant $\tilde{C}$ that is independent of $x\in \overline{\hat{\Omega}}$. The claim follows from Lebesgue's dominated convergence.

{\bf Second case: $\frac{p-1}{p}<\theta<1$ and $\hat{u}\in C^{1,\alpha} (\hat{\Omega}) \cap L^\infty (\hat{\Omega})$ for some $\theta p - (p-1) <\alpha\leq 1$.} We distinguish two subcases. 

First, if $x\not\in \overline{{\mathbf S} (x)}$, then $x$ has positive distance from ${\mathbf S} (x)$, and therefore the kernel $k$ is integrable in ${\mathbf S} (x)$. It follows that the function $y\mapsto k(x,y) |\hat{u}(x)-\hat{u}(y)|^{p-1}$ is integrable in ${\mathbf S} (x)$ by the boundedness of $\hat{u}$. The corresponding estimates hold uniformly in a neighbourhood of $x$, and the local uniform convergence follows.

Second, if $x\in {\mathbf S} (x)$, then, as ${\mathbf S}$ is open, there exists $r>0$ such that $B(x,r)\subseteq {\mathbf S}(x)$, and the radius $r$ can even be chosen uniformly in a neighbourhood of $x$. For every $\varepsilon \in ]0,r]$, 
\begin{align*}
  \int_{y\in{\mathbf S}(x) \atop |x-y| > \varepsilon} k(x,y) |\hat{u}(x) - \hat{u}(y)|^{p-2} (\hat{u}(x)-\hat{u}(y)) \;\ud y & = \int_{y\in{\mathbf S}(x) \atop |x-y| > r} k(x,y) |\hat{u}(x) - \hat{u}(y)|^{p-2} (\hat{u}(x)-\hat{u}(y)) \;\ud y \\
  & \phantom{= \ } + \int_{r> |x-y| > \varepsilon} k(x,y) |\hat{u}(x) - \hat{u}(y)|^{p-2} (\hat{u}(x)-\hat{u}(y)) \;\ud y .
\end{align*}
 The first integral on the right hand side converges absolutely and does not depend on $\varepsilon >0$, so we only have to understand the limit as $\varepsilon \to 0+$ of the second integral on the right hand side. As $\hat{u}\in C^{1,\alpha} (\hat{\Omega})$, one has
\[
\hat{u} (x) - \hat{u} (y) = \nabla \hat{u} (x) \cdot (x-y) + r(x,y) 
\]
for some rest term satisfying
\[
|r(x,y)| = C\, |x-y|^{1+\alpha} . 
\]
Then 
\begin{align*}
\lefteqn{|\hat{u}(x)-\hat{u}(y)|^{p-2} ( \hat{u}(x)-\hat{u}(y) )} \\
& = |\nabla \hat{u}(x) \cdot (x-y) + r(x,y)|^{p-2} ( \nabla \hat{u}(x) \cdot (x-y) + r(x,y) ) \\
& = |\nabla \hat{u}(x) \cdot (x-y)|^{p-2} \nabla \hat{u}(x) \cdot (x-y) \\
& \phantom{=\ } + ( |\nabla \hat{u}(x) \cdot (x-y) + r(x,y)|^{p-2} - |\nabla \hat{u}(x) \cdot (x-y)|^{p-2} ) \nabla \hat{u}(x) \cdot (x-y) \\
& \phantom{=\ } + |\nabla \hat{u}(x) \cdot (x-y) + r(x,y)|^{p-2} r(x,y) .
\end{align*}
By antisymmetry of the integrand and by symmetry of the domain of integration, for every $\varepsilon \in ]0,r]$, 
\begin{align*}
  \int_{r> |x-y| > \varepsilon} k(x,y) |\nabla \hat{u}(x) \cdot (x-y)|^{p-2} \nabla \hat{u}(x) \cdot (x-y) \;\ud y & = 0 ,
\end{align*}
so that
\begin{align*}
\lefteqn{\int_{r> |x-y| > \varepsilon} k(x,y) |\hat{u}(x) - \hat{u}(y)|^{p-2} (\hat{u}(x)-\hat{u}(y)) \;\ud y} \\
& = \int_{r> |x-y| > \varepsilon} k(x,y) ( |\nabla \hat{u}(x) \cdot (x-y) + r(x,y)|^{p-2} - |\nabla \hat{u}(x) \cdot (x-y)|^{p-2} ) \nabla \hat{u}(x) \cdot (x-y) \;\ud y \\
& \phantom{=\ } + \int_{r> |x-y| > \varepsilon} k(x,y) |\nabla \hat{u}(x) \cdot (x-y) + r(x,y)|^{p-2} r(x,y) \;\ud y .
\end{align*}
For the integrands of the two integrals on the right hand side one has the estimate
\begin{align*}
 \lefteqn{k(x,y) \big| |\nabla \hat{u}(x) \cdot (x-y) + r(x,y)|^{p-2} - |\nabla \hat{u}(x) \cdot (x-y)|^{p-2} \big| \, | \nabla \hat{u}(x) \cdot (x-y) |} \\
 & + k(x,y) |\nabla \hat{u}(x) \cdot (x-y) + r(x,y)|^{p-2} |r(x,y)| \\
 & \leq C \, \frac{|x-y|^{p-3} \, |r(x,y)| \, |x-y|}{|x-y|^{N+\theta p}} + C \, \frac{|x-y|^{p-2} |r(x,y)|}{|x-y|^{N+\theta p}} \\
 & \leq C \, \frac{|x-y|^{p-1+\alpha}}{|x-y|^{N+\theta p}} ,
\end{align*}
so that, by the assumption on $\alpha$, they are absolutely integrable over $B(x,r)$.
\end{proof}

The last line in \eqref{eq.nonlocal.weak.example} is a weak form of a Neumann boundary condition on $\partial_{\mathbf S}\hat{\Omega}$. The equality $\mathcal N_{\partial_{\mathbf S}\hat{\Omega}}^{2-2\theta,p}\hat{u} = 0$ on $\partial_{\mathbf S}\hat{\Omega}$ {\em in the weak sense} by definition means, for every $\hat{v} \in W^{\theta,{\mathbf S}}_{p,2} (\hat{\Omega},\Omega)$, 
 \begin{multline} \label{eq.operator.n}
    \frac12 \iint_{{\mathbf S}} k(x,y) |\hat{u} (x) - \hat{u} (y)|^{p-2} (\hat{u} (x) - \hat{u} (y)) (\hat{v}(x) - \hat{v}(y)) \; \ud\mu (y) \; \ud\mu (x) = \\ 
    \int_{\hat{\Omega}} P.V.\int_{{\mathbf S}(x)} k(x,y) |\hat{u} (x) - \hat{u} (y)|^{p-2} (\hat{u} (x) - \hat{u} (y)) \; \ud\mu (y) \hat{v} (x) \; \ud\mu (x);
 \end{multline}
compare with \eqref{eq.nonlocal.weak.identification.bc} in the abstract problem (Theorem \ref{thm.identification.2}). We discuss this boundary condition in more detail in Example \ref{ex.9.3} below.

By Theorem \ref{thm.identification.2}, a function $u\in L^2 (\Omega )$ is a solution of the stationary problem \eqref{eq.nonlocal.weak.example} if and only if $(u,f)\in\partial_j\E$, that is, if and only if there is an elliptic extension $\hat{u}\in W^{\theta,{\mathbf S}}_{p,2} (\hat{\Omega} ,\Omega )$ of $u$, which is a weak solution of the first two equations in \eqref{eq.nonlocal.weak.example} and which satisfies the boundary condition in \eqref{eq.nonlocal.weak.example} (last equation) in a weak sense. By Corollary \ref{cor.neumann}, the nonlocal evolution equation \eqref{eq.nonlocal.weak.example.cp} admits for every initial value $u_0\in L^2 (\Omega )$ and for every right hand side $f\in L^2_{loc} ([0,\infty [ ; L^2 (\Omega ))$ a unique global solution $u\in C (\R_+ ; L^2 (\Omega ))\cap W^{1}_{2,loc} (]0,\infty [; L^2 (\Omega ))$. The solution is a strong solution on $]0,\infty [$, that is, $u\in W^1_{2,loc} ([0,\infty [;L^2 (\Omega))$ and $(u(t) , f(t) -\dot{u}(t))\in\dom{\partial_j\E}$ for almost every $t\geq 0$. For almost every $t>0$, the function $u(t)\in L^2 (\Omega )$ admits an elliptic extension $\hat{u} (t)\in W^{\theta,{\mathbf S}}_{p,2} (\hat{\Omega} ,\Omega )$ and the equalities in \eqref{eq.nonlocal.weak.example.cp} are satisfied in a weak sense. In particular, the negative $j$-subgradient $-\partial_j\E$ generates a strongly continuous semigroup $S^N$ of contractions on $L^2 (\Omega )$. We study the comparison principle, the maximum principle and ultracontractivity of this semigroup in the following sections.
\end{example}

\begin{remark}
 {\em In the special case when
 \[
 k(x,y)=\frac{C_{N,\theta,p}}{|x-y|^{N+\theta p}} ,
 \]
 where the constant $C_{N,\theta,p}$ is given by
 \begin{equation}\label{CNSP}
C_{N,\theta,p}:=\frac{\theta 2^{2\theta-1}\Gamma\left(\frac{p\theta+p+N-2}{2}\right)}{\pi^{\frac N2}\Gamma(1-\theta)} ,
 \end{equation}
 then the energy $\E$ is given by
 \begin{align*}\label{func-E}
 \mathcal E (\hat{u})=\frac{C_{N,\theta,p}}{2p} \iint_{\mathbf{S}} \frac{|\hat{u}(x)-\hat{u}(y)|^p}{|x-y|^{N+\theta p}}\;\ud x\;\ud y \quad (\hat{u}\in W^{\theta}_{p,2}(\hat{\Omega},\Omega)) .
 \end{align*}
Notice that for general $p\in ]1,\infty[$, the constant $C_{N,\theta,p}$ given in \eqref{CNSP} has been introduced in Warma \cite{War-NODEA} to ensure that for smooth functions $\hat{u}\in {\mathcal D} (\hat{\Omega}) = C_c^\infty (\hat\Omega)$ we have at the limit
\begin{equation*}
 \lim_{\theta\uparrow 1^-}\frac{C_{N,\theta,p}}{2p}\int_{\hat\Omega}\int_{\hat\Omega}\frac{|\hat{u}(x)-\hat{u}(y)|^p}{|x-y|^{N+\theta p}}\;\ud y\;\ud x=\frac{1}{p}\int_{\hat\Omega}|\nabla \hat{u}|^p\;\ud x,  
\end{equation*}
that is, at the limit we recover the energy of the well-known $p$-Laplace operator on $\hat{\Omega}$; see also Foghem \cite[Section 9]{Fo25}.

Throughout the rest of the paper, since we are not taking limits as $\theta\uparrow 1^-$, the constant $C_{N,\theta,p}$ is irrelevant. We have included it for completeness and for researchers who are interested to do further research in this direction.}
\end{remark}

\begin{example}[\bf Graphs] \label{ex.3.2}
 Let $\hat{\Omega}$ be a countable set equipped with a weighted counting measure with weight $m : \hat{\Omega} \to ]0,\infty [$, let $\Omega\subseteq\hat{\Omega}$ be a nonempty subset, and let $(b,0)$ be a {\em graph} on $\hat{\Omega}$, that is, $b: \hat{\Omega} \times\hat{\Omega} \to [0,\infty [$ is a symmetric function (see Keller, Lenz \& Wojciechowski \cite[Chapter 1, Definition 1.1]{KeLeWo21}). We say that two vortices $x$, $y\in\hat{\Omega}$ are {\em connected by an edge} (we write $x\sim y$ in this case) if $b(x,y) >0$, and we say that $(b,0)$ is a {\em locally finite graph} if for every $x\in\hat{\Omega}$ the set $\{y\in\hat{\Omega}\st y\sim x\}$ is finite. 

Let $p\in ]1,\infty [$ and $\ell^2_m (\Omega ) := \left\{ u : \Omega \to\R \st \sum_{x\in\Omega} m(x) |u(x)|^2\right\}$. For given $f\in \ell^2_m (\Omega )$, consider the stationary problem of finding $u\in \ell^2_m (\Omega )$ such that
\begin{equation} \label{eq.nonlocal.graph.example}
\begin{cases}
     \multicolumn{2}{l}{\text{there exists an extension of } u \text{ to a function } \hat{u}\in W^{(b,0)}_{p,2} (\hat{\Omega},\Omega ), \text{ that is, } \hat{u}|_{\Omega} = u , \text{ such that}} \\
    \displaystyle \sum_{y\in\hat{\Omega} \atop y\sim x} b(x,y) |\hat{u} (x) - \hat{u} (y)|^{p-2} (\hat{u} (x) - \hat{u} (y)) = f(x) & \text{for } x\in\Omega , \\
   \displaystyle \sum_{y\in\hat{\Omega} \atop y\sim x} b(x,y) |\hat{u} (x) - \hat{u} (y)|^{p-2} (\hat{u} (x) - \hat{u} (y)) = 0 & \text{for } x\in\hat{\Omega}\setminus\Omega , \\
   \displaystyle \sum_{x,y\in\hat{\Omega} \atop x\sim y} b(x,y) |\hat{u} (x) - \hat{u} (y)|^{p-2} (\hat{u} (x) - \hat{u} (y)) (\hat{v}(x) - \hat{v}(y)) = & \\
   \phantom{\ud x} \displaystyle = \sum_{x\in\hat{\Omega}} \sum_{y\in\hat{\Omega}\atop y\sim x} b(x,y) |\hat{u} (x) - \hat{u} (y)|^{p-2} (\hat{u} (x) - \hat{u} (y)) \hat{v}(x) & \text{for every } \hat{v} \in W^{(b,0)}_{p,2} (\hat{\Omega},\Omega ) . 
\end{cases}
\end{equation}
Consider also for a given initial value $u_0 \in \ell^2_m (\Omega )$ and a right hand side $f\in L^2_{loc} ([0,\infty [; \ell^2_m (\Omega ))$ the associated evolution problem
\begin{equation} \label{eq.nonlocal.graph.example.cp}
\begin{cases}
   \displaystyle  \partial_t \hat{u} (t,x) + \sum_{y\in\hat{\Omega} \atop y\sim x} b(x,y) |\hat{u} (t,x) - \hat{u} (t,y)|^{p-2} (\hat{u} (t,x) - \hat{u} (t,y)) = f(t,x) \quad &\text{for } (t,x)\in ]0,\infty [ \times \Omega , \\
   \displaystyle \sum_{y\in\hat{\Omega} \atop y\sim x} b(x,y) |\hat{u} (t,x) - \hat{u} (t,y)|^{p-2} (\hat{u} (t,x) - \hat{u} (t,y)) \; \ud\mu (y) = 0 \quad &\text{for } (t,x)\in ]0,\infty [ \times (\hat{\Omega}\setminus\Omega) , \\
   \displaystyle \sum_{x,y\in\hat{\Omega} \atop x\sim y} b(x,y) |\hat{u} (t,x) - \hat{u} (t,y)|^{p-2} (\hat{u} (t,x) - \hat{u} (t,y)) (\hat{v}(x) - \hat{v}(y)) = & \\
   \phantom{\ud x} \displaystyle = \sum_{x\in\hat{\Omega}} \sum_{y\in\hat{\Omega} \atop y\sim x} b(x,y) |\hat{u} (t,x) - \hat{u} (t,y)|^{p-2} (\hat{u} (t,x) - \hat{u} (t,y)) \hat{v} (x) & \text{for } t\in ]0,\infty [, \; \forall \hat{v} \in W^{(b,0)}_{p,2} (\hat{\Omega} , \Omega  , \\
   \displaystyle \hat{u} (0,x) = u_0 (x)&\text{for } x\in\Omega.
\end{cases}
\end{equation}
Here, 
 \begin{align*}
 W^{(b,0)}_{p,2} (\hat{\Omega} , \Omega ) & = \left\{ \hat{u} : \hat{\Omega} \to \R \st \hat{u}|_\Omega \in\ell^2_m (\Omega ) \text{ and } \sum_{x\in\hat{\Omega}} \sum_{y\in\hat{\Omega}} b(x,y) |\hat{u} (x) - \hat{u} (y)|^p <\infty \right\} .
 \end{align*}

Both problems \eqref{eq.nonlocal.graph.example} and \eqref{eq.nonlocal.graph.example.cp} fit into the setting of this article, when one defines the kernel $\Phi : \hat{\Omega}\times\hat{\Omega}\times\R \to [0,\infty [$ by $\Phi (x,y,s) = \frac{1}{p} \frac{b(x,y)}{m(x)m(y)} |s|^p$. This kernel satisfies the standard condition \eqref{cond.standard}, and it is continuously partially differentiable with respect to the third variable. The space $W^{(b,0)}_{p,2} (\hat{\Omega} , \Omega)$ is equal to the abstract space $W^{\Phi ,2} (\hat{\Omega} ,\Omega )$.  

Assume that every vortex $y\in\hat{\Omega}\setminus\Omega$ is {\em connected} to a vortex $x\in\Omega$ {\em by a finite path} in the sense that there exists a finite family $(x_i)_{0\leq i\leq n}$ in $\hat{\Omega}$ such that $x_0 = y$, $x_n = x$ and $x_i \sim x_{i+1}$ for every $0\leq i\leq n-1$. Then the set ${\mathbf S}_\Phi := \{ (x,y)\in\hat{\Omega}\times\hat{\Omega} \st x\sim y\}$
satisfies the thickness condition \eqref{cond.connected}. In fact, start with singletons $A^{(0)}_n = \{ x_n\}\subseteq\Omega$ such that $\bigcup_n A^{(0)}_n = \Omega$. Having chosen $A^{(j)}_n$ for some $j\geq 0$ and every $n\in\N$, let $B^{(j+1)}_n = \{ y\in\hat{\Omega} \st y \sim x_n\}$ be the set of vortices which are connected to $x_n$ by an edge, and let $A^{(j+1)}_n = \{y_n\}$ be singletons such that $\bigcup_n B^{(j+1)}_n = \bigcup_n A^{(j+1)}_n$. By the assumption that every vortex $y\in\hat{\Omega}\setminus\Omega$ is connected to a vortex $x\in\Omega$ by a finite path it follows that $\bigcup_{j\in\N}\bigcup_{n\in\N} A^{(j)}_n = \hat{\Omega}$. By Lemma \ref{lem.banach}, $W^{(b,0)}_{p,2} (\hat{\Omega} ,\Omega )$ is a Banach space.

The energy $\E_{(b,0)} : W^{(b,0)}_{p,2} (\hat{\Omega} , \Omega ) \to [0,\infty ]$ given by
 \begin{align*}
 \E_{(b,0)} (\hat{u} ) & = \frac{1}{2p} \sum_{x\in\hat{\Omega}} \sum_{y\in\hat{\Omega}} b(x,y) |\hat{u} (x) - \hat{u} (y)|^p \\
 & = \frac{1}{2p} \sum_{x,y\in\hat{\Omega} \atop x\sim y} b(x,y) |\hat{u} (x) - \hat{u} (y)|^p 
 \end{align*}
 is G\^ateaux differentiable (Lemma \ref{lem.identification}) and, for every $\hat{u}$, $\hat{v}\in W^{(b,0)}_{p,2} (\hat{\Omega} , \Omega )$, 
 \[
 \E_{(b,0)}'(\hat{u}) \hat{v} = \frac 12\sum_{x\in\hat{\Omega}} \sum_{y\in\hat{\Omega} \atop y\sim x} b(x,y) |\hat{u}(x) - \hat{u} (y)|^{p-2} (\hat{u}(x) - \hat{u} (y)) (\hat{v}(x) - \hat{v} (y)) . 
 \]

By Theorem \ref{thm.identification.2}, a function $u\in \ell^2_m (\Omega )$ is a solution of the stationary problem \eqref{eq.nonlocal.graph.example} if and only if $(u,f)\in\partial_j\E_{(b,0)}$, that is, there exists an elliptic extension $\hat{u}\in W^{(b,0)}_{p,2} (\hat{\Omega},\Omega )$ such that the first two equations hold. By Corollary \ref{cor.neumann}, the nonlocal evolution equation \eqref{eq.nonlocal.graph.example.cp} admits for every initial value $u_0\in \ell^2_m (\Omega )$ and for every right hand side $f\in L^2_{loc} ([0,\infty [ ; \ell^2_m (\Omega ))$ a unique global solution $u\in C (\R_+ ; \ell^2_m (\Omega ))\cap W^{1}_{2,loc} (]0,\infty [; \ell^2_m (\Omega ))$. The solution is a strong solution on $]0,\infty [$, that is, $(u(t) , f(t) -\dot{u}(t))\in\dom{\partial_j\E_{(b,0)}}$ for almost every $t> 0$. For almost every $t>0$, the function $u(t)\in L^2 (\Omega )$ admits an elliptic extension $\hat{u} (t)\in W_{p,2}^{(b,0)} (\hat{\Omega} ,\Omega )$ and the equalities in \eqref{eq.nonlocal.graph.example.cp} are satisfied. In particular, the negative $j$-subgradient $-\partial_j\E_{(b,0)}$ generates a strongly continuous semigroup $S^N$ of contractions on $\ell^2_m (\Omega )$.\\

Depending on the choice of $\Omega$, $\hat{\Omega}$ and the graph $(b,0)$, this example includes graph $p$-Laplace operators with Neumann exterior conditions, or Neumann boundary conditions (but what is the boundary of a graph?), or both.
\end{example}

\section{Comparison principles, maximum principles and nonlinear Dirichlet forms}\label{NDF}

Let $(\Omega ,{\mathcal A}, \mu )$ be a $\sigma$-finite measure space. We call a proper, lower semicontinuous, convex, densely defined function $\E : L^2 (\Omega ) \to ]-\infty , \infty ]$ a {\em (nonlinear) Dirichlet form} if the semigroup $S$ generated by its negative subgradient is {\em submarkovian}, that is, by definition, if the semigroup is {\em order preserving} and {\em $L^\infty$-contractive}. Here, a semigroup $S$ on $L^2 (\Omega )$ is called {\em order preserving} if, for every $u_0$, $v_0\in {L^2(\Omega)}$,
\begin{equation} \label{eq.order.preserving}
 u_0 \leq v_0 \quad \Rightarrow \quad \forall t\geq 0 : S(t)u_0 \leq S(t) v_0 ,
\end{equation}
and it is {\em $L^\infty$-contractive} if, for every $u_0$, $v_0\in {L^2(\Omega)}$ and for every $t\geq 0$,
\begin{equation} \label{eq.linfty.contractive}
    \| S(t) u_0 - S(t)v_0\|_{L^\infty (\Omega )} \leq \| u_0 - v_0\|_{L^\infty (\Omega )} .
\end{equation}
In the preceding inequality we interpret the $L^\infty$-norms as equal to $\infty$ if the functions under the norms are not in $L^\infty (\Omega )$. The semigroup being order preserving expresses a comparison principle, while the $L^\infty$-contractivity expresses a maximum principle. 

\begin{remark}
    Note that in our definition, Dirichlet forms are necessarily densely defined on $L^2 (\Omega )$, and therefore the semigroups they generate are defined on the entire space $L^2 (\Omega )$. It is possible to define Dirichlet forms which are not necessarily densely defined, but then one has to be careful with the definitions of order preserving and $L^\infty$-contractive semigroups and with the Beurling-Deny criteria (see below). We prefer to avoid any subtleties in this direction.
\end{remark}

When $\E$ is quadratic, then the subgradient is a linear selfadjoint, positive semidefinite operator, and the semigroup is a linear, strongly continuous contraction semigroup of selfadjoint operators. Quadratic Dirichlet forms (or, equivalently, symmetric, bilinear Dirichlet forms) are studied in the monographs by Fukushima \cite{FuOsTa11}, Ma \& R\"ockner \cite{MaRo92} and Bouleau \& Hirsch \cite{BoHi91}. The notion of nonlinear Dirichlet forms has been coined by B\'enilan \& Picard \cite{BePi79a}, was further studied in Barth\'elemy \cite{By96} and Cipriani \& Grillo \cite{CiGr03}, and recently in Claus \cite{Cl21,Cl23}, Brigati \cite{Br23}, Brigati \& Hartarsky \cite{BrHa24s}, Brigati \& Dello Schiavo \cite{BrDS25}, Puchert \cite{Pu25b} and Schmidt \& Zimmermann \cite{ScZi25}. In particular, the classical Beurling-Deny criteria characterizing order preserving and $L^\infty$-contractive semigroups have been generalized to the nonlinear situation: by \cite[Theorem 3.6]{CiGr03}, the semigroup $S$ generated by $-\partial \E$ is order preserving if and only if, for every $u$, $v\in L^2 (\Omega )$, 
\begin{equation} \label{eq-OP}
 \E \left(\frac 12\left( u+u\wedge v\right)\right ) + \E \left(\frac 12\left(v+u\vee v\right)\right) \leq \E (u) + \E (v),
\end{equation}
and it is $L^\infty$-contractive if and only if, for every $u$, $v\in L^2 (\Omega )$ and for every $\alpha >0$,
\begin{equation}\label{non-expan}
    \E(v+p_\alpha (u-v))+\E(u-p_\alpha (u-v))\leq \E(u)+\E(v) ,
\end{equation}
where
\[ 
p_\alpha (\tau):=\frac 12 \left[(\tau+\alpha)^+-(\tau-\alpha)^-\right].
\]
A function $p : \R\to\R$ is a {\em normal contraction} if $p$ is increasing, $p(0) = 0$ and if $p$ is Lipschitz continuous with Lipschitz constant $\leq 1$. The functions $p_\alpha$ defined above are normal contractions, and also the function $p(s) := \frac12 s^+$ is a normal contraction. For this normal contraction, one has
\begin{align*}
& u - \frac12 (u-v)^+ = \frac12 (u + u\wedge v), \text{ and} \\
& v + \frac12 (u-v)^+ = \frac12 (v + u\vee v) . 
\end{align*}
Hence, if for every $u$, $v\in L^2 (\Omega )$ and for every normal contraction $p$, 
\begin{equation} \label{eq.energy.contraction.p} 
    \E (u- p(u-v)) + \E (v+p(u-v)) \leq \E (u) + \E (v) ,
\end{equation}
then $\E$ is a Dirichlet form. Actually, the converse is true: if $\E$ is a Dirichlet form, then the  inequality \eqref{eq.energy.contraction.p} holds for every $u$, $v\in L^2 (\Omega )$ and every normal contraction $p$ (see e.g. B\'enilan \& Picard \cite{BePi79a}, Brigati \& Hartarsky \cite{BrHa24s}, Puchert \cite{Pu25}). 

The following result is taken from \cite[Theorem 3.6]{CiGr03}.

\begin{theorem}
Let $\E : L^2 (\Omega ) \to ]-\infty , +\infty ]$ be a Dirichlet form and let $S$ be the semigroup generated on $L^2 (\Omega )$ by $-\partial\mathcal E$. Then, for every $q\in [1,\infty[$, the semigroup $S$ extends from $L^2(\Omega) \cap L^q (\Omega )$ to a strongly continuous, order preserving contraction semigroup on $L^q (\Omega )$. The semigroup also extends to an order preserving semigroup of contractions on $\overline{L^2 (\Omega ) \cap L^\infty(\Omega)}^{L^\infty(\Omega)}$ which, however, is not always strongly continuous.
\end{theorem}

\begin{remark}
By B\'enilan \& Crandall \cite{BeCr91} (see in particular Lemma 7.1 therein), for every $N$-function $\Psi$ the semigroup $S$ extends from $L^2(\Omega) \cap L^\Psi (\Omega )$ to a strongly continuous, order preserving contraction semigroup on the Orlicz space $L^\Psi (\Omega )$, when this space is equipped with the Minkowski norm. More generally, for every normal Banach space $X\subseteq L^1(\Omega) + L^\infty (\Omega )$ the semigroup $S$ extends from $L^2 (\Omega ) \cap X$ to a strongly continuous, order preserving contraction semigroup on $X$ (see  \cite{BeCr91} for the definition of a normal Banach space).
\end{remark}

\begin{theorem}\label{theo-subm}
 Let $\Phi : \hat{\Omega} \times \hat{\Omega} \times \R \to [0,\infty ]$ be a function satisfying the standard condition \eqref{cond.standard}, let $W^{\Phi ,2} (\hat{\Omega} , \Omega )$ be the space defined in \eqref{eq.space.w}, let $\E : W^{\Phi ,2} (\hat{\Omega} , \Omega ) \to [0,\infty]$ be the energy defined in \eqref{energy.e}, let $j:W^{\Phi ,2} (\hat{\Omega} ,\Omega ) \to L^2 (\Omega )$ be the restriction operator, and assume that $j$ has dense range. Let $\E^{L^2}$ be the energy on $L^2 (\Omega )$ associated with $\E$ and $j$, that is, $\partial_j \E = \partial\E^{L^2}$. Then the energy $\E^{L^2}$ is a densely defined Dirichlet form, and therefore the semigroup $S^N$ generated by its negative subgradient is order preserving and $L^\infty$-contractive. 
\end{theorem}

The proof of this theorem relies on the following two fundamental lemmas for functions of one or two real variables. For the first lemma, see Barth\'elemy \cite[Lemme 4.2]{By96} or Velez Santiago \& Warma \cite[Lemma 3.3]{VSWa10}.

\begin{lemma}\label{Bar-Lem}
Let $F:\;\RR^2\to(-\infty,+\infty]$ be a convex and lower semicontinuous function with effective domain $\dom{F}:=\{(u,v)\in\RR^2 \st F(u,v)<\infty\}$. Assume that $U:=\mbox{Int}(\dom{F})$ is nonempty. Then, the following assertions are equivalent.
\begin{enumerate}
\item[(i)] For all $(u_0,u_1)$, $(v_0,v_1)\in\RR^2$ with $(u_0-v_0)(u_1-v_1)<0$,
\begin{equation}\label{eq-nec}
F(u_0,v_1)+F(v_0,u_1)\le F(u_0,u_1)+F(v_0,v_1).
\end{equation}
\item[(ii)] $\displaystyle \frac{\partial^2F}{\partial u_0\partial u_1}\le 0$ in $\mathcal D'(U)$ and, for every $(u_0,u_1)$, $(v_0,v_1)\in \dom{F}$ with $(u_0-v_0)(u_1-v_1)<0$,
\begin{equation}\label{eq-suf}
 (u_0, v_1), (v_0,u_1)\in \dom{F}.
\end{equation}
\end{enumerate}
\end{lemma}

\begin{lemma} \label{lem.normal.contraction}
Let $\Phi : \R \to [0,\infty ]$ be a convex function such that its effective domain $\dom{\Phi}$ contains an open, nonempty interval. Then the following assertions hold.
\begin{itemize}
\item[(a)] For every $u_0$, $u_1$, $v_0$, $v_1\in\R$ and for every normal contraction $p:\R\to\R$,
\begin{multline*}
\Phi ((u_0 - p(u_0 - v_0)) - (u_1 - p(u_1 - v_1))) - \Phi ((v_0 + p(u_0 - v_0)) - (v_1 + p(u_1 - v_1))) \\
\leq \Phi (u_0 - u_1) + \Phi (v_0 - v_1) .
\end{multline*}
\item[(b)] For every $u_0$, $u_1$, $v_0$, $v_1\in\R$,
\begin{multline} \label{eq.proof.order.preserving.a}
   \Phi (\frac12 (u_0 + u_0 \wedge v_0) - \frac12 (u_1 + u_1 \wedge v_1  ))  + \Phi ( \frac12 (u_0 + u_0 \vee v_0) - \frac12 (u_1 + u_1 \vee v_1))  \\
\leq \Phi (u_0 - u_1) + \Phi (v_0 - v_1) .
\end{multline}
\item[(c)] For every $u_0$, $u_1$, $v_0$, $v_1\in\R$ and for every $\alpha >0$,
\begin{multline} \label{eq.proof.nonexpan.a}
   \Phi((u_0 - p_\alpha (u_0-v_0))-(u_1 - p_\alpha (u_1-v_1)) + \Phi((v_0 + p_\alpha (u_0-v_0))-(v_1 - p_\alpha (u_1-v_1))) \\ 
   \leq \Phi (u_0 -u_1) + \Phi (v_0 -v_1) . 
\end{multline}
\end{itemize}
\end{lemma}

\begin{proof}
The function $\Phi$ is convex and therefore continuous in the interior of the effective domain and lower semicontinuous on the entire effective domain. Outside the effective domain, the function takes the value $\infty$, so that altogether $\Phi$ is lower semicontinuous. As a consequence, the function $F:\R^2 \to [0,\infty ]$, $F (u,v) := \Phi (u-v)$ is convex and lower semicontinuous, too. Since the effective domain of $\Phi$ has nonempty interior in $\R$, the effective domain of $F$ has nonempty interior in $\R^2$. On the interior of the effective domain,  one has $\frac{\partial^2 F}{\partial u \partial v} \leq 0$ in the sense of distributions because of the convexity of $\Phi$. Therefore, the function $F$ satisfies the assumptions and the condition (ii) of Lemma \ref{Bar-Lem}. By Lemma \ref{Bar-Lem}, the function $F$ also satisfies the condition (i) of Lemma \ref{Bar-Lem}.

(a) Let $p:\R \to\R$ be a normal contraction, and let $u_0$, $u_1$, $v_0$, $v_1\in\R$. In order to prove the inequality in (a), we define $\lambda_0$, $\lambda_1 \in \R$ by
\begin{equation*}
\lambda_0 := \begin{cases} \frac{p(u_0 - v_0)}{u_0 - v_0} & \text{ if } u_0\not= v_0 , \\ 0 & \text{ if } u_0 = v_0 , \end{cases} \text{ and } \\
\lambda_1 := \begin{cases} \frac{p(u_1 - v_1)}{u_1 - v_1} & \text{ if } u_1\not= v_1 , \\ 0 & \text{ if } u_1 = v_1 . \end{cases}
\end{equation*}
As $p$ is increasing and Lipschitz continuous with Lipschitz constant $\leq 1$, and as $p(0)=0$, then $\lambda_0$, $\lambda_1\in [0,1]$. Moreover, 
\begin{equation*}\label{conv}
u_i- p (u_i - v_i) = \lambda_i v_i + (1-\lambda_i) u_i \text{ and } v_i+p(u_i - v_i)=\lambda_i u_i+(1-\lambda_i) v_i ,\; i=0,1.
\end{equation*}
From here and from the convexity of $\Phi$ follows 
\begin{align}\label{A1}
& \Phi((u_0 - p(u_0 - v_0)) -(u_1 - p(u_1-v_1))) \notag\\
= & \Phi(\lambda_0 v_0+(1-\lambda_0) u_0 - \lambda_1 v_1 - (1-\lambda_1) u_1) \notag\\
\le & \lambda_0\Phi(v_0-\lambda_1 v_0 -(1-\lambda_1) u_1)
+ (1-\lambda_0) \, \Phi(u_0-\lambda_1 v_1-(1-\lambda_1) u_1) \notag\\
\le & \lambda_0\lambda_1\Phi(v_0-v_1) + \lambda_0 (1-\lambda_1) \Phi(v_0-u_1) \notag\\
& +(1-\lambda_0) \lambda_1 \Phi(u_0-v_1) + (1-\lambda_0) (1-\lambda_1) \Phi(u_0-u_1)
\end{align}
and
\begin{align}\label{A2}
& \Phi((v_0 + p(u_0-v_0)) - (v_1 + p(u_1-v_1))) \notag\\
= & \Phi(\lambda_0 u_0 + (1-\lambda_0) v_0 - \lambda_1 u_1 -(1-\lambda_1) v_1) \notag\\
\le & \lambda_0\Phi (u_0-\lambda_1 u_1-(1-\lambda_1) v_1)
+ (1-\lambda_0) \Phi(v_0 - \lambda_1 u_1 -(1-\lambda_1) v_1) \notag\\
\le & \lambda_0 \lambda_1 \Phi(u_0-u_1) +\lambda_0(1-\lambda_1) \Phi(u_0-v_1) \notag\\
& + (1-\lambda_0) \lambda_1 \Phi(v_0-u_1) + (1-\lambda_0) (1-\lambda_1) \Phi(v_0 -v_1) .
\end{align}
Using Lemma \ref{Bar-Lem}, we obtain that
\begin{align}\label{A3}
\Phi (v_0-u_1) + \Phi (u_0-v_1) \le \Phi (u_0-u_1) + \Phi (v_0-v_1).
\end{align}
Combining \eqref{A1}-\eqref{A3}, it follows that
\begin{align}
 &  \Phi ((u_0 - p(u_0-v_0)) - (u_1 - p(u_1-v_1)) + \Phi((v_0 - p(u_0-v_0)) -(v_1 - p(u_1-v_1))) \notag\\
 \le & \lambda_0\lambda_1 \, \Big[\Phi(v_0-v_1) + \Phi(u_0-u_1) \Big] + \lambda_0(1-\lambda_1) \, \Big[ \Phi (u_0-u_1) + \Phi (v_0-v_1) \Big] \notag\\
 & + (1-\lambda_0) \lambda_1\, \Big[\Phi (u_0-u_1) + \Phi (v_0-v_1)\Big] \notag\\
 & + (1-\lambda_0) (1-\lambda_1) \, \Big[ \Phi (u_0-u_1) + \Phi (v_0-v_1) \Big]\notag\\
 = & \Phi (u_0-u_1) + \Phi (v_0-v_1),
\end{align}
and we have shown (a).
The inequality in (b) follows from (a) by taking the normal contraction $p(s) = \frac12 s^+$. Similarly, the inequality in (c) follows from (a) by taking the normal contractions $p_\alpha$. 
\end{proof}

\begin{proof}[\bf Proof of Theorem \ref{theo-subm}]
We first prove that for every $\hat{u}$, $\hat{v}\in W^{\Phi ,2} (\hat{\Omega} , \Omega )$ and for every normal contraction $p:\R\to\R$,
\begin{equation} \label{e.p}
 \E ( \hat{u} - p(\hat{u} - \hat{v})) + \E (\hat{v} + p(\hat{u} - \hat{v})) \leq \E (\hat{u}) + \E (\hat{v}) .
\end{equation}
In fact, this inequality is a simple consequence of Lemma \ref{lem.normal.contraction} (a) from which it follows that for every $\hat{u}$, $\hat{v}\in W^{\Phi ,2} (\hat{\Omega} , \Omega )$, for every normal contraction $p:\R\to\R$, and for every $x$, $y\in\hat{\Omega}$,
\begin{align*}
    \Phi (x,y,(\hat{u} (x) - p(\hat{u} (x) - \hat{v}(x))) - (\hat{u} (y) - p(\hat{u} (y) - \hat{v} (y))))\\
    + \Phi (x,y,(\hat{v} (x) + p(\hat{u} (x) - \hat{v} (x))) - (\hat{v} (y) + p(\hat{u} (y) - \hat{v}(y)))) \\
\leq \Phi (x,y,\hat{u}(x) - \hat{u} (y)) + \Phi (x,y,\hat{v} (x) - \hat{v}(y)) .
\end{align*}
The inequality \eqref{e.p} follows from here by integrating over $\hat{\Omega}\times\hat{\Omega}$ with respect to the product measure $\mu\otimes\mu$.

We deduce that for every $u$, $v\in L^2 (\Omega )$ and for every normal contraction $p:\R\to\R$, 
\begin{equation} \label{el2.p}
    \E^{L^2} (u - p(u-v)) + \E^{L^2} (v+p(u-v)) \leq \E^{L^2} (u) + \E^{L^2} (v) .
\end{equation}
Indeed,  let $u$, $v\in L^2 (\Omega )$. If $\E^{L^2} (u) = \infty$ or $\E^{L^2} (v) = \infty$, then the preceding inequality \eqref{el2.p} is obviously true. So let us assume that both $u$ and $v$ belong to the effective domain of $\E^{L^2}$. Let $\varepsilon >0$. By the representation \eqref{eq.energy.h} of the energy $\E^{L^2}$, there exist $\hat{u}$ and $\hat{v}$ in $W^{\Phi ,2} (\hat{\Omega} , \Omega )$ such that 
\begin{equation}
\begin{cases}
 j(\hat{u}) = u,\, j(\hat{v} ) = v , \\
 \E (\hat{u}) \leq \E^{L^2} (u) + \varepsilon \text{ and } \E (\hat{v}) \leq \E^{L^2} (v) + \varepsilon .
\end{cases}
\end{equation}
Note that due to the special form of the restriction operator $j$,
\[
j(\hat{u} - p(\hat{u} - \hat{v})) = u - p(u-v) \text{ and } j(\hat{v} + p(\hat{u} - \hat{v})) = v + p(u-v) .
\]
By the representation \eqref{eq.energy.h}, by the inequality \eqref{e.p} and by the choice of the functions $\hat{u}$ and $\hat{v}$,
\begin{align*}
  \E^{L^2} (u - p(u,v)) + \E^{L^2} (v+p(u-v)) & \leq \E (\hat{u} - p(\hat{u} - \hat{v})) + \E (\hat{v} + p(\hat{u} - \hat{v})) \\
  & \leq \E (\hat{u} ) - \E (\hat{v}) \\
  & \leq \E^{L^2} (u) + \E^{L^2} (v) + 2\varepsilon .
\end{align*}
Since $\varepsilon >0$ is arbitrary, the inequality \eqref{el2.p} follows. As mentioned around the inequality \eqref{eq.energy.contraction.p}, the inequality \eqref{el2.p} for every $u$, $v\in L^2 (\Omega )$ and every normal contraction $p:\R\to\R$ implies that $\E^{L^2}$ is a Dirichlet form. 
\end{proof}

\section{Interlude: Regular Dirichlet forms and the capacity induced by the space $W^{\Phi ,2} (\hat{\Omega} , \Omega )$}\label{sec5-cap}

The material in this section basically follows the work of Claus \cite{Cl21,Cl23}. It is useful for our study in the next sections.
We assume that $(\hat{\Omega} , {\mathcal B} (\hat{\Omega}) ,\mu )$ is a $\sigma$-finite Hausdorff topological measure space, that is, $\hat{\Omega}$ is a Hausdorff topological space, ${\mathcal B} (\hat{\Omega} )$ is the Borel $\sigma$-algebra on $\hat{\Omega}$, and $\mu$ is a $\sigma$-finite Borel measure which has full support, that is, there is no nonempty open subset $U\subseteq\hat{\Omega}$ with $\mu (U) = 0$. 

Let $\Phi : \hat{\Omega} \times \hat{\Omega} \times\R \to [0,\infty ]$ be a function satisfying the standard condition \eqref{cond.standard}, and let the space $W^{\Phi ,2} (\hat{\Omega} , \Omega )$ and the energy $\E$ be defined as in \eqref{eq.space.w} and \eqref{energy.e}, respectively. 

We let $\tilde{\Omega}$ be a compactification of $\hat{\Omega}$, that is, $\tilde{\Omega}$ is compact and there is a continuous embedding $\hat{\Omega}\to\tilde{\Omega}$ with dense range. Functions in $C(\tilde{\Omega})$ are completely determined by their values on $\hat{\Omega}$. We call \(C(\tilde{\Omega})\) a {\em core} of \(\E\) if \(W^{\Phi ,2} (\hat{\Omega} , \Omega ) \cap C(\tilde{\Omega})\) is dense in \(W^{\Phi ,2} (\hat{\Omega} , \Omega )\), and we call \(C(\tilde{\Omega})\) {\em regular} if 
\begin{align*}
    \overline{\langle W^{\Phi ,2} (\hat{\Omega} , \Omega ) \cap C(\tilde{\Omega}), \{1\} \rangle}^{\|\cdot\|_{\infty}}=\langle \overline{W^{\Phi ,2} (\hat{\Omega} , \Omega ) \cap C(\tilde{\Omega})}^{\|\cdot\|_{\infty}}, \{1\} \rangle= C(\tilde{\Omega}) .
\end{align*}
We call \(\mathcal{E}\) {\em regular} (with respect to $\tilde{\Omega}$), if $C(\tilde{\Omega})$ is a regular core. Of course, this definition depends on the compactification $\tilde{\Omega}$. The aim always is to choose $\tilde{\Omega}$ (and $C(\tilde{\Omega})$) as small as possible, so that $C(\tilde{\Omega})$ still is a regular core. Whatever the compactification, the measure $\mu$ can be extended to $\tilde{\Omega}$ by assuming that $\tilde{\Omega}\setminus\hat{\Omega}$ is a null set. We denote this extended measure again by $\mu$. 

For every Borel set $B\subseteq\tilde{\Omega}$ we define
\[
\Cap_\E (B) := \inf \{ \| \hat{w}\|_{W^{\Phi , 2} (\hat{\Omega} , \Omega )}^2 + \mu ((B\cap\hat{\Omega})\setminus\Omega ) \st \hat{w} \geq 1 \,\, \mu-\text{almost everywhere on an open set } U \supseteq B \} .
\]
The function $\Cap_\E$ is called the {\em capacity} associated with the energy $\E$ or, as the energy $\E$ is determined by the function $\Phi$ and as $\Phi$ determines the space, the capacity associated with the space $W^{\Phi ,2} (\hat{\Omega} ,\Omega )$. Note carefully that the capacity is defined for subsets of the compactification $\tilde{\Omega}$! 

The following lemma has been shown in \cite{Cl23} in the setting $\Omega = \hat{\Omega}$, but it also holds in the general case $\Omega\subseteq\hat{\Omega}$ and every admissible choice of $\tilde{\Omega}$.

\begin{lemma}[{\cite{Cl23}}] \label{lem.capacity}
The following assertions hold.
 \begin{itemize}
     \item[(a)] For every  Borel sets $B\subseteq\tilde{B}\subseteq\tilde{\Omega}$, 
     \[
     \Cap_\E (B) \leq \Cap_\E (\tilde{B}) .
     \]
     \item[(b)] For every Borel sets $B_1$, $B_2\subseteq\tilde{\Omega}$, 
     \[
     \Cap_\E (B_1\cup B_2) \leq \Cap_\E (B_1) + \Cap_\E (B_2) .
     \]
     \item[(c)] For every sequence $(B_n)$ of Borel subsets of $\tilde{\Omega}$, 
     \[
     \Cap_\E (\bigcup_n B_n) \leq \sum_n \Cap_\E (B_n).
     \]
     \item[(d)] For every decreasing sequence $(K_n)$ of compact subsets of $\tilde{\Omega}$, 
     \[
     \Cap_\E (\bigcap_n K_n) = \inf_n \Cap_\E (K_n) .
     \]
     \item[(e)] For every measurable set $B\subseteq\tilde{\Omega}$, 
     \[
     \mu (B) \leq \Cap_\E (B) \leq \mu ((B\cap\hat{\Omega})\cup\Omega ) .
     \]
 \end{itemize}
\end{lemma}

\begin{proof}
As $\Omega$ is not necessarily equal to $\hat{\Omega}$, the setting is not exactly the same as in \cite{Cl23}; for example, the restriction operator $j$ is in general not injective from the Dirichlet space $W^{\Phi ,2} (\hat{\Omega} , \Omega )$ into $L^2 (\Omega )$. This may explain the additional term $\mu ((B\cap\hat{\Omega})\setminus\Omega )$ in the definition of the capacity. However, the proofs from \cite{Cl23} can be easily adapted. For the assertions (a) and (b), see the proof of \cite[Lemma 5.2]{Cl23}, and for assertion (c), see the proof of \cite[Lemma 5.3]{Cl23}. For assertion (d), see the proof of \cite[Lemma 5.4]{Cl23}. The first inequality in assertion (e) follows from noting that every function $\hat{w}\in W^{\Phi ,2} (\hat{\Omega} , \Omega )$ which is greater or equal than $1$ on an open set $U\supseteq B$ satisfies $\| \hat{w}\|_{W^{\Phi , 2}}^2 \geq \| \hat{w}\|_{L^2 (\Omega )}^2 \geq \mu (B\cap\Omega)$, and therefore $\Cap_\E (B) \geq \mu (B\cap\Omega) + \mu ((B\cap\hat{\Omega})\setminus\Omega ) = \mu (B\cap\hat{\Omega})$; if there is no such function, then $\Cap_\E (B) = \infty$, and the first inequality is trivially true.  The second inequality in assertion (e) is trivially true if $\mu (\Omega ) = \infty$. If $\mu (\Omega ) <\infty$, it suffices to note that the constant function $1$ is an element in $W^{\Phi ,2} (\hat{\Omega} , \Omega )$, that is, it is an admissible function in the definition of the capacity and $\| 1\|_{W^{\Phi ,2}}^2 = \mu (\Omega )$. As a consequence, $\Cap_\E (B) \leq \mu (\Omega ) + \mu ((B\cap\hat{\Omega})\setminus\Omega ) = \mu ((B\cap\hat{\Omega})\cup \Omega )$. 
\end{proof}

We call a function $\hat{w} : \tilde{\Omega} \to \R$ {\em quasi-continuous} if for every $\varepsilon >0$ there exists an open set $U\subseteq\tilde{\Omega}$ such that $\Cap_\E (U) \leq \varepsilon$ and  $\hat{w}$ is continuous on $\tilde{\Omega} \setminus U$. A function $\hat{w}\in L^0 (\hat{\Omega})$ is quasi-continuous if there exists a representative which has a quasi-continuous extension to $\tilde{\Omega}$. 

We say that a set $B\subseteq\tilde{\Omega}$ is {\em polar} if $\Cap_\E (B) = 0$, and some property $P(x)$ depending on $x\in\tilde{\Omega}$ holds {\em quasi everywhere} if there exists a polar set $B\subseteq\tilde{\Omega}$ such that $P(x)$ holds for every $x\in\tilde{\Omega}\setminus B$. 

\begin{lemma}[{\cite[Corollary 6.4]{Cl23}}] \label{lem.quasi.almost}
Two quasi-continuous representatives of functions in $W^{\Phi,2} (\hat{\Omega},\Omega )$ are equal almost everywhere if and only if they are equal quasi everywhere. 
\end{lemma}

\begin{lemma}[{\cite[Theorem 6.8]{Cl23}}] \label{lem.subsequence}
  Every convergent sequence in $W^{\Phi ,2} (\hat{\Omega},\Omega )$ admits a subsequence which converges pointwise quasi everywhere. 
\end{lemma}

\begin{lemma}[{\cite[Theorem 6.8 and Corollary 6.9]{Cl21}}] \label{lem.quasi.continuous}
    Every function in $\overline{W^{\Phi ,2} (\hat{\Omega} , \Omega ) \cap C (\tilde{\Omega} )}^{\|\cdot\|_{W^{\Phi ,2}(\hat\Omega,\Omega)}}$ is quasi-con\-ti\-nuous. If $W^{\Phi , 2} (\hat{\Omega} , \Omega ) \cap C (\tilde{\Omega} )$ is dense in $W^{\Phi , 2} (\hat{\Omega} , \Omega )$, then every element of $W^{\Phi , 2} (\hat{\Omega} , \Omega )$ is quasi-continuous.  
\end{lemma}

We say that a Borel measure $\nu$ on $\tilde{\Omega}$ is {\em absolutely continuous with respect to the capacity $\Cap_\E$} if $\nu (B) = 0$ for every polar Borel set $B\subseteq\tilde{\Omega}$. By Lemma \ref{lem.capacity} (e), the measure $\mu$ on $\tilde{\Omega}$ is absolutely continuous with respect to the capacity $\Cap_\E$, and hence every measure which is absolutely continuous with respect to $\mu$ is also absolutely continuous with respect to the capacity $\Cap_\E$. However, there may also be measures $\nu$ which are singular with respect to $\mu$, but absolutely continuous with respect to the capacity $\Cap_\E$.

\begin{example} \label{ex.5.1}
    If $W^{\Phi ,2} (\hat{\Omega} ,\Omega )\subseteq C (\tilde{\Omega} )$ with (necessarily) continuous embedding, then $\| \hat{w}\|_{C(\tilde{\Omega})} \leq C\, \| \hat{w}\|_{W^{\Phi ,2}(\hat\Omega,\Omega)}$ for every $\hat{w}\in W^{\Phi ,2} (\hat{\Omega} , \Omega )$. This implies that every singleton $B = \{ x\}$ has positive capacity, and hence the empty set is the only polar subset of $\tilde{\Omega}$. In this case, every Borel measure $\nu$ on $\tilde{\Omega}$ is absolutely continuous with respect to the capacity. 
\end{example}

\begin{example} \label{ex.5.2.a}
    In the context of Example \ref{ex.3.1}, that is, when $\Omega\subseteq \hat{\Omega} \subseteq \R^N$ are open subsets with Lipschitz continuous boundary, we let $\tilde{\Omega}$ be the closure of $\hat{\Omega}$ in the one-point compactification of $\R^N$. This natural compactification $\tilde{\Omega}$ of $\hat{\Omega}$ contains the Euclidean boundary of $\hat{\Omega}$, a wishful situation in applications. When $\theta\in ]0,1[$ and $p\in ]1,\infty [$, then the space 
    \[
     W^{\theta ,{\mathbf S}}_{p,2} (\hat{\Omega} ,\Omega ) = \left\{ \hat{w}\in L^0 (\hat{\Omega} ) \st \hat{w}|_\Omega \in L^2 (\Omega ) \text{ and } \iint_{{\mathbf S}} \frac{| \hat{w}(x) - \hat{w}(y)|^p}{|x-y|^{N+\theta p}} \; \ud\mu (y) \; \ud\mu (x) <\infty \right\} 
    \]
    embeds continuously into $C(\tilde{\Omega})$ as soon as $\theta p >N$ (compare, for example, with \cite[Chapter 1]{Gris}, even when $\Omega\subsetneq\hat{\Omega})$. In this case, every Borel measure $\nu$ on $\tilde{\Omega}$ is absolutely continuous with respect to the capacity. This includes, for example, the surface measure on the boundary $\partial\hat{\Omega}$, that is, the $(N-1)$-dimensional Hausdorff measure $\mathcal H_{N-1}$ on the boundary.
\end{example}

\begin{example} \label{ex.5.2}
We have the following situation regarding the boundary of $\hat{\Omega}$. We discuss the case $\Omega=\hat{\Omega}\subseteq\mathbb R^N$. We say that 
 the boundary $\partial\hat{\Omega}$ is a $d$-set, for some $d\in ]0,N]$, if there are two constants $0<C_1\le C_2$ such that for every $x\in \partial\hat{\Omega}$ and $r\in ]0,1]$,
\begin{equation}\label{d-set}
 C_1r^d\le \mathcal H_d( \partial\hat{\Omega}\cap B(x,r))\le C_2r^d, 
\end{equation}
where $\mathcal H_{d}$ is the restriction to $\partial\hat{\Omega}$ of the $d$-dimensional Hausdorff measure $\mathcal H_{d}$.

Let $p\in ]1,2]$ and $\theta\in ]0,1[$. We have the following two situations. 

\begin{enumerate}
\item Assume that $\hat{\Omega}$ is an $N$-set in the sense of \eqref{d-set}. Then, by Jonsson \& Wallin \cite[Theorem 1, pp 103]{JW1984}, $\hat{\Omega}$ has the extension property in the sense of \eqref{S0} below. Moreover, by Warma \cite[Section 4]{Wa15} the Borel measure $\mathcal H_N$ on $\tilde \Omega$ is absolutely continuous with respect to the capacity.

\item  Now assume that $\hat{\Omega}$ has the extension property in the sense of \eqref{S0} below and that $\partial \hat{\Omega}$ is a $d$-set in the sense of \eqref{d-set} for some $d\in]0,N[$. 
Then, $\mathcal H_{d}$ is absolutely continuous with respect to the capacity. 
This follows from the fact that if $p\in]1,2]$ and $0<\beta:=\theta-\frac{N-d}{p}<1$, then we have the continuous embedding
\begin{equation}\label{CE-1}
W^{\theta ,{\mathbf S}}_{p,2} (\hat{\Omega} )\hookrightarrow L^{p^\star}(\partial\hat{\Omega},\mathcal H_d),    
\end{equation}
where $p^\star:=\frac{pd}{d-\beta p}>p$. In particular, if $\mathcal H_d(\partial\hat{\Omega})<\infty$, then we have the continuous embedding
\begin{equation}\label{CE-2}
W^{\theta ,{\mathbf S}}_{p,2} (\hat{\Omega} )\hookrightarrow L^r(\partial\hat{\Omega},\mathcal H_d), \;\forall\; r\in [1,p^\star].   
\end{equation}
We refer to Jonsson \& Wallin \cite[Theorem 1, pp 103]{JW1984} or Danielli,  Garofalo \& Nhieu \cite[Theorem 11.1]{DGN} for the proofs of \eqref{CE-1} and \eqref{CE-2} and more details on this topic. Observe that the assumptions $\theta\in ]0,1[$ and $0<\theta-\frac{N-d}{p}<1$ imply that $0<\frac{N-d}{p}<\theta<1$. This is not a restriction on $\theta$, since if $0<\theta\le \frac{N-d}{p}<1$, then $\Cap_\E(\partial\hat{\Omega})=0$ so that functions in $W^{\theta ,{\mathbf S}}_{p,2} (\hat{\Omega} )$ are zero quasi-everywhere on $\partial\hat{\Omega}$ (see Warma \cite[Section 4]{Wa15}).
\end{enumerate}

\begin{itemize}
\item If $\hat{\Omega}$ has a Lipschitz continuous boundary $\partial\hat{\Omega}$, then it has the extension property and
$d$ in  \eqref{d-set} is the Hausdorff dimension of $\partial\hat{\Omega}$ which 
is given by $d=N-1$. Hence, $\mathcal H_{N-1}$ is absolutely continuous with respect to the capacity.

\item If $\partial\hat{\Omega}$ is the von Koch curve in dimension $N=2$, then $\hat{\Omega}$ has the extension property and $d$ in \eqref{d-set} is  the Hausdorff dimension of $\partial\hat{\Omega}$ which is given by $d:=\ln(4)/\ln(3)$. Hence, $\mathcal H_{d}$ is absolutely continuous with respect to the capacity. This is also the case for several domains with a special fractal geometry.

\item  But let us observe that if $\hat{\Omega}=]0,1[\setminus {\bf C}$, where ${\bf C}$ is the Cantor set, then the Hausdorff dimension of the boundary ${\bf C}$ is given by $d:=\ln(3)/\ln(2)$. It is known that {\bf C} is a $d$-set, but $\hat{\Omega}$ does not have the extension property. 
In addition,  the restriction to {\bf C} of the $d$-dimensional Hausdorff measure $\mathcal H_d$ is not absolutely continuous with respect to the capacity (see Arendt \& Warma \cite{ArWa-Pot} for the case $p=2$). Other interesting examples in this direction are contained in \cite{ArWa-Pot}.
\end{itemize}
\end{example}

\section{Potentials, Robin exterior conditions and Robin boundary conditions}\label{sec.robin}

Let $(\hat{\Omega} , {\mathcal B} (\hat{\Omega}) , \mu)$ be a $\sigma$-finite Hausdorff topological measure space, as in the previous section, let $\Phi : \hat{\Omega}\times\hat{\Omega} \times\R \to [0,\infty ]$ be a function satisfying the standard condition \eqref{cond.standard}, let $W^{\Phi ,2} (\hat{\Omega} , \Omega )$ be the space defined in \eqref{eq.space.w}, and let $\E : W^{\Phi ,2} (\hat{\Omega} ,\Omega )\to [0,\infty]$ be the energy defined in \eqref{energy.e}, $j:W^{\Phi ,2} (\hat{\Omega} ,\Omega )\to L^2 (\Omega )$ the restriction operator and $\partial_j\E$ the $j$-subgradient on $L^2 (\Omega )$. Let $\tilde{\Omega}$ be a compactification of $\hat{\Omega}$, and assume that the associated energy function $\E$ is regular with respect to $\tilde{\Omega}$, that is, $C (\tilde{\Omega})$ is a regular core. In particular, by Lemma \ref{lem.quasi.continuous}, every element of $W^{\Phi,2} (\hat{\Omega},\Omega )$ is quasi-continuous (has a quasi-continuous extension) on $\tilde{\Omega}$. 

In this situation, let $B:\tilde{\Omega}\times\R \to [0 ,\infty]$ be a function such that 
\begin{equation}\label{cond.b}
\begin{cases}
     B(x,0) = 0 \text{ for every } x\in\tilde{\Omega} , \\
B(x,\cdot ) \text{ is lower semicontinuous and convex for every } x\in\tilde{\Omega} , \\
B(\cdot , s) \text{ is measurable for every } s\in\R ,
\end{cases}
\end{equation}
and let $\nu$ be a Borel measure on $\tilde{\Omega}$ which is absolutely continuous with respect to the capacity $\Cap_\E$. Define the energy $\E_{B,\nu} : W^{\Phi ,2} (\hat{\Omega} , \Omega ) \to [0 ,\infty ]$ by
\[
\E_{B,\nu} (\hat{u}) := \E (\hat{u}) + \int_{\tilde{\Omega}} B(x,\hat{u}(x)) \; \ud\nu (x) ,
\]
so that $\E_{B,\nu}$ is an additive perturbation of the energy $\E$. In this definition, the integral is to be understood with respect to a quasi-continuous representative of $\hat{u}$. The quasi-continuous representative is defined on $\tilde{\Omega}$. 

\begin{theorem} \label{thm.eb.lsc}
Under the assumptions in this section, the energy $\E_{B,\nu}$ is well-defined, proper, lower semicontinuous, convex and $j$-elliptic. The $j$-subgradient $\partial_j\E_{B,\nu}$ on $L^2 (\Omega )$ is equal to the subgradient $\partial\E_{B,\nu}^{L^2}$ of a proper, lower semicontinuous and convex function $\E_{B,\nu}^{L^2} : L^2 (\Omega ) \to [0,\infty ]$. The negative $j$-subgradient generates a strongly continuous semigroup $S^{B,\nu}$ of (nonlinear) contractions on $\overline{j(\dom{\E_{B,\nu}})}^{L^2(\Omega)}$. More precisely, for every $u_0\in\overline{j(\dom{\E_{B,\nu}})}^{L^2(\Omega)}$ and every $f\in L^2 (0,T;L^2 (\Omega ))$ the abstract gradient system 
\begin{equation} \label{cp.j.b}
    \dot{u} + \partial_j\E_{B,\nu} (u) \ni f \text{ on } [0,\infty [ , \qquad u(0) = u_0 ,
\end{equation}
admits a unique strong solution $u\in C([0,\infty [ ; L^2 (\Omega )) \cap W^{1}_{2,loc} (]0,\infty [ ; L^2 (\Omega ))$.
\end{theorem}

\begin{proof}
By Theorem \ref{thm.e.lsc}, the function $\E$ is continuous, convex and $j$-elliptic on $W^{\Phi ,2} (\hat{\Omega} ,\Omega )$. It suffices to show that the function ${\mathcal B} : \hat{u} \mapsto \int_{\hat{\Omega}} B(x,\hat{u} (x))\; \ud\nu (x)$ is well-defined, proper, lower semicontinuous, and convex on $W^{\Phi ,2} (\hat{\Omega} ,\Omega )$. In order to see that ${\mathcal B}$ is well-defined, it suffices to note that any two quasi-continuous representatives of some function $\hat{u}\in W^{\Phi,2} (\hat{\Omega},\Omega )$ coincide quasi everywhere (Lemma \ref{lem.quasi.almost}), and that $\nu$ does not charge any polar sets. The function $\mathcal B$ is proper as ${\mathcal B} (0) = 0$ by the first line in condition \eqref{cond.b}. The function $\mathcal B$ is convex by the convexity of $B$ in the second variable. In order to prove lower semicontinuity, let $(\hat{u}_n)$ be a convergent sequence in $W^{\Phi ,2} (\hat{\Omega},\Omega )$, and $\hat{u}:=\lim_{n\to\infty} \hat{u}_n$. By Lemma \ref{lem.banach} and Lemma \ref{lem.subsequence}, $(\hat{u}_n)$ has a subsequence (again denoted by $(\hat{u}_n)$) which converges quasi everywhere, that is, there exists a polar set $A\subseteq\tilde{\Omega}$ such that $\lim_{n\to\infty} \hat{u}_n (x) = \hat{u} (x)$ for every $x\in\tilde{\Omega}\setminus A$. Now, as $\nu$ is absolutely continuous with respect to the capacity, by the continuity of $B$ in the second variable and by Fatou's lemma,
\begin{align*}
    {\mathcal B} (\hat{u}) & = \int_{\tilde{\Omega}} B(x,\hat{u}(x))\; \ud\nu (x) \\
    & = \int_{\tilde{\Omega}\setminus A} B(x,\hat{u}(x))\; \ud\nu (x) \\
    & \leq \int_{\tilde{\Omega}\setminus A} \liminf_{n\to\infty} B(x,\hat{u}_n(x))\; \ud\nu (x) \\
    & \leq \liminf_{n\to\infty} \int_{\tilde{\Omega}\setminus A} B(x,\hat{u}_n (x))\; \ud\nu (x) \\
    & = \liminf_{n\to\infty} \int_{\tilde{\Omega}} B(x,\hat{u}_n (x))\; \ud\nu (x) \\
    & = \liminf_{n\to\infty} {\mathcal B} (\hat{u}_n) ,
\end{align*}
which proves that $\mathcal B$ is lower semicontinuous. 
\end{proof}

\begin{theorem}\label{Sub}
Under the assumptions in this section, and if $j(\dom{\E_{B,\nu}})$ is dense in $L^2 (\Omega )$, the energy $\E_{B,\nu}^{L^2}$ is a Dirichlet form, that is, the semigroup $S^{B,\nu}$ on $L^2 (\Omega )$ is order preserving and $L^\infty$-contractive ($=$ submarkovian).
\end{theorem}

\begin{proof}
 We prove that for every $\hat{u}$, $\hat{v}\in W^{\Phi ,2} (\hat{\Omega} , \Omega )$ and for every normal contraction $p:\R\to\R$ 
 \begin{equation} \label{e.b.p}
 \E_{B,\nu} ( \hat{u} - p(\hat{u} - \hat{v})) + \E_{B,\nu} (\hat{v} + p(\hat{u} - \hat{v})) \leq \E_{B,\nu} (\hat{u}) + \E_{B,\nu} (\hat{v}) .
\end{equation}
A similar inequality was already proved for the energy $\E$ instead of $\E_{B,\nu}$, that is, when $B=0$ (see  \eqref{e.p}). Let $\hat{u}$, $\hat{v}\in W^{\Phi ,2} (\hat{\Omega} , \Omega )$ and let $p:\R\to\R$ be a normal contraction. Let 
\begin{equation*}
\lambda (x) := \begin{cases} \frac{p(\hat{u} (x) - \hat{v} (x))}{\hat{u} (x) - \hat{v} (x)} & \text{ if } \hat{u} (x) \not= \hat{v} (x) , \\ 0 & \text{ if } \hat{u} (x) = \hat{v} (x) . \end{cases} 
\end{equation*}
As $p$ is increasing and Lipschitz continuous with Lipschitz constant $\leq 1$, and as $p(0)=0$, then $\lambda (x)\in [0,1]$. Moreover, 
\begin{equation*} 
\hat{u} (x) - p (\hat{u} (x) - \hat{v} (x)) = \lambda (x) \hat{v} (x) + (1-\lambda (x)) \hat{u} (x) \text{ and } \hat{v} (x) +p(\hat{u} (x) - \hat{v} (x) )=\lambda (x) \hat{u} (x) +(1-\lambda (x)) \hat{v} (x) .
\end{equation*}
Hence, as a consequence of the inequality \eqref{e.p} and of the convexity of $B$ in the second variable,
\begin{align*}
    \E_{B,\nu} ( \hat{u} - p(\hat{u} - \hat{v})) + \E_{B,\nu} (\hat{v} + p(\hat{u} - \hat{v})) & = \E ( \hat{u} - p(\hat{u} - \hat{v})) + \E (\hat{v} + p(\hat{u} - \hat{v})) \\
    & \phantom{=\ } + \int_{\tilde{\Omega}} B(x,\lambda (x) \hat{v} (x) + (1-\lambda (x)) \hat{u} (x)) \; \ud\nu (x) \\
    &\phantom{=\ }+ \int_{\tilde{\Omega}} B(x,\lambda (x) \hat{u} (x) +(1-\lambda (x)) \hat{v} (x)) \; \ud\nu (x)\\
    & \leq \E (\hat{u} ) + \E (\hat{v}) + \int_{\tilde{\Omega}} ( \lambda (x) B(x,\hat{v} (x)) +  (1-\lambda (x)) B(x,\hat{u} (x))) \; \ud\nu (x) \\
    & \phantom{=\ } + \int_{\tilde{\Omega}} ( \lambda (x) B(x,\hat{u} (x)) +  (1-\lambda (x)) B(x,\hat{v} (x))) \; \ud\nu (x) \\
    & = \E_{B,\nu} (\hat{u}) + \E_{B,\nu} (\hat{v}) ,
\end{align*}
which proves \eqref{e.b.p}. Proceeding similarly as in the proof of Theorem \ref{theo-subm}, one deduces that the energy $\E_{B,\nu}^{L^2}$ satisfies the analogue of inequality \eqref{e.b.p} for every $u$, $v\in L^2 (\Omega )$ and every normal contraction $p:\R\to\R$, and this latter inequality implies that $\E_{B,\nu}^{L^2}$ is a Dirichlet form. 
\end{proof}

Let us identify the $j$-subgradient of $\E_{B,\nu}$ in the case when $B$ and its conjugate $B^*$ both satisfy the $\Delta_2$-condition (the integrable function appearing in the $\Delta_2$-condition has to be integrable with respect to the measure $\nu$). In this case, the Musielak-Orlicz space $L^B (\tilde{\Omega} , \nu )$ is a reflexive Banach space and $\mathcal B \hat{u} = \int_{\tilde{\Omega}} B(x , \hat{u} (x ))\;\ud\nu (x)$ (see also the proof of Theorem \ref{thm.eb.lsc}) is defined everywhere on $L^B (\tilde{\Omega}, \nu )$. 

\begin{lemma} \label{lem.identification.b.1}
Assume that $\Phi$ and $B$ are continuously partially differentiable in the last variables, that is, there exist functions $\varphi : \hat{\Omega} \times\hat{\Omega}\times\R\to\R$ and $\beta : \tilde{\Omega} \times \R \to\R$ such that the condition \eqref{cond.phi.derivative} (for the function $\varphi$) holds as well as the condition that 
\begin{equation}  \label{cond.B.derivative}
\begin{cases}
 \beta (\cdot , s) \text{ is measurable on } \tilde{\Omega} \text{ for every } s\in\R, \\
\beta (x,0 ) = 0 \text{ for every } x\in \tilde{\Omega } , \\
 \beta (x,\cdot ) \text{ is continuous, increasing and odd for every } x\in\tilde{\Omega} , \text{ and} \\
\displaystyle B (x,s) = \int_0^s \beta (x,\tau ) \; \ud\tau \text{ for every }  x\in\tilde{\Omega} , \, s\in\R . 
\end{cases}
\end{equation}
Then $\E_{B,\nu}$ is G\^ateaux differentiable on $W^{\Phi ,2} (\hat{\Omega} , \Omega ) \cap L^B (\tilde{\Omega} ,\nu )$, for every $\hat{u}$, $\hat{v}\in W^{\Phi ,2} (\hat{\Omega} ,\Omega )\cap L^B (\tilde{\Omega} , \nu )$, the functions
\begin{align*}
& \hat{\Omega}\times\hat{\Omega} \to \R, \quad \varphi (x,y,\hat{u}(x) - \hat{u}(y))(\hat{v} (x)-\hat{v}(y)) \quad \text{and} \\
& \tilde{\Omega} \to \R , \quad x\mapsto \beta (x,\hat{u}(x))\hat{v}(x) 
\end{align*}
are integrable, and 
\begin{equation} \label{eq.e.gateaux.B}
\E_{B,\nu}'(\hat{u}) \hat{v} = \frac12 \int_{\hat{\Omega}} \int_{\hat{\Omega}} \varphi (x,y,\hat{u} (x) - \hat{u} (y)) (\hat{v} (x) - \hat{v} (y)) \; \ud\mu (y) \; \ud\mu (x) + \int_{\tilde{\Omega}} \beta (x,\hat{u}(x)) \hat{v} (x) d\nu (x) .
\end{equation}
\end{lemma}

\begin{proof}
Proceed as in the proof of Lemma \ref{lem.identification}.
\end{proof}

Given $\hat{f}\in{\mathcal D}_a$, we say that $\hat{u}\in W^{\Phi ,2} (\hat{\Omega} , \Omega ) \cap L^B (\tilde{\Omega} ,\nu )$ is a {\em weak solution} of
\[
\text{P.V.}\int_{\hat{\Omega}} \varphi (x,y,\hat{u} (x) - \hat{u} (y)) \; \ud\mu (y) + \beta (x,\hat{u}(x)) \frac{\ud\nu}{\ud\mu} (x) = \hat{f} (x) \qquad (x\in\tilde{\Omega})
\]
if for every $\hat{v}\in{\mathcal D}$,
\[
\frac12 \int_{\hat{\Omega}} \int_{\hat{\Omega}} \varphi (x,y,\hat{u} (x) - \hat{u} (y)) (\hat{v} (x) - \hat{v} (y)) \; \ud\mu (y) \; \ud\mu (x) + \int_{\tilde{\Omega}} \beta (x,\hat{u}(x)) \hat{v} (x) \;\ud\nu (x) = \int_{\hat{\Omega}} \hat{f}(x) \hat{v} (x) \; \ud\mu (x) .
\]
The notation $\frac{\ud\nu}{\ud\mu}$ is of course only formal, and the meaning is made precise in the preceding line. However, if $\nu$ was absolutely continuous with respect to $\mu$, then $\frac{\ud\nu}{\ud\mu}$ would be the Radon-Nikodym derivative. 

\begin{theorem} \label{thm.identification.b.2}
Assume that $\Phi$ and $B$ are continuously partially differentiable in the last variables, that is, there exist functions $\varphi : \hat{\Omega} \times\hat{\Omega}\times\R\to\R$ and $\beta : \tilde{\Omega} \times \R \to\R$ such that the conditions \eqref{cond.phi.derivative} and \eqref{cond.B.derivative} hold. Assume that ${\mathcal D} \subseteq W^{\Phi ,2} (\hat{\Omega} , \Omega ) \cap L^B (\tilde{\Omega} ,\nu )$, and that ${\mathcal D}$ separates the points of the associate space ${\mathcal D}_a$. 

Then for every $u$, $f\in L^2 (\Omega )$ one has $(u,f)\in\partial_j\E_{B,\nu}$ if and only there exists $\hat{u}\in W^{\Phi ,2} (\hat{\Omega},\Omega ) \cap L^B (\tilde{\Omega},\nu)$ such that $\hat{u}|_{\Omega} = u$ and $\hat{u}$ is a weak solution of
\begin{equation} \label{eq.nonlocal.weak.converse.B}
\begin{cases}
    \displaystyle \text{P.V.}\int_{\hat{\Omega}} \varphi (x,y,\hat{u} (x) - \hat{u} (y)) \; \ud\mu (y) + \beta (x,\hat{u}(x)) \;\frac{\ud\nu}{\ud\mu} (x) = f(x) & \text{ for } x\in \Omega , \\
    \displaystyle \text{P.V.}\int_{\hat{\Omega}} \varphi (x,y,\hat{u} (x) - \hat{u} (y)) \; \ud\mu (y) + \beta (x,\hat{u} (x)) \;\frac{\ud\nu}{\ud\mu} (x) = 0 & \text{ for } x\in \tilde{\Omega}\setminus\Omega ,
\end{cases}
\end{equation}
and, moreover, for every $\hat{v}\in W^{\Phi ,2} (\hat{\Omega},\Omega ) \cap L^B (\tilde{\Omega},\nu )$,
\begin{equation} \label{eq.nonlocal.weak.converse.B.bc}
 \begin{split}
& \displaystyle \frac12 \int_{\hat{\Omega}} \int_{\hat{\Omega}} \varphi (x,y,\hat{u}(x)-\hat{u}(y)) (\hat{v} (x) -\hat{v} (y)) \; \ud\mu(y) \; \ud\mu(x) + \int_{\tilde{\Omega}} \beta (x,\hat{u}(x)) \hat{v} (x) \; \ud \nu (x)  \\
    & \phantom{\ud x} \displaystyle = \int_{\hat{\Omega}} P.V.\int_{\hat{\Omega}} \varphi (x,y,\hat{u} (x) - \hat{u} (y)) \; \ud\mu (y) \hat{v} (x) \; \ud\mu (x).
\end{split}
\end{equation}
Every such $\hat{u}$ is an elliptic extension of $u$ and satisfies 
\[
\hat{u} = \argmin_{j(\hat{w}) = u} (\E_{B,\nu} (\hat{w}) - \int_\Omega f \hat{w} \; \ud\mu ) = \argmin_{j(\hat{w}) = u} \E_{B,\nu} (\hat{w}) .
\]
The elliptic extension is unique when $\Phi (x,y,\cdot )$ is strictly convex for every $(x,y)\in {\mathbf S}_\Phi$.
\end{theorem}

\begin{proof}
    Proceed as in the proof of Theorem \ref{thm.identification.2}.
\end{proof}

\begin{remark}
Under the conditions of Theorem \ref{thm.identification.b.2} and when $\Phi (x,y,\cdot )$ is strictly convex for every $(x,y)\in {\mathbf S}_\Phi$, for every $(u,f)\in\partial_j\E_{B,\nu}$ there exists a unique elliptic extension $\hat{u}$ of $u$ such that \eqref{eq.nonlocal.weak.converse.B} and \eqref{eq.nonlocal.weak.converse.B.bc} hold. Unfortunately, there are in the general case no representation formulas for this elliptic extension.  
The situation is better in the case of quadratic energies (so that $\partial_j\E_{B,\nu}$ is a linear operator). 
More precisely, when $\displaystyle\Phi(x,y,s)=k(x,y) \frac{s^2}{2}$ and $\displaystyle B(x,s)=\beta(x) \frac{s^2}{2}$ for some functions $k:\hat{\Omega}\times\hat{\Omega}\to [0,\infty]$ and $\beta : \tilde{\Omega}\to [0,\infty [$, when $\nu = \mu$, and when 
\[
{\mathbf S}_\Phi = \{ (x,y)\in\hat{\Omega}\times\hat{\Omega} \st k(x,y) \not= 0\} \stackrel{!}{=} (\hat{\Omega}\times\Omega) \cup (\Omega\times\hat{\Omega}) ,
\]
letting 
\[
\rho(x):=\int_{\Omega}k(x,y) \;\ud y \qquad (x\in \hat\Omega\setminus\Omega),
\]
we have the following situation.
\begin{itemize}
    \item For Neumann exterior conditions, that is, $B=0$, the unique elliptic extension $\hat{u}_N$ is given by
    \[
    \hat{u}_N(x) =
    \begin{cases}
    \displaystyle u(x)\qquad &\;\quad\mbox{ if }  x\in \Omega,\\[2mm]
    \displaystyle\frac{1}{\rho(x)} \int_{\Omega} k(x,y) u(y)\;\ud\mu (y)  &\;\quad\mbox{ if } x\in \hat\Omega\setminus \Omega.
    \end{cases}
    \]
    \item For Robin exterior conditions (see Example \ref{ex.6.1} below), the unique elliptic extension $\hat{u}_R$ is given by
    \[
    \hat{u}_R(x)=
    \begin{cases}
    \displaystyle u(x)\qquad &\;\quad\mbox{ if }  x\in \Omega,\\[2mm]
    \displaystyle\frac{1}{\rho(x)+\beta(x)} \int_{\Omega} k(x,y) u(y)\;\ud\mu (y)   &\;\quad\mbox{ if } x\in \hat\Omega\setminus \Omega .
    \end{cases}
    \]
\end{itemize}
 These representation formulas follow easily from the second equation in \eqref{eq.nonlocal.weak.converse.B}; we refer to \cite{Cl21,ClWa20} for more details on these explicit elliptic extensions (even though these works consider only the case $\Omega\subseteq\hat{\Omega} = \R^N$). 
\end{remark}

\begin{example}[\bf Potentials and Robin complementary conditions] \label{ex.6.1}
In general, let us define the functions $B_{\Omega}$, $B_{\tilde{\Omega}\setminus\Omega} : \tilde{\Omega} \times\R \to [0,\infty ]$ by $B_{\Omega} (x,s):= B(x,s) \, 1_{\Omega} (x)$ and $B_{\tilde{\Omega}\setminus\Omega} (x,s) = B(x,s) \, 1_{\tilde{\Omega}\setminus\Omega} (x)$, so that $B=B_\Omega +B_{\tilde{\Omega}\setminus\Omega}$. Choose a Borel measure $\nu$ on $\tilde{\Omega}$ which is absolutely continuous with respect to the capacity $\Cap_\E$. The function $B_\Omega$ (or more precisely, the term $\hat{u}\mapsto \int_{\tilde{\Omega}} B_\Omega (x,\hat{u} (x)) \;\ud\nu (x)$ in the energy) can be seen as the restriction of $B$ to $\Omega$, and its contribution to the $j$-subgradient $\partial_j \E_{B,\nu}$ acts as a nonlinear, zero order perturbation of the $j$-subgradient $\partial_j \E$ in $L^2 (\Omega)$; see the first equality in \eqref{eq.nonlocal.weak.converse.B}. The function $B_{\tilde{\Omega}\setminus\Omega}$ (or more precisely, the term $\hat{u}\mapsto \int_{\tilde{\Omega}} B_{\tilde{\Omega}\setminus\Omega} (x,\hat{u} (x)) \;\ud\nu (x)$ in the energy) can be seen as the restriction of $B$ to $\tilde{\Omega}\setminus\Omega$, and its contribution to the $j$-subgradient $\partial_j\E_{B,\nu}$ acts as a nonlinear, zero order perturbation of the Neumann exterior condition (in $\hat{\Omega}\setminus\Omega$) and of the Neumann boundary condition (on $\tilde{\Omega}\setminus\hat{\Omega}$ or some other boundary); see the second equality in \eqref{eq.nonlocal.weak.converse.B} and \eqref{eq.nonlocal.weak.converse.B.bc}. We call the second equation in \eqref{eq.nonlocal.weak.converse.B} {\em Robin exterior condition} and equation \eqref{eq.nonlocal.weak.converse.B.bc} {\em Robin boundary condition}, and both together form the {\em Robin complementary condition}.

The smallest possible function $B$ satisfying the condition \eqref{cond.b}, namely the function $B=0$, yields the energy defined in Section \ref{sec.neumann}, no matter the choice of $\nu$, and its negative $j$-subgradient generates the semigroup $S^N$ from Section \ref{sec.neumann}. In the case $B=0$, the associated Robin exterior and boundary condition is just the Neumann exterior and boundary conditions. 
\end{example}

\begin{example}[\bf Dirichlet complementary condition] \label{ex.6.2}
In the context of the preceding example, the largest possible function $B$ satisfying the condition \eqref{cond.b} and acting only on the exterior of $\Omega$ (that is, $B_\Omega =0$!) is the function given by
\begin{equation} \label{eq.bd}
B_D (x,s) := \begin{cases}
    0 & \text{ if } s=0 \text{ or } x\in\Omega , \\
    \infty & \text{ otherwise.}
\end{cases}
\end{equation}
Choose now a Borel measure $\nu$ on $\tilde{\Omega}$ which is absolutely continuous with respect to the capacity, and consider the energy $\E_{B_D,\nu}$. The effective domain of this energy function is the space
\[ 
\mathring{W}^{\Phi ,2} (\hat{\Omega} , \Omega ) := \{ \hat{u}\in W^{\Phi ,2} (\hat{\Omega} , \Omega ) \st \hat{u} = 0 \; \nu\text{-almost everywhere in } \tilde{\Omega}\setminus\Omega \} .
\]
The complementary condition of the associated $j$-subgradient, namely the condition $\hat{u} = 0$ $\nu$-almost everywhere in $\tilde{\Omega}\setminus\Omega$, is called the {\em Dirichlet complementary condition}. Of course, this Dirichlet complementary condition depends on the choice of the measure $\nu$, and for appropriate choices of $\nu$, one may have the Dirichlet complementary condition on a part of $\tilde{\Omega}\setminus\Omega$, and the Neumann complementary condition on an other part. For example, when $\nu=\mu$, then the Dirichlet complementary condition means that $\hat{u}=0$ $\mu$-almost everywhere on $\hat{\Omega}\setminus\Omega$. This complementary condition becomes relevant when $\hat{\Omega}\setminus\Omega$ has nonzero measure, and we call it {\em Dirichlet exterior condition}. Another example is when $\nu$ is an appropriate "boundary measure" on $\tilde{\Omega}\setminus\hat{\Omega}$ or on the boundary of $\Omega$ (the measure of course has to be absolutely continuous with respect to the capacity $\Cap_\E$). Then we call the Dirichlet complementary condition {\em Dirichlet boundary condition}. Of course, one may have Dirichlet conditions on one part of the boundary and Neumann conditions on  other part of the boundary and on $\hat{\Omega}\setminus\Omega$.

The semigroup generated by the negative $j$-subgradient $-\partial_j \E_{B_D,\nu}$ is denoted by $S^{D,\nu}$. As a special case of $S^{B,\nu}$ (for $B=B_D$), if $j(\mathring{W}^{\Phi ,2} (\hat{\Omega} ,\Omega))$ is dense in $L^2 (\Omega )$, then $\E_{B_D,\nu}$ is a Dirichlet form and the semigroup $S^{D,\nu}$ is an order preserving, $L^\infty$-contractive ($=$ submarkovian) semigroup of contractions on $L^2 (\Omega )$.
\end{example}

\begin{example}[\bf Graphs] \label{ex.6.3}
    We take the situation of Example \ref{ex.3.2}, that is, $\hat{\Omega}$ is a countable set equipped with a weighted counting measure with weight $m:\hat{\Omega}\to ]0,\infty [$, $\Omega\subseteq\hat{\Omega}$ is a nonempty subset, $(b,0)$ is a locally finite graph on $\hat{\Omega}$, and every vortex $y\in\hat{\Omega}$ is connected to a vortex $x\in\Omega$ by a path of finite length. However, now we consider a full graph $(b,c)$ where in addition to the symmetric function $b:\hat{\Omega}\times\hat{\Omega}\to [0,\infty [$ one has a function $c:\hat{\Omega}\to [0,\infty [$ (see Keller, Lenz \& Wojciechowski \cite[Chapter 1, Definition 1.1]{KeLeWo21}); the case $c=0$ was considered in Example \ref{ex.3.2}. For $p\in ]1,\infty [$ one now considers the energy $\E_{(b,c)} : W^{(b,0)}_{p,2} (\hat{\Omega} , \Omega ) \to [0,\infty]$ given by,
    \[
    \E_{(b,c)} (\hat{u}) = \frac{1}{2p} \sum_{x\in\hat{\Omega}} \sum_{y\in\hat{\Omega}} b(x,y) |\hat{u} (x) - \hat{u} (y)|^p + \frac{1}{p} \sum_{x\in\hat{\Omega}} c(x) |u(x)|^p . 
    \]
    Let $c_\Omega = c\, 1_\Omega$ and $c_{\hat{\Omega}\setminus\Omega} := c\, 1_{\hat{\Omega}\setminus\Omega}$, so that $c = c_\Omega + c_{\hat{\Omega}\setminus\Omega}$. The contribution of $c_\Omega$ to the $j$-subgradient $\partial_j\E_{(b,c)}$ acts as a potential in $\Omega$ while the contribution of $c_{\hat{\Omega}\setminus\Omega}$ can be seen as a nonlinear perturbation of the Neumann exterior condition resulting in a Robin exterior condition.
\end{example}

\section{Domination} \label{sec-dom}

When we have two proper, lower semicontinuous, convex, densely defined functions $\E^{(1)}$, $\E^{(2)} : L^2 (\Omega ) \to ]-\infty ,\infty]$, and when we denote the semigroups generated by their negative subgradients by $S^{(1)}$ and $S^{(2)}$, respectively, then we say that $S^{(1)}$ is {\em dominated} by $S^{(2)}$ if, for every $t\geq 0$ and every $u_0\in L^2 (\Omega )$,
\[
| S^{(1)} (t) u_0 | \leq S^{(2)} (t) |u_0| ,
\]
and we say that $S^{(1)}$ is {\em totally dominated} by $S^{(2)}$ if  $S^{(1)}$ is dominated by both $S^{(2)}$ and the semigroup $\tilde{S}^{(2)}$ which is given by $\tilde{S}^{(2)} (t) u_0 := - S^{(2)} (t) (-u_0)$. 

Note that the semigroup $\tilde{S}^{(2)}$ is generated by the negative subgradient of the function $\tilde{\E}^{(2)} : L^2 (\Omega ) \to ]-\infty , \infty ]$ given by $\tilde{\E}^{(2)} (u) := \E^{(2)} (-u)$. Also domination and total domination are completely characterized in terms of the energies $\E^{(1)}$ and $\E^{(2)}$, at least when the dominating semigroup is order preserving. More precisely,  by Barth\'elemy  \cite[Th\'eor\`eme 3.3]{By96}, if the semigroup $S^{(2)}$ is order preserving, then the semigroup $S^{(1)}$ is dominated by the semigroup $S^{(2)}$ if and only if, for every $u$, $v\in L^2 (\Omega )$,
\begin{equation} \label{eq.char.domination}
 \E^{(1)} \left( \frac12 (|u| + |u|\wedge v)^+ \cdot\sgn(u) \right) + \E^{(2)} \left(\frac12 ( v + |u|\vee v)^+ \right) \leq \E^{(1)} (u)  + \E^{(2)} (v) . 
\end{equation}
The resulting characterization for total domination follows from \eqref{eq.char.domination} and by replacing in addition $\E^{(2)}$  by $\tilde{\E}^{(2)}$.

\begin{theorem}\label{theo-7.1}
 Let $(\hat{\Omega} ,{\mathcal B} (\hat{\Omega} ), \mu)$ be a $\sigma$-finite Hausdorff topological measure space, $\Phi : \hat{\Omega} \times \hat{\Omega} \times \R \to [0,\infty]$ satisfy the standard condition \eqref{cond.standard}, let $W^{\Phi ,2} (\hat{\Omega} ,\Omega)$ be the space defined in \eqref{eq.space.w}, let $\E : W^{\Phi ,2} (\hat{\Omega} ,\Omega )\to [0,\infty ]$ be the energy defined in \eqref{energy.e}, and let $j: W^{\Phi ,2} (\hat{\Omega} ,\Omega) \to L^2 (\Omega )$ be the restriction operator. Let $B^{(1)}$, $B^{(2)} : \tilde{\Omega}\times\R \to [0,\infty ]$ be two functions satisfying the condition \eqref{cond.b}, let $\nu^{(1)}$ and $\nu^{(2)}$ be two Borel measures on $\tilde{\Omega}$ which are absolutely continuous with respect to the capacity $\Cap_\E$, let $\E_{B^{(1)},\nu^{(1)}}$ and $\E_{B^{(2)},\nu^{(2)}}$ be the associated energies on $W^{\Phi ,2} (\hat{\Omega} , \Omega )$, and let $S^{B^{(1)},\nu^{(1)}}$ and $S^{B^{(2)},\nu^{(2)}}$ be the semigroups generated by their respective negative $j$-subgradients. 
If, for every measurable set $A\subseteq\tilde{\Omega}$,
 \begin{equation}
 \begin{cases}
\displaystyle \int_A B^{(2)} (x,|s|) \; \ud\nu^{(2)} (x) \leq \int_A B^{(1)} (x,s) \; \ud\nu^{(1)} (x) \text{ for every } s\in\R , \\
\displaystyle s\mapsto  \int_A B^{(1)} (x,s) \; \ud\nu^{(1)} (x) - \int_A B^{(2)} (x,s) \; \ud\nu^{(2)} (x) \text{ is increasing on } [0,\infty [, \text{ and} \\
\displaystyle s\mapsto \int_A B^{(1)} (x,-s) \; \ud\nu^{(1)} (x) - \int_A B^{(2)} (x,s) \; \ud\nu^{(2)} (x) \text{ is increasing on } [0,\infty [, 
\end{cases}
 \end{equation}
 then $S^{B^{(1)},\nu^{(1)}}$ is dominated by $S^{B^{(2)},\nu^{(2)}}$. If, in addition, 
  \begin{equation}
      \begin{cases}
\displaystyle \int_A B^{(2)} (x,\pm s) \; \ud\nu^{(2)} (x) \leq \int_A B^{(1)} (x,s) \; \ud\nu^{(1)} (x) \text{ for every } s\in\R , \\
\displaystyle s\mapsto  \int_A B^{(1)} (x,s) \; \ud\nu^{(1)} (x) - \int_A B^{(2)} (x,-s) \; \ud\nu^{(2)} (x) \text{ is increasing on } [0,\infty [, \text{ and} \\
\displaystyle s\mapsto \int_A B^{(1)} (x,-s) \; \ud\nu^{(1)} (x) - \int_A B^{(2)} (x,-s) \; \ud\nu^{(2)} (x) \text{ is increasing on } [0,\infty [, 
\end{cases}
 \end{equation}
 then $S^{B^{(1)},\nu^{(1)}}$ is totally dominated by $S^{B^{(2)},\nu^{(2)}}$. 
 
 In particular, if $B$ is any function satisfying the condition \eqref{cond.b}, if $\nu$ is any Borel measure on $\tilde{\Omega}$ which is absolutely continuous with respect to the capacity $\Cap_\E$, if $\E_{B,\nu}$ is the associated energy on $W^{\Phi ,2} (\hat{\Omega} , \Omega )$, and if $S^{B,\nu}$ is the semigroup generated by the negative $j$-subgradient of $\E_{B,\nu}$, then $S^{B,\nu}$ is totally dominated by $S^N$ (Neumann exterior conditions), and $S^{D,\nu}$ (Dirichlet exterior condition) is totally dominated by $S^{B,\nu}$. Moreover, if $\nu^{(1)}$ and $\nu^{(2)}$ are two Borel measures which are absolutely continuous with respect to the capacity $\Cap_\E$ and $\nu^{(2)} \leq \nu^{(1)}$, then $S^{D,\nu^{(1)}}$ is totally dominated by $S^{D,\nu^{(2)}}$.
\end{theorem}

The proof of the theorem uses the following result.

\begin{lemma} \label{lem.domination}
Let $(\hat{\Omega},{\mathcal A}, \mu )$ be a $\sigma$-finite measure space, and let $\Phi : \hat{\Omega}\times\hat{\Omega} \times\R \to [0,\infty]$ satisfy the condition \eqref{cond.phi} and the $\Delta_2$-condition. Then for every $x$, $y\in\hat\Omega$, and for every $u_0$, $u_1$, $v_0$, $v_1\in\R$,
\begin{multline}  \label{eq.proof.domination}
    \Phi (x,y, \frac12 (|u_0| + |u_0|\wedge v_0)^+ \sign u_0 - \frac12 (|u_1| + |u_1|\wedge v_1)^+ \sign u_1 ) + \Phi (x,y,|u_0| \vee v_0 - |u_1| \vee v_1 ) \\
    \leq \Phi (x,y,u_0 -u_1) + \Phi (x,y,v_0 -v_1) .
\end{multline}    
\end{lemma}

\begin{proof}
Fix $x$, $y\in\hat{\Omega}$. If $\Phi (x,y,s) = \infty$ for every $s\in\R$, $s\not= 0$, then the inequality \eqref{eq.proof.domination} is satisfied for every $u_0$, $u_1$, $v_0$, $v_1\in\R$. So assume that $\Phi$ is not equal to this function. As $\Phi$ satisfies the $\Delta_2$-condition, $\Phi (x,y,\cdot )$ takes only values in $\R$. As $\Phi (x,y,\cdot )$ is convex and takes only finite values, it is a continuous function. 

On the space $L^2 (\{ x,y\} )$ ($=\R^2$) we consider the energy $\E : L^2 (\{ x,y\}) \to\R$ given by $\E (u) := \Phi (x,y,u(x) - u(y))$. This energy is convex and continuous. By Theorem \ref{theo-subm}, it is in addition a Dirichlet form. The semigroup $S$ generated by $-\partial\E$ is order preserving (and $L^\infty$-contractive) on $L^2 (\{ x,y\})$. As a consequence, $S$ totally dominates itself. Indeed, for every $u\in L^2 (\{ x,y\})$ one has $-|u|\leq u \leq |u|$ and $-|u|\leq -u \leq |u|$, and as $S$ is order preserving, this implies for every $t\geq 0$,
\[
S(t) (-|u|) \leq S(t) u \leq S(t) |u| \text{ and } S(t) (-|u|) \leq S(t) (-u) \leq S(t) |u| .
\]
Now the fact that $\Phi(x,y,\cdot )$ is an even function implies that $\E$ is even, and therefore $S(t) (-u) = - S(t)u$. Altogether this gives the total domination. Now the characterization of total domination gives the inequality \eqref{eq.proof.domination} for the fixed $x$, $y\in\hat{\Omega}$ and every $u_0$, $u_1$, $v_0$, $v_1\in\R$. As $x$ and $y$ were arbitrary, this proves the claim.
\end{proof}

\begin{proof}[\bf Proof of Theorem \ref{theo-7.1}]
We show that for every $\hat{u}$, $\hat{v}\in W^{\Phi ,2} (\hat{\Omega}, \Omega )$,
\begin{equation} \label{e1.e2}
   \E_{B^{(1)},\nu^{(1)}} \left( \frac12 (|\hat{u}| + |\hat{u}|\wedge \hat{v})^+ \sgn \hat{u} \right) + \E_{B^{(2)},\nu^{(2)}} \left( \frac12 (\hat{v} + |\hat{u}| \vee \hat{v})^+ \right) \leq \E_{B^{(1)},\nu^{(1)}} (\hat{u}) + \E_{B^{(2)},\nu^{(2)}} (\hat{v}) . 
\end{equation}
Observe that the semigroup $S^N$ totally dominates itself because $0$ is an equilibrium point of $S^N$, $S^N$ is order preserving (by Theorem \ref{theo-subm}), and $S^N (t) (-u_0) = -S^N (t) u_0$ (the energy $\E$ is even). 

By Lemma \ref{lem.domination}, for every $\hat{u}$, $\hat{v}\in W^{\Phi ,2} (\hat{\Omega}, \Omega )$ and for every $x$, $y\in\hat{\Omega}$,
\begin{multline*} 
    \Phi \left(x,y, \frac12 (|\hat{u} (x)| + |\hat{u} (x)|\wedge \hat{v} (x))^+ \sign \hat{u} (x) - \frac12 (|\hat{u} (y)| + |\hat{u} (y)|\wedge \hat{v} (y))^+ \sign \hat{u} (y)\right ) \\
    + \Phi \left(x,y, \frac12 (\hat{v} (x) + |\hat{u} (x)| \vee \hat{v}(x))^+ - \frac12 (\hat{v} (y) + |\hat{u} (y)| \vee \hat{v} (y))^+\right ) \\
    \leq \Phi (x,y, \hat{u}(x) -\hat{u}(y)) + \Phi (x,y, \hat{v}(x) -\hat{v}(y)) .
\end{multline*}
Integrating this inequality over $(x,y)\in \hat{\Omega}\times\hat{\Omega}$ with respect to the product measure $\mu\otimes\mu$ implies that for every $\hat{u}$, $\hat{v}\in W^{\Phi ,2} (\hat{\Omega}, \Omega )$,
\begin{equation} \label{e.dominates.e}
   \E \left( \frac12 (|\hat{u}| + |\hat{u}|\wedge \hat{v})^+ \sgn \hat{u} \right) + \E \left( \frac12 (\hat{v} + |\hat{u}| \vee \hat{v})^+ \right) \leq \E (\hat{u}) + \E (\hat{v}) . 
\end{equation}
Next, for every $\hat{u}$, $\hat{v}\in W^{\Phi ,2} (\hat{\Omega}, \Omega )$,
\begin{align*}
  & \int_{\tilde{\Omega}} B^{(1)} \left(x, \frac12 (|\hat{u} (x)| + |\hat{u} (x) |\wedge \hat{v} (x))^+ \sgn \hat{u} (x)\right ) \; \ud\nu^{(1)} (x) + \int_{\tilde{\Omega}} B^{(2)} \left(x, \frac12 (\hat{v} (x) + |\hat{u} (x) |\vee \hat{v} (x))^+ \right) \; \ud\nu^{(2)} (x)  \\
   & =: I_1 + I_2 + I_3 + I_4 ,
\end{align*}
where for $k=1$, $2$, $3$, $4$
\[
I_k = \int_{\tilde{\Omega}_k} B^{(1)} \left(x, \frac12 (|\hat{u} (x)| + |\hat{u} (x) |\wedge \hat{v} (x))^+ \sgn \hat{u} (x) \right) \; \ud\nu^{(1)} (x) + \int_{\tilde{\Omega}_k} B^{(2)} \left(x, \frac12 (\hat{v} (x) + |\hat{u} (x) |\vee \hat{v} (x))^+ \right) \; \ud\nu^{(2)} (x) 
\]
and
\begin{align*}
    \tilde{\Omega}_1 & := \{ x\in\tilde{\Omega} \st |\hat{u} (x)| \leq \hat{v} (x) \} , \\
    \tilde{\Omega}_2 & := \{ x\in\tilde{\Omega} \st -|\hat{u} (x)| \geq \hat{v} (x) \text{ and } (\hat{u} (x) , \hat{v} (x)) \not= (0,0) \} , \\
    \tilde{\Omega}_3 & := \{ x\in\tilde{\Omega} \st \hat{u} (x) > |\hat{v} (x)| \} , \text{ and}\\
    \tilde{\Omega}_4 & := \{ x\in\tilde{\Omega} \st \hat{u} (x) < - |\hat{v} (x)| \} .
\end{align*}
Note that $\tilde{\Omega}$ is the disjoint union of the $\tilde{\Omega}_k$. For $x\in\tilde{\Omega}_1$ one has $\frac12 (|\hat{u} (x)| + |\hat{u} (x) |\wedge \hat{v} (x))^+ \sgn \hat{u} (x) = \hat{u} (x)$ and $\frac12 (\hat{v} (x) + |\hat{u} (x) |\vee \hat{v} (x))^+ = \hat{v} (x)$, so that 
\[
I_1 = \int_{\tilde{\Omega}_1} B^{(1)} (x,\hat{u} (x)) \; \ud\nu^{(1)} (x) + \int_{\tilde{\Omega}_1} B^{(2)} (x,\hat{v} (x)) \; \ud\nu^{(2)} (x) .
\]
For $x\in\tilde{\Omega}_2$ one has $\frac12 (|\hat{u} (x)| + |\hat{u} (x) |\wedge \hat{v} (x))^+ \sgn \hat{u} (x) = 0$ and $\frac12 (\hat{v} (x) + |\hat{u} (x) |\vee \hat{v} (x))^+ = 0$. Since $B^{(1)} (x,0) = B^{(2)} (x,0)=0$, and since $B^{(1)}$ and $B^{(2)}$ are positive, we have
\[
I_2 = 0 \leq \int_{\tilde{\Omega}_2} B^{(1)} (x,\hat{u} (x)) \; \ud\nu^{(1)} (x) + \int_{\tilde{\Omega}_2} B^{(2)} (x,\hat{v} (x)) \; \ud\nu^{(2)} (x) .
\]
For $x\in\tilde{\Omega}_3$ one has $\frac12 (|\hat{u} (x)| + |\hat{u} (x) |\wedge \hat{v} (x))^+ \sgn \hat{u} (x) = \frac12 (\hat{u} (x) + \hat{v} (x))$ and $\frac12 (\hat{v} (x) + |\hat{u} (x) |\vee \hat{v} (x))^+ = \frac12 (\hat{u} (x) + \hat{v} (x))$. By convexity of $B^{(1)}$ and $B^{(2)}$, we have that
\begin{align*}
I_3 & \leq \frac12 \int_{\tilde{\Omega}_3} B^{(1)} (x,\hat{u} (x)) \; \ud\nu^{(1)} (x) + \frac12 \int_{\tilde{\Omega}_3} B^{(1)} (x,\hat{v} (x)) \; \ud\nu^{(1)} (x) \\
& \phantom{\leq \ } + \frac12 \int_{\tilde{\Omega}_3} B^{(2)} (x,\hat{u} (x)) \; \ud\nu^{(2)} (x) + \frac12 \int_{\tilde{\Omega}_3} B^{(2)} (x,\hat{v} (x)) \; \ud\nu^{(2)} (x) \\
& \leq \frac12 \int_{\tilde{\Omega}_3} B^{(1)} (x,\hat{u} (x)) \; \ud\nu^{(1)} (x) + \frac12 \int_{\tilde{\Omega}_3} B^{(1)} (x,\hat{u} (x)) \; \ud\nu^{(1)} (x) \\
& \phantom{\leq \ } + \frac12 \int_{\tilde{\Omega}_3} B^{(2)} (x,\hat{v} (x)) \; \ud\nu^{(2)} (x) + \frac12 \int_{\tilde{\Omega}_3} B^{(2)} (x,\hat{v} (x)) \; \ud\nu^{(2)} (x) \\
& = \int_{\tilde{\Omega}_3} B^{(1)} (x,\hat{u} (x)) \; \ud\nu^{(1)} (x) + \int_{\tilde{\Omega}_3} B^{(2)} (x,\hat{v} (x)) \; \ud\nu^{(2)} (x) .
\end{align*}
For the second estimate we have used the assumption.

In a similar way one proves that
\[
I_4 \leq \int_{\tilde{\Omega}_4} B^{(1)} (x,\hat{u} (x)) \; \ud\nu^{(1)} (x) + \int_{\tilde{\Omega}_4} B^{(2)} (x,\hat{v} (x)) \; \ud\nu^{(2)} (x) .
\]
Taking the estimates for the $I_k$ together, it follows that 
\begin{align*}
  & \int_{\tilde{\Omega}} B^{(1)} (x, \frac12 (|\hat{u} (x)| + |\hat{u} (x) |\wedge \hat{v} (x))^+ \sgn \hat{u} (x) ) \; \ud\nu^{(1)} (x) \\
   &  + \int_{\tilde{\Omega}} B^{(2)} (x, \frac12 (\hat{v} (x) + |\hat{u} (x) |\vee \hat{v} (x))^+ ) \; \ud\nu^{(2)} (x)  \\
   &  \leq \int_{\tilde{\Omega}} B^{(1)} (x,\hat{u} (x)) \; \ud\nu^{(1)} (x) + \int_{\tilde{\Omega}} B^{(2)} (x,\hat{v} (x)) \; \ud\nu^{(2)} (x) .
\end{align*}
This estimate together with  \eqref{e.dominates.e} yields our first claim \eqref{e1.e2}. Using again that $j$ is simply the restriction operator, and in particular that $j( \frac12 (|\hat{u}| + |\hat{u}|\wedge \hat{v})^+ \sgn \hat{u} ) = \frac12 (|j(\hat{u})| + |j(\hat{u})|\wedge j(\hat{v}))^+ \sgn j(\hat{u})$ and $j (\frac12 (\hat{v} + |\hat{u}|\vee \hat{v} )^+) = \frac12 (j(\hat{v}) + |j(\hat{u})|\vee j(\hat{v}))^+$, it follows from \eqref{e1.e2} and the representation \eqref{energy.e} that for every $u$, $v\in L^2 (\Omega )$,
\begin{equation} \label{e1.e2.l2}
   \E_{B^{(1)},\nu^{(1)}}^{L^2} \left( \frac12 (|u| + |u|\wedge v)^+ \sgn u \right) + \E_{B^{(2)},\nu^{(2)}}^{L^2} \left( \frac12 (v + |u| \vee v)^+ \right) \leq \E_{B^{(1)},\nu^{(1)}}^{L^2} (u) + \E_{B^{(2)},\nu^{(2)}}^{L^2} (v) . 
\end{equation}
It follows from here and from the characterization \eqref{eq.char.domination} that the semigroup $S^{B^{(1)},\nu^{(1)}}$ is dominated by the semigroup $S^{B^{(2)},\nu^{(2)}}$. 

Now, if $B : \tilde{\Omega} \times\R\to [0,\infty ]$ is any function satisfying the condition \eqref{cond.b}, and if $\nu$ is any Borel measure on $\tilde{\Omega}$ which is absolutely continuous with respect to the capacity $\Cap_\E$, then the semigroup $S^{B,\nu}$ generated by the negative $j$-subgradient of the energy $\E_{B,\nu}$ is totally dominated by the semigroup $S^N$ (apply the first part of the statement with $B^{(1)} = B$, $\nu^{(1)} = \nu$, $B^{(2)} = 0$ and $\nu^{(2)} = \mu$), and the semigroup $S^{D,\nu}$ is totally dominated by $S^{B,\nu}$ (apply the first part of the statement with $B^{(1)} = B_D$, $\nu^{(1)} = \nu$, $B^{(2)} = B$ and $\nu^{(2)} = \nu$). 

The final statement about the two semigroups coming from Dirichlet conditions follow from the first part of the statement applied with $B^{(1)} = B^{(2)} = B_D$ and the measures $\nu^{(1)}$ and $\nu^{(2)}$.
\end{proof}

\begin{remark}\label{rem-to-dom}
In the situation of the second part of Theorem \ref{theo-7.1} we write
\begin{equation*}
S^{D,\nu}\preceq S^{B,\nu}\preceq S^N
\end{equation*}
to say that the semigroup $S^{B,\nu}$ is sandwiched between the semigroups $S^{D,\nu}$ and $S^N$, that is, $S^{B,\nu}$ is dominated by $S^N$, and $S^{B,\nu}$ dominates $S^{D,\nu}$.
This property  for linear semigroups generated by the Laplace operators with Dirichlet boundary conditions, Neumann boundary conditions and Robin boundary conditions has been previously studied in Arendt \& Warma \cite{ArWa2003}. The case of nonlinear semigroups generated by the $p$-Laplace operator ($p\in ]1,\infty[$) has been investigated in Chill \& Warma \cite{ChWa-1,ChWa19}. For more general situations of nonlinear semigroups generated by Dirichlet forms, see Claus \cite{Cl21} and Chill \& Claus \cite{ChCl25}.
\end{remark}

\section{Ultracontractivity for Dirichlet and Neumann exterior conditions} \label{sec-UlDN}

Before we give the main concern of this section, we recall that the ultracontractivity properties of the $p$-Laplace operator $p\in [2,\infty[$ with Dirichlet and/or Neumann boundary conditions have been intensively studied in the works by Cipriani \cite{Cip1994,Cip1994-2}, Cipriani \& Grillo \cite{CipGri}, Coulhon \& Hauer \cite{CoHa2016} and their references. For Wentzell (dynamic) boundary conditions we refer to the paper by Warma \cite{War-Ultra} and the references therein.

In this section, we take the situation of Examples \ref{ex.3.1} and \ref{ex.6.2}, that is, $\Omega\subseteq\hat{\Omega}\subseteq\R^N$ ($N\ge 1$) are two open, nonempty sets equipped with the Lebesgue measure. We assume in this section that $\Omega$ is bounded. Let $p\in ]1,\infty [$, and the function $\Phi :\hat{\Omega}\times\hat{\Omega}\times\R \to [0,\infty]$ be given by
  \begin{equation}\label{Phi}
 \Phi (x,y,s) = \frac{1}{p} k(x,y) |s|^p 1_{{\mathbf S}} (x,y)
 \end{equation}
 where the kernel $k:\hat{\Omega}\times\hat{\Omega} \to [0,\infty]$ satisfies the condition \eqref{cond.k} for some constants $C_1$, $C_2 >0$ and $\theta\in ]0,1[$, and ${\mathbf S} \subseteq\hat{\Omega}\times\hat{\Omega}$ satisfies the thickness condition \eqref{cond.connected}. 

Throughout the remainder of the paper, we simply write $\|\cdot\|_q:=\|\cdot\|_{L^q(\Omega)}$ for the norm in $L^q(\Omega)$ ($q\in [1,\infty ]$). Moreover, for every $p$, $r\in [1,\infty[$ and $\theta\in ]0,1[$ we let
\begin{align*}\label{Wsp}
W^{\theta,{\mathbf S}}_{p,r}(\hat{\Omega},\Omega ) := \left\{ \hat{u}\in L^0(\hat\Omega) \st \hat{u}|_\Omega \in L^r (\Omega ) \text{ and } \iint_{{\mathbf S}} \frac{|\hat{u}(x)-\hat{u}(y)|^p}{|x-y|^{N+\theta p}}\;\ud y\;
\ud x<\infty\right\}
\end{align*}
which endowed with the norm
\begin{equation}\label{norm-N}
\| \hat{u}\|_{W^{\theta,{\mathbf S}}_{p,r}(\hat{\Omega},\Omega)} := \|\hat{u}|_\Omega\|_{r} +\left(\iint_{{\mathbf S}} \frac{|\hat{u}(x)-\hat{u}(y)|^p}{|x-y|^{N+\theta p}}\;dy\;dx\right)^{1/p},
\end{equation}
is a Banach space (proceed as in Lemma \ref{lem.identification}). 

When  $\hat{\Omega} = \Omega$, we simply denote $W^{\theta,{\mathbf S}}_{p,p}(\hat{\Omega},\Omega )$ by $W^{\theta}_{p}(\Omega)$, and $\mathring{W}^{\theta,{\mathbf S}}_{p,p}(\hat{\Omega},\Omega)$ by $\mathring{W}^{\theta}_{p}(\Omega)$. In this case  ${\mathbf S}=\Omega\times\Omega$.

\subsection{Dirichlet exterior and boundary conditions}\label{sec-UlD}
 
We start with the situation of Example \ref{ex.6.2}. We consider Dirichlet exterior conditions, that is, we take the function $B_D$ from \eqref{eq.bd},
and we let $\nu$ be the sum of the Lebesgue measure on $\hat{\Omega}\setminus\Omega$ and the restriction of the $(N-1)$-dimensional Hausdorff measure $\mathcal H_{N-1}$ to the boundary $\partial\hat{\Omega}$. We assume that $\hat{\Omega}$ has Lipschitz boundary, so that the $(N-1)$-dimensional Hausdorff measure is absolutely continuous with respect to the capacity $\Cap_\E$. Consider the energy $\E_{B_D,\nu} : W^{\theta,{\mathbf S}}_{p,2} (\hat{\Omega}, \Omega ) \to [0,\infty]$ given by 
 \[
 \E_{B_D,\nu} (\hat{u} ) = \frac{1}{2p} \iint_{{\mathbf S}} k(x,y) \, |\hat{u} (x) - \hat{u} (y)|^p \; \ud y \; \ud x + \int_{\hat{\Omega}\setminus\Omega} B_D (x,\hat{u} (x))\; \ud x +\int_{\tilde{\Omega}\setminus\hat{\Omega}}B_D(x,\hat{u}(x))\;d\mathcal H_{N-1}(x),
 \]
so that the effective domain is the Sobolev space $\mathring{W}^{\theta,{\mathbf S}}_{p,2} (\hat{\Omega} , \Omega )$. Here,
\begin{equation*}\label{Wsp0}
\mathring{W}^{\theta,{\mathbf S}}_{p,r}(\hat{\Omega},\Omega) := \left\{u\in W_{p,r}^{\theta,{\mathbf S}}(\hat{\Omega},\Omega) \st u=0 \;\nu\text{-almost everywhere in } \tilde{\Omega}\setminus\Omega\right\}
\end{equation*}
which is a closed subspace of ${W_{p,r}^{\theta,{\mathbf S}}(\hat{\Omega},\Omega)}$. If $r=p$, then $\mathring{W}^{\theta,{\mathbf S}}_{p,p}(\hat{\Omega},\Omega)$ endowed with the norm
\begin{align}\label{norm}
\|u\|_{\mathring{W}^{\theta}_{p,p}(\hat{\Omega},\Omega)}=\left(\iint_{{\mathbf S}}\frac{|u(x)-u(y)|^p}{|x-y|^{N+\theta p}}\;\ud x\;\ud y\right)^{1/p}
\end{align}
is a Banach space (${\Omega}$ is bounded!) and the norm in \eqref{norm} is equivalent to the one given in \eqref{norm-N} with $r=p$. 

As a consequence,
  \[
 \E_{B_D,\nu} (\hat{u}) = \begin{cases}
     \displaystyle \frac{1}{2p} \iint_{{\mathbf S}} k(x,y) \, |\hat{u} (x) - \hat{u} (y)|^p \; \ud y \; \ud x & \text{if } \hat{u}\in \mathring{W}^{\theta,{\mathbf S}}_{p,2} (\hat{\Omega} , \Omega ) , \\[2mm]
     \infty & \text{otherwise}.
     \end{cases}
 \]

We let $j$ be the restriction operator from $\mathring{W}^{\theta,{\mathbf S}}_{p,2} (\hat{\Omega} , \Omega )$ into $L^2 (\Omega )$, and we denote by $S^{D,\nu}$ the semigroup on $L^2 (\Omega )$ generated by the negative $j$-subgradient of $\E_{B_D,\nu}$ (Theorem \ref{Sub}). It follows from the G\^ateaux differentiability of $\mathcal E_{B_D}$ on the space $\mathring{W}^{\theta,{\mathbf S}}_{p,2} (\hat{\Omega} ,\Omega )$ and the definition of the $j$-subgradient $\partial_j\E_{B_D,\nu}$ that for every $(u,f)\in\partial_j\E_{B_D,\nu}$ there exists an elliptic extension $\hat{u}\in \mathring{W}^{\theta,{\mathbf S}}_{p,2} (\hat{\Omega} ,\Omega )$ of $u$, such that for every $\hat{v}\in \mathring{W}^{\theta,{\mathbf S}}_{p,2}(\hat{\Omega},\Omega)$ the following equality holds:
\begin{align}\label{INt-parts}
\int_{\Omega} f \hat{v} \;\ud x = \iint_{{\mathbf S}}k(x,y)|\hat{u}(x)-\hat{u}(y)|^{p-2}(\hat{u}(x)-\hat{u}(y))(\hat{v}(x)-\hat{v} (y)\; \ud y \; \ud x.
\end{align}
Also remark that the restriction operator $j$ maps $W_{p,p}^{\theta,{\mathbf S}}(\hat{\Omega},\Omega)$ into $W_p^{\theta}(\Omega)$.\\

The aim of this subsection is to prove an $L^q-L^\infty$-H\"older continuity property of the semigroup $S^{D,\nu}$ generated by the negative $j$-subgradient $-\partial_j\mathcal E_{B_D,\nu}$. Throughout the remainder of this subsection, without any mention we assume that $\theta\in ]0,1[$ and $p\in [2,\infty[$.
More precisely, we have the following theorem which is the main result of this subsection.

\begin{theorem}\label{theo-ultra}
 Take the assumptions of this subsection, and let $S^{D,\nu}$ be the semigroup generated by the negative $j$-subgradient $-\partial_j\mathcal E_{B_D,\nu}$.
 Then, for every $q\in [1,\infty[$ there is a constant $C=C(N,p,q,\theta )>0$ such that for every $t>0$ and every $u_0$, $v_0\in L^2(\Omega)\cap L^q(\Omega)$ one has the estimates
 \begin{equation}\label{est-65-3}
 \|S^{D,\nu}(t)u_0-S^{D,\nu}(t)v_0\|_\infty \le C\frac{|\Omega|^\alpha}{t^{\beta}}\|u_0-v_0\|_q^\gamma,
 \end{equation}
where
 \begin{equation}\label{abg2}
 \begin{cases}
\alpha= \alpha(p,q):= \frac{N-\theta p}{N}\Big[1-\left(\frac{q}{q+p-2}\right)^{\frac{N}{\theta p}}\Big],\\
\beta=\beta(p,q):=\frac{1}{p-2}\Big[1-\left(\frac{q}{q+p-2}\right)^{\frac{N}{\theta p}}\Big],\\
 \gamma=\gamma(p,q):=\left(\frac{q}{q+p-2}\right)^{\frac{N}{\theta p}},
 \end{cases}
 \end{equation}
 if $2\le p<\frac{N}{\theta}$, and
 \begin{equation}\label{abg2-b}
 \begin{cases}
\alpha= \alpha(p,q):= \frac{\theta p-N}{N}\Big[1-\left(\frac{q}{q+p-2}\right)^{\frac{\theta p}{N}}\Big],\\
\beta=\beta(p,q):=\frac{1}{p-2}\Big[1-\left(\frac{q}{q+p-2}\right)^{\frac{\theta p}{N}}\Big],\\
 \gamma=\gamma(p,q):=\left(\frac{q}{q+p-2}\right)^{\frac{\theta p}{N}},
 \end{cases}
 \end{equation}
 if $2<\frac{N}{\theta}\le p<\infty$. 
 \end{theorem}

\begin{remark}
{\em In Theorem \ref{theo-ultra} we did not include the case $q=\infty$, since in that case $(\alpha,\beta,\gamma)=(0,0,1)$, and \eqref{est-65-3} is exactly the $L^\infty$-contractivity of the semigroup. Also observe that if $\frac{N}{\theta}= p$, then $ \displaystyle(\alpha,\beta,\gamma)=\left(0, 1,\frac{q}{q+p-2}\right)$.
}
\end{remark}

In order to prove the theorem, we need some preparations. This will be done in a series of lemmas. 
We start with a logarithmic Sobolev inequality, where we use
that there is a constant $C=C(N,\theta,p,\Omega)>0$ such that for every $\hat{u}\in \mathring{W}^{\theta,{\mathbf S}}_{p,p}(\hat{\Omega},\Omega)$ the estimate
\begin{equation}\label{sob-emb}
   \|\hat{u}\|_{q}:=  \|\hat{u}|_{\Omega}\|_{q}\le C\|\hat{u}\|_{\mathring{W}^{\theta,{\mathbf S}}_{p,p}(\hat{\Omega},\Omega)}
\end{equation}
holds (that is, we have the continuous inclusion $\mathring{W}^{\theta,{\mathbf S}}_{p,p}(\hat{\Omega},\Omega)\subset L^q(\Omega)$), with
\begin{equation}\label{eq-q}
q=
\begin{cases}
\frac{Np}{N-\theta p}  \qquad &\mbox{ if } N>\theta p,\\
\infty &\mbox{ if } N<\theta p,\\
q<\infty &\mbox{ if } N=\theta p.
\end{cases}    
\end{equation}

\begin{lemma}[\bf Logarithmic Sobolev inequality for $\mathring{W}^{\theta,{\mathbf S}}_{p,p}(\hat{\Omega},\Omega)$]\label{LSI}
For every $\varepsilon>0$ and $\hat{u}\in \mathring{W}^{\theta,{\mathbf S}}_{p,p}(\hat{\Omega},\Omega)$, the logarithmic Sobolev inequality
\begin{align}\label{LSI-2-1}
  \int_\Omega|\hat{u}|^p\log|\hat{u}|^p\;\ud x&-\left(\int_{\Omega}|\hat{u}|^p\;\ud x\right)\log\left(\int_\Omega|\hat{u}|^p\;\ud x\right)\notag\\
  \le &C(N,p,\theta)\left[-\|\hat{u}\|_p^p\log(\varepsilon)+\varepsilon C\iint_{{\mathbf S}}k(x,y)|\hat{u}(x)-\hat{u}(y)|^p\;\ud y\;\ud x  \right]
\end{align}
holds true, with
\begin{equation}\label{C-SLI}
C(N,p,\theta)=
\begin{cases}
\frac{N}{\theta p^2}=\frac{N}{\theta p}\frac{1}{p}\;&\mbox{ if } 2\le  p<\frac{N}{\theta}, \\
\frac{\theta}{N}=\frac{\theta p}{N}\frac{1}{p}&\mbox{ if } \frac{N}{\theta}\le  p<\infty ,
\end{cases}
\end{equation}
and where $C\geq 0$ in \eqref{LSI-2-1} is the constant from the inequality \eqref{sob-emb}.
\end{lemma}

\begin{proof}
First, observe that, since $p\ge 2$, the continuous inclusion $\mathring{W}^{\theta,{\mathbf S}}_{p,p}(\hat{\Omega},\Omega)\subseteq \mathring{W}^{\theta,{\mathbf S}}_{p,2}(\hat{\Omega},\Omega)$ holds.

Second, using the homogeneity and since
\[ 
\iint_{{\mathbf S}}\frac{||\hat{u}|(x)-|\hat{u}|(y)|^p}{|x-y|^{N+\theta p}}\;\ud y\;\ud x\le \iint_{{\mathbf S}}\frac{|\hat{u}(x)-\hat{u}(y)|^p}{|x-y|^{N+\theta p}}\;\ud y\;\ud x
\]
for all $\hat{u}\in \mathring{W}^{\theta,{\mathbf S}}_{p,p}(\hat{\Omega},\Omega)$, it suffices to prove the lemma for non-negative functions $\hat{u}\in \mathring{W}^{\theta,{\mathbf S}}_{p,p}(\hat{\Omega},\Omega)$ such that $\|\hat{u}\|_{p}=1$. For this, let $\hat{u}\in \mathring{W}^{\theta,{\mathbf S}}_{p,p}(\hat{\Omega},\Omega)$ be non-negative with $\|\hat{u}\|_{p}=1$, so that the weighted Lebesgue measure $\hat{u}(x)^p\;\ud x$ is a probability measure. 

\begin{itemize}
\item {\bf The case $2\le  p<\frac N\theta$}:
let $p_\theta:=\frac{Np}{N-\theta p}$ and $q:=\frac{\theta p^2}{N-\theta p}$ so that $p+q=p_\theta$.
Then, for every $\varepsilon>0$,
\begin{align*}
    \int_\Omega\log(\hat{u})\hat{u}^p\;\ud x =&  \frac{1}{q}\int_\Omega\log(\hat{u}^q)\hat{u}^p\;\ud x\\
    \le & \frac 1q\log\int_{\Omega}\hat{u}^{p+q}\;\ud x =\frac 1q\log\int_{\Omega}\hat{u}^{p_\theta}\;\ud x\\
    \le & \frac{N}{\theta p^2}\log\|\hat{u}\|_{\frac{Np}{N-\theta p}}^p\\
    \le& \frac{N}{\theta p^2}\left(-\log(\varepsilon)+\var \|\hat{u}\|_{\frac{Np}{N-\theta p}}^p\right)\\
    \le &\frac{N}{\theta p^2}\left(-\log(\var) +\varepsilon C\iint_{{\mathbf S}}\frac{|\hat{u}(x)-\hat{u}(y)|^p}{|x-y|^{N+\theta p}}\;\ud y\;\ud x\right)\\
     \le &\frac{N}{\theta p^2}\left(-\log(\var) +\varepsilon C\iint_{{\mathbf S}}k(x,y)|\hat{u}(x)-\hat{u}(y)|^p\;\ud y\;\ud x\right) ,
\end{align*}
where we used Jensen's inequality (since the logarithmic function is concave), the fact that $\log(t)<t$,  the Sobolev embedding \eqref{sob-emb} and the condition \eqref{cond.k}. We have shown \eqref{LSI-2-1} with $C(N,p,\theta)$ given in \eqref{C-SLI}.

\item {\bf The case $\frac N\theta<  p<\infty$}: in this case we know that for every $r\in [1,\infty]$ there is a constant $C>0$ such that
\begin{equation}\label{SI}
\|\hat{u}\|_{r}\le C\|\hat{u}\|_{\mathring{W}^{\theta,{\mathbf S}}_{p,p}(\hat{\Omega},\Omega)} . 
\end{equation}
Taking $r=p_\theta=\frac{\theta p^2}{\theta p-N}$ in \eqref{SI},  $q:=\frac{Np}{\theta p-N}$ and proceeding as in Step 1 we obtain the estimate \eqref{LSI-2-1} with the constant $C(N,p,\theta)$ given in \eqref{C-SLI}. 

\item {\bf The case $\frac N\theta= p$}: this case follows easily with the necessary modification.
\end{itemize}
This completes the proof of the lemma.
\end{proof}

\begin{lemma}\label{comp}
Let  $p,r\in \lbrack 2,\infty [$
and let $\mathcal{F}$ be given by%
\begin{equation*}
\mathcal{F}\left( \hat{u},\hat{v}\right) =\iint_{{\mathbf S}}k\left( x,y\right)|\hat{u}(x)-\hat{u}(y)|^{p-2}(\hat{u}(x)-\hat{u}(y))(\hat{v}(x)-\hat{v}(y)) \ud y\;\ud x,
\end{equation*}%
for some non-negative kernel $k:\hat\Omega\times \hat\Omega\rightarrow
[0,\infty]$. Then,
\begin{equation}
C_{r,p}\mathcal{F}(|\hat{u}|^{\frac{r-2}{p}}\hat{u},|\hat{u}|^{\frac{r-2}{p}}\hat{u})\leq \mathcal{F}
(\hat{u},\left\vert \hat{u}\right\vert^{r-2}\hat{u}),  \label{comp-energy}
\end{equation}%
for all functions $\hat{u}$ for which the terms in \eqref{comp-energy} make sense, and where
\begin{equation}\label{eq-crp}
C_{r,p}:=\left( r-1\right) \left( \frac{p}{p+r-2}\right) ^{p}.
\end{equation}
\end{lemma}

\begin{proof}
Let the function $g:\mathbb{R%
}\times \mathbb{R}\rightarrow \mathbb{R}$ be given by%
\begin{equation}  \label{func-g}
g\left( z,t\right) =\left\vert z-t\right\vert ^{p-2}\left( z-t\right) \left(
\left\vert z\right\vert ^{r-2}z-\left\vert t\right\vert ^{r-2}t\right)
-C_{r,p}\left\vert \left\vert z\right\vert ^{\frac{r-2}{p}}z-\left\vert
t\right\vert ^{\frac{r-2}{p}}t\right\vert ^{p}.
\end{equation}
Using the definition of $\mathcal{F}$, we first notice that (\ref%
{comp-energy}) is proven once we show that
\begin{equation}  \label{eq-GW}
g\left( z,t\right) \ge 0 \text{ for every } \left( z,t\right) \in \mathbb{R} \times \R .
\end{equation}
Second, we mention that it is easy to verify that
\begin{align*}
g(z,t)=g(t,z),\;\;g(z,0)\ge 0,\; \;g(0,t)\ge 0\;\mbox{ and }\; g(z,t)=g(-z,-t).
\end{align*}
Therefore, without any restriction, we may assume that $z\ge t$. Then,
\begin{equation*}
g\left( z,t\right) =(z-t)^{p-1} \left( \left\vert z\right\vert
^{r-2}z-\left\vert t\right\vert ^{r-2}t\right) -C_{r,p}\left\vert \left\vert
z\right\vert ^{\frac{r-2}{p}}z-\left\vert t\right\vert ^{\frac{r-2}{p}%
}t\right\vert ^{p}.
\end{equation*}%
A simple calculation shows that
\begin{align*}  
\frac{p}{r+p-2}\left[\left\vert z\right\vert ^{\frac{r-2}{p}}z-\left\vert
t\right\vert ^{\frac{r-2}{p}}t\right]=\int_t^z|\tau|^{\frac{r-2}{p}}\;d\tau.
\end{align*}
Since the function $\mathbb{R}\to\mathbb{R}:\;\;\tau\mapsto |\tau|^{p}$ is convex, then using 
Jensen's inequality, it follows from the preceding equality that
\begin{align*}
C_{r,p}\left\vert \left\vert z\right\vert ^{\frac{r-2}{p}}z-\left\vert
t\right\vert ^{\frac{r-2}{p}}t\right\vert ^{p}&=(r-1)\left|\frac{p}{r+p-2}%
\left[\left\vert z\right\vert ^{\frac{r-2}{p}}z-\left\vert t\right\vert ^{%
\frac{r-2}{p}}t\right]\right|^p \\
&=(r-1)\left|\int_t^z|\tau|^{\frac{r-2}{p}}\;\ud\tau\right|^p \\
&=(r-1)(z-t)^p\left|\int_t^z|\tau|^{\frac{r-2}{p}}\;\frac{\ud\tau}{z-t}%
\right|^p \\
&\le (r-1)(z-t)^p\int_t^z|\tau|^{r-2}\;\frac{\ud\tau}{z-t} \\
&=(r-1)(z-t)^{p-1}\int_t^z|\tau|^{r-2}\;d\tau \\
&=(z-t)^{p-1} \left(\left\vert z\right\vert ^{r-2}z-\left\vert t\right\vert
^{r-2}t\right).
\end{align*}
We have shown \eqref{eq-GW} and this completes the proof of the lemma.
\end{proof}

We recall the following well-known algebraic inequalities: for every $r\in [2,\infty[$ there is a constant $C\in ]0,1]$ such that for every $a$, $b\in\R^N$ the following inequalities hold:
\begin{equation}\label{in-p}
    (a-b)\cdot\Big(a|a|^{r-2}-b|b|^{r-2}\Big)\ge C|a-b|^r,
\end{equation}
and
\begin{equation}\label{in-p-2}
    \Big|a|a|^{r-2}-b|b|^{r-2}\Big|\ge C|a-b|^{r-1}.
\end{equation}
Using these estimates we obtain the following result.

 \begin{lemma} \label{lem.F}
Let $p$, $r\in[2,\infty[$ and $F:\R^4\to\R$ be the function given by
\[
F(a,b,c,d):=\Big[|a|^{p-2}a-|b|^{p-2}b\Big]\Big[|c|^{r-2}c-|d|^{r-2}d\Big].
\]
If
\begin{equation}\label{M1}
a-b=c-d ,
\end{equation}
then
\begin{equation}\label{prod}
F(a,b,c,d)\ge C^2|a-b|^{p+r-2} ,
\end{equation}
where $C\in ]0,1]$ is the constant appearing in \eqref{in-p}.
 \end{lemma}

 \begin{proof}
 Assume that \eqref{M1} holds. Then, we have two cases.
\begin{itemize}
\item If $a=b$ (hence, $c=d$ by \eqref{M1}), then
$$F(a,b,c,d)=0=C^2 |a-b|^{r+p-2}.$$
\item If $a\ne b$ (hence, $c\ne d$ by \eqref{M1}), then using \eqref{in-p} and \eqref{M1} we have that
\begin{align*}
 F(a,b,c,d)=&\Big[|a|^{p-2}a-|b|^{p-2}b\Big]\Big[|c|^{r-2}c-|d|^{r-2}d\Big]\\
 =&\frac{(a-b)^2\Big[|a|^{p-2}a-|b|^{p-2}b\Big]\Big[|c|^{r-2}c-|d|^{r-2}d\Big]}{(a-b)^2}\\
=& \frac{(a-b)\Big[|a|^{p-2}a-|b|^{p-2}b\Big](c-d)\Big[|c|^{r-2}c-|d|^{r-2}d\Big]}{(a-b)^2}\\
\ge& C^2\frac{|a-b|^p|c-d|^r}{(a-b)^2}=C^2\frac{|a-b|^{p+r}}{(a-b)^2}=C^2|a-b|^{p+r-2},
\end{align*}
\end{itemize}
where $C\in ]0,1]$ is the constant appearing in \eqref{in-p}. 
 \end{proof}

Next, we show that the norm of the difference of two solutions of the Cauchy problem \eqref{CP} associated with the $j$-subgradient $\partial_j\E_{B_D,\nu}$ is differentiable.

 \begin{lemma}\label{lem-deri-2}
Let $u_0$, $v_0\in L^\infty(\Omega)\subseteq L^2 (\Omega )$, let $u:=S^{D,\nu}(\cdot)u_0$ and $v:=S^{D,\nu}(\cdot)v_0$ be the associated strong solutions of \eqref{CP} (with $\partial\E$ replaced by $\partial_j\E_{B_D,\nu}$). For any $r\in [2,\infty[$ consider the function $f_r:]0,\infty[\to ]0,\infty[$ defined by
\[
f_r(t):=\int_{\Omega}|u(t,x)-v(t,x)|^r\;\ud x.
\]
Then, $f_r$ is differentiable for almost every $t\ge 0$ and
\begin{equation}\label{derivative-2}
\begin{split}
\frac{\ud}{\ud t} f_r(t) & = -r\int_{\Omega}|(u-v)(t,x)|^{r-2}(u-v)(t,x)\Big[ - \dot{u} (t,x) + \dot{v} (t,x)\Big]\;\ud x \\
& \leq -rC^2 \iint_{\mathbf S} k(x,y)|(\hat{u}-\hat{v})(t,x)-(\hat{u}-\hat{v})(t,y)|^{p+r-2}\;\ud y \; \ud x , 
\end{split}
\end{equation}
where $C\in ]0,1]$ is the constant appearing in \eqref{in-p}, and where $\hat{u}(t,\cdot)$ and $\hat{v} (t,\cdot )$ are elliptic extensions of $u(t,\cdot )$ and $v(t,\cdot )$, respectively.
\end{lemma}

\begin{proof}
Let $u_0$, $v_0\in L^\infty(\Omega)$, and let $u:=S^{D,\nu}(\cdot)u_0$ and $v:=S^{D,\nu}(\cdot)v_0$ be the associated strong solutions of \eqref{CP} (with $\partial\E$ replaced by $\partial_j\E_{B_D,\nu}$). Observe that in that case, for almost every $t\in ]0,\infty [$, the elliptic extensions $\hat{u}(t)$, $\hat{v}(t)$ of $u(t)$, $v(t)$ belong to $\mathring{W}^{\theta,{\mathbf S}}_{p,p}(\hat{\Omega},\Omega)$ for almost every $t\ge 0$.
It follows from the $L^\infty$-contractivity of the semigroup $S^{D,\nu}$ (Theorem \ref{Sub}) and the fact that the origin is a fixed point of the semigroup $S^{D,\nu}$ that $u$, $v\in L^\infty(]0,\infty [\times\Omega)$. Since $\Omega$ is bounded this implies that $f_r$ is well-defined. Since $u(\cdot,x)$ and $v(\cdot,x)$ are differentiable almost everywhere on $]0,\infty [$, we have that $f_r$ is differentiable almost everywhere. Using Leibniz's rule, we get
\begin{align*}
\frac{\ud}{\ud t} f_r(t)=&r\int_{\Omega}|(u-v)(t,x)|^{r-1}\sgn(u-v)(t,x)\Big[\dot u(t,x)-\dot v(t,x)\Big]\;\ud x 
\end{align*}
and we have shown the first line in \eqref{derivative-2}. As $u$ and $v$ are strong solutions of the Cauchy problem \eqref{CP} with $\partial\E$ replaced by $\partial_j\E_{B_D,\nu}$, one has $(u(t),-\dot{u} (t))$, $(v(t),-\dot{v}(t))\in\partial_j\E_{B_D,\nu}$ for almost every $t\geq 0$. This, the formula \eqref{INt-parts} and the preceding equality imply that
\begin{align}\label{J1}
    \frac{\ud}{\ud t} f_r(t)=-r\iint_{\mathbf S} k(x,y) & \Big[|\hat{u}(t,x)-\hat{u}(t,y)|^{p-2}(\hat{u}(t,x)-\hat{u}(t,y))-|\hat{v}(t,x)-\hat{v}(t,y)|^{p-2}(\hat{v}(t,x)-\hat{v}(t,y)\Big]\times\notag\\
   & \Big[(\hat{u}-\hat{v})(t,x)|(\hat{u}-\hat{v})(t,x)|^{r-2}-(\hat{u}-\hat{v})(t,y)|(\hat{u}-\hat{v})(t,y)|^{r-2}\Big]\;\ud y\;\ud x.
\end{align}
Applying Lemma \ref{lem.F} with
\[
a:=\hat{u}(t,x)-\hat{u}(t,y),\; b:=\hat{v}(t,x)-\hat{v}(t,y),\; c:=(\hat{u}-\hat{v})(t,x)\;\text{ and } d:=(\hat{u}-\hat{v})(t,y),
\]
we get 
\begin{align}\label{J2}
&\Big[|\hat{u}(t,x)-\hat{u}(t,y)|^{p-2}(\hat{u}(t,x)-\hat{u}(t,y))-|\hat{v}(t,x)-\hat{v}(t,y)|^{p-2}(\hat{v}(t,x)-\hat{v}(t,y)\Big]\times\notag\\
   & \Big[(\hat{u}-\hat{v})(t,x)|(\hat{u}-\hat{v})(t,x)|^{r-2} - (\hat{u}-\hat{v})(t,y)|(\hat{u}-\hat{v})(t,y)|^{r-2}\Big]\notag\\
   &\ge C^2 |(\hat{u}-\hat{v})(t,x)-(\hat{u}-\hat{v})(t,y)|^{p-1}|(\hat{u}-\hat{v})(t,x)-(\hat{u}-\hat{v})(t,y)|^{r-1}\notag\\
    &= C^2 |(\hat{u}-\hat{v})(t,x)-(\hat{u}-\hat{v})(t,y)|^{p+r-2}.
\end{align}
Combining \eqref{J1}-\eqref{J2} gives the second line in \eqref{derivative-2}.
\end{proof}

Next, we show further differentiability properties of the norm of the difference of two solutions to \eqref{CP}.

\begin{lemma}\label{lem-deri-log-22}
Let $u_0$, $v_0\in L^\infty(\Omega)\subseteq L^2 (\Omega )$, let $u:=S^{D,\nu}(\cdot)u_0$ and $v:=S^{D,\nu}(\cdot)v_0$ be the associated strong solutions of \eqref{CP} (with $\partial\E$ replaced by $\partial_j\E_{B_D,\nu}$). Let $r:[0,\infty [\to [2,\infty [$ be any increasing differentiable function. Then, for a.e. $t\ge 0$, 
\begin{align}\label{derivative-log-2}
\frac{\ud}{\ud t}\log\|(u -v)(t,\cdot)\|_{r(t)}  =& \frac{r'(t)}{r(t)}\int_{\Omega}\frac{|(u-v)(t,x)|^{r(t)}}{\|(u-v)(t,\cdot)\|_{r(t)}^{r(t)}}\log\left(\frac{|(u-v)(t,x)|}{\|(u-v)(t,\cdot)\|_{r(t)}}\right)\;\ud x\notag\\
&-\frac{1}{\|(u -v)(t,\cdot)\|_{r(t)}^{r(t)}}\int_{\Omega}|(u-v)(t,x)|^{r(t)-2}(u-v)(t,x)\Big[-\dot{u}(t,x) + \dot{v} (t,x) \Big]\; \ud x \notag \\
\leq & \frac{r'(t)}{r(t)} \int_{\Omega}\frac{|(u-v)(t,x)|^{r(t)}}{\|(u-v)(t,\cdot)\|_{r(t)}^{r(t)}}\log\left(\frac{|(u-v)(t,x)|}{\|(u-v)(t,\cdot)\|_{r(t)}}\right)\;\ud x\notag\\
&-\frac{C^2}{\|(u -v)(t,\cdot)\|_{r(t)}^{r(t)}} \iint_{{\mathbf S}}k(x,y)|(\hat{u}-\hat{v})(t,x) - (\hat{u}-\hat{v})(t,y)|^{p+r(t)-2}\;\ud y \; \ud x
\end{align}
where $C\in ]0,1]$ is the constant appearing in \eqref{in-p}, and  $\hat{u}(t,\cdot)$ and $\hat{v} (t,\cdot )$ are elliptic extensions of $u(t,\cdot )$ and $v(t,\cdot )$, respectively.
\end{lemma}

\begin{proof}
Since for almost every $t\ge 0$,
\begin{align*}
\frac{\ud}{\ud t}\|(u -v)(t,\cdot)\|_{r(t)}^{r(t)}=r'(t)\frac{\partial}{\partial r}\|(u -v)(t,\cdot)\|_r^r\Big|_{r=r(t)}+\frac{\partial}{\partial t}\|(u -v)(t,\cdot)\|_r^r\Big|_{r=r(t)}
\end{align*}
and
\begin{align*}
 \frac{\partial}{\partial r}\|(u -v)(t,\cdot)\|_r^r\Big|_{r=r(t)}=\int_{\Omega}|(u-v)(t,x)|^{r(t)}\log(|(u-v)(t,x)|)\;\ud x, 
\end{align*}
it follows from \eqref{derivative-2} that
\begin{align}\label{66-2}
 \frac{\ud}{\ud t}\|(u -v)(t,\cdot)\|_{r(t)}^{r(t)} = & r'(t)\int_{\Omega}|(u-v)(t,x)|^{r(t)}\log(|(u-v)(t,x)|)\;\ud x \notag\\
 &-r(t)\int_{\Omega}|(u-v)(t,x)|^{r(t)-2}(u-v)(t,x)\Big[-\dot{u} (t,x)+\dot{v} (t,x)\Big]\;\ud x. 
\end{align}
Since
\begin{align*}
\frac{\ud}{\ud t}\log(\|(u -v)(t,\cdot)\|_{r(t)})=-\frac{r'(t)}{r(t)}\log(\|(u -v)(t,\cdot)\|_{r(t)})+\frac{1}{r(t)\|(u-v)(t,\cdot)\|_{r(t)}^{r(t)}}\frac{\ud}{\ud t}\|(u-v)(t,\cdot)\|_{r(t)}^{r(t)},
\end{align*}
it follows from \eqref{66-2} that
\begin{align*}
  \frac{\ud}{\ud t}\log(\|(u-v)(t,\cdot)\|_{r(t)})  =& -\frac{r'(t)}{r(t)}\log(\|(u-v)(t,\cdot)\|_{r(t)})\\
  &- \frac{1}{\|(u-v)(t,\cdot)\|_{r(t)}^{r(t)}}\int_{\Omega}|(u-v)(t,x)|^{r(t)-2}(u-v)(t,x)\Big[-\dot{u}(t,x)+ \dot{v} (t,x)\Big]\;\ud x\\
  &+\frac{r'(t)}{r(t)\|((u -v)(t,\cdot)\|_{r(t)}^{r(t)}}\int_{\Omega}|(u-v)(t,x)|^{r(t)}\log|(u-v)(t,x)|\;\ud x\\
  =&\frac{r'(t)}{r(t)}\int_{\Omega}\frac{|(u-v)(t,x)|^{r(t)}}{\|(u -v)(t,\cdot)\|_{r(t)}^{r(t)}}\log\left(\frac{|(u-v)(t,x)|}{\|(u-v)(t,\cdot)\|_{r(t)}}\right)\;\ud x\\
  &-\frac{1}{\|(u-v)(t,\cdot)\|_{r(t)}^{r(t)}}\int_{\Omega}|(u-v)(t,x)|^{r(t)-2}(u-v)(t,x)\Big[-\dot{u}(t,x) +\dot{v} (t,x)\Big]\;\ud x
\end{align*}
and we have shown the first line in \eqref{derivative-log-2}. Noting that $(u(t),-\dot{u}(t))$, $(v(t),-\dot{v}(t))\in\partial_j\E_{B_D,\nu}$ for almost every $t\geq 0$, and recalling \eqref{INt-parts}, then the second line in \eqref{derivative-log-2} follows from the first line similarly as in the proof of Lemma \ref{lem-deri-2}. 
\end{proof}

\begin{lemma}\label{lem-deri-log-23}
 Let $u_0$, $v_0\in L^\infty(\Omega)\subset L^2(\Omega)$, let $u:=S^{D,\nu}(\cdot)u_0$ and $v:=S^{D,\nu}(\cdot)v_0$ be the associated strong solutions of \eqref{CP}  (with $\partial\E$ replaced by $\partial_j\E_{B_D,\nu}$). Let $r:[0,\infty [\to [2,\infty [$ be any increasing differentiable function. Then, there is a constant $C>0$ such that for every $\varepsilon>0$ and for a.e. $t\ge 0$, 
\begin{align}\label{AA02}
&\frac{\ud}{\ud t}\log\|(u -v)(t,\cdot)\|_{r(t)}  \le \frac{r'(t)}{r(t)}\int_{\Omega}\frac{|(u-v)(t,x)|^{r(t)}}{\|(u -v)(t,\cdot)\|_{r(t)}^{r(t)}}\log\left(\frac{|(u-v)(t,x)|}{\|(u -v)(t,\cdot)\|_{r(t)}}\right)\;\ud x\notag\\
&-C_{r(t),p}\frac{\|(u -v)(t,\cdot)\|_{r(t)+p-2}^{r(t)+p-2}}{\|(u -v)(t,\cdot)\|_{r(t)}^{r(t)}} \frac{\theta p^2(r(t)+p-2)}{pNC\varepsilon}\int_{\Omega}\frac{|(u-v)(t,x)|^{r(t)+p-2}}{\|(u -v)(t,\cdot)\|_{r(t)}^{r(t)}} \log\left(\frac{|(u-v)(t,x)|}{\|(u -v)(t,\cdot)\|_{r(t)+p-2}}\right)\;\ud x \notag\\
&-C_{r(t),p}\frac{\|(u -v)(t,\cdot)\|_{r(t)+p-2}^{r(t)+p-2}}{\|(u -v)(t,\cdot)\|_{r(t)}^{r(t)}}\frac{\log\varepsilon}{C\varepsilon}
\end{align}
where 
$\displaystyle C_{r(t),p}:=\left( r(t)-1\right) \left( \frac{p}{p+r(t)-2}\right) ^{p}.$
\end{lemma}

\begin{proof}
Let
\[
U(t,x):=\frac{|(\hat{u}-\hat{v})(t,x)|^{\frac{r(t)+p-2}{p}}}{\|(u -v)(t,\cdot)\|_{r(t)+p-2}^{\frac{r(t)+p-2}{p}}} ,
\]
where $\hat{u} (t,\cdot )$ and $\hat{v} (t,\cdot )$ are elliptic extensions of $u(t,\cdot )$ and $v(t,\cdot )$, respectively. Recall that all the norms are on $\Omega$.

Then, for almost every $t\ge 0$, we have that  $U(t,\cdot)\in \mathring{W}^{\theta,{\mathbf S}}_{p,p}(\hat{\Omega},\Omega)$, $U(t,\cdot)\ge 0$ a.e. in $\Omega$ and $\|U(t,\cdot)\|_p^p=1$. It follows from Lemma \ref{LSI} that there is a constant $C>0$ such that for every $\varepsilon>0$,  
\begin{align}\label{AA12}
\frac{\theta p^2}{N}\int_{\Omega}|U(t,x)|^p\log(U(t,x))\;\ud x\le \varepsilon C\iint_{\mathbf S}k(x,y)|U(t,x)-U(t,y)|^{p}\;\ud y\;\ud x - \log(\varepsilon).
\end{align}
Notice that
\begin{align}\label{AA22}
&\iint_{\mathbf S}k(x,y)|U(t,x)-U(t,y)|^{p}\;\ud y\;\ud x \notag\\
=&\frac{1}{\|(u -v)(t,\cdot)\|_{r(t)+p-2}^{r(t)+p-2}} \iint_{\mathbf S} k(x,y)\left||(\hat{u}-\hat{v})(t,x)|^{\frac{r(t)+p-2}{p}}-|(\hat{u}-\hat{v}) (t,y)|^{\frac{r(t)+p-2}{p}}\right|^p\;\ud y\;\ud x\notag\\
=&\frac{1}{\|(u -v)(t,\cdot)\|_{r(t)+p-2}^{r(t)+p-2}} \iint_{\mathbf S} k(x,y) \left||(\hat{u}-\hat{v})(t,x)|^{\frac{r(t)-2}{p}}|(\hat{u}-\hat{v})(t,x)|-|(\hat{u}-\hat{v})(t,y)|^{\frac{r(t)-2}{p}}|(\hat{u}-\hat{v})(t,y)|\right|^p\;\ud y\;\ud x
\end{align}
and
\begin{align}\label{AA32}
& \int_{\Omega}|U(t,x)|^p\log(U(t,x))\;\ud x\notag\\
= & \frac{r(t)+p-2}{p\|(u -v)(t,\cdot)\|_{r(t)+p-2}^{r(t)+p-2}} \int_{\Omega}|(u-v)(t,x)|^{r(t)+p-2}\log\left(\frac{|(u-v)(t,x)|}{\|(u -v)(t,\cdot)\|_{r(t)+p-2}}\right)\; \ud x.
\end{align}
Using \eqref{AA22}, \eqref{AA32} and \eqref{comp-energy} in Lemma \ref{comp} we can deduce from \eqref{AA12} that there is a constant $C>0$ such that for every $\varepsilon>0$, 
\begin{align*}
&\frac{\theta p(r(t)+p-2)}{N\|(u -v)(t,\cdot)\|_{r(t)+p-2}^{r(t)+p-2}} \int_{\Omega} |(u-v)(t,x)|^{r(t)+p-2}\log\left(\frac{|(u-v)(x,t)|}{\|(u -v)(t,\cdot)\|_{r(t)+p-2}}\right)\;\ud x\\
\le & \frac{\varepsilon C}{\|(u -v)(t,\cdot)\|_{r(t)+p-2}^{r(t)+p-2}} \iint_{\mathbf S} k(x,y)\left||(\hat{u}-\hat{v})(t,x)|^{\frac{r(t)-2}{p}} |(\hat{u}-\hat{v})(t,x)|-|(\hat{u}-\hat{v})(t,y)|^{\frac{r(t)-2}{p}} |(\hat{u}-\hat{v})(t,y)|\right|^p\;\ud y\;\ud x
-\log\varepsilon\\
\le & \frac{\varepsilon C}{C_{r(t),p}\|(u -v)(t,\cdot)\|_{r(t)+p-2}^{r(t)+p-2}} \iint_{\mathbf S} k(x,y)|(\hat{u}-\hat{v})(t,x) - (\hat{u}-\hat{v})(t,y)|^{r(t)+p-2}\;\ud y\;\ud x-\log\varepsilon.
\end{align*}
This implies that
\begin{align*}
  -&\iint_{\mathbf S} k(x,y) |(\hat{u}-\hat{v})(t,x)-(\hat{u}-\hat{v})(t,y)|^{r(t)+p-2}\;\ud y\;\ud x\\
  \le &-\frac{\theta p(r(t)+p-2)C_{r,p}}{NC\varepsilon}\int_{\Omega}|(u-v)(t,x)|^{r(t)+p-2}\log\left(\frac{|(u-v)(t,x)|}{\|(u -v)(t,\cdot)\|_{r(t)+p-2}} \right)\;\ud x\\
  &-C_{r(t),p}\|(u -v)(t,\cdot)\|_{r(t)+p-2}^{r(t)+p-2} \frac{\log\varepsilon}{C\varepsilon}.
\end{align*}
Using  the preceding estimate we get that for every $\varepsilon>0$,
\begin{align*}
&\frac{\ud}{\ud t}\log\|(u -v)(t,\cdot)\|_{r(t)}  \le \frac{r'(t)}{r(t)} \int_{\Omega} \frac{|(u-v)(t,x)|^{r(t)}}{\|(u -v)(t,\cdot)\|_{r(t)}^{r(t)}} \log\left(\frac{|(u-v)(t,x)|}{\|(u -v)(t,\cdot)\|_{r(t)}}\right) \;\ud x\\
&-\frac{\theta p(r(t)+p-2)CC_{r(t),p}\|(u -v)(t,\cdot)\|_{r(t)+p-2}^{r(t)+p-2}}{N\varepsilon\|(u -v)(t,\cdot)\|_{r(t)}^{r(t)}} \int_{\Omega} \frac{|(u-v)(t,x)|^{r(t)+p-2}}{\|(u -v)(t,\cdot)\|_{r(t)+p-2}^{r(t)+p-2}} \log\left(\frac{|(u-v)(t,x)|}{\|(u -v)(t,\cdot)\|_{r(t)+p-2}} \right)\;\ud x \\
&-C_{r(t),p} \frac{\|(u -v)(t,\cdot)\|_{r(t)+p-2}^{r(t)+p-2}}{\|(u -v)(t,\cdot)\|_{r(t)}^{r(t)}} \frac{\log\varepsilon}{C\varepsilon},
\end{align*}
and we have shown \eqref{AA02}.
\end{proof}
 
Finally, we show that the norm of the difference of two solutions to \eqref{CP} satisfies a suitable ordinary differential inequality.

\begin{lemma}\label{lem-diff-ine-2}
 Let $u_0$, $v_0\in L^\infty(\Omega)\subset L^2(\Omega)$, let $u:=S^{D,\nu}(\cdot)u_0$ and $v:=S^{D,\nu}(\cdot)v_0$ be the associated strong solutions of \eqref{CP} (with $\partial\E$ replaced by $\partial_j\E_{B_D,\nu}$). Let $r:[0,\infty[\to [2,\infty[$ be any increasing differentiable function. Then, for almost every $t\ge 0$,
\begin{equation}\label{eq-diff-ine2}
    \frac{\ud}{\ud t}\log\|(u -v)(t,\cdot)\|_{r(t)} \le -\mathcal P(t) \log\|(u -v)(t,\cdot)\|_{r(t)} -\mathcal Q(t)
\end{equation}
where
\begin{equation}\label{eq-P}
    \mathcal P(t):=\frac{r'(t)}{r(t)}\frac{N}{\theta p}\frac{p-2}{r(t)+p-2}
\end{equation}
and 
\begin{equation}\label{eq-Q}
    \mathcal Q(t):=\frac{N}{\theta p}\frac{r'(t)}{r(t)}\frac{1}{r(t)+p-2}\log\left(\frac{r(t)}{r'(t)}\frac{C_{r(t),p}\theta pC(r(t)+p-2)}{N}\right)
  -\frac{r'(t)}{r(t)}\left(\frac{N}{\theta p}-1\right)\frac{p-2}{r(t)+p-2}\log|\Omega|.
\end{equation}
\end{lemma}

\begin{proof}
We give the full proof of the case $N>\theta p$. The other cases follow similarly with very minor modifications.
We define the function $\mathcal K:[1,\infty[\times L^\infty(\Omega)\to [0,\infty]$ by
\begin{align*}
\mathcal K(q,w):=\int_\Omega \frac{|w(x)|^q}{\|w\|_q^q}\log\left(\frac{|w(x)|}{\|w\|_q}\right)\;\ud x.
\end{align*}
It is clear that for every fixed $w$, the mapping $q\mapsto \log\|w\|_q^q$ is convex and 
\begin{align*}
\frac{\ud}{\ud q}\log\|w\|_q^q=\mathcal K(q,w)+\log\|w\|_q.
\end{align*}
Moreover, the mapping $q\mapsto \frac{\ud}{\ud q}\log\|w\|_q^q$ is non-decreasing. Therefore, for every $q_2\ge q_1\ge 1$ we have that
\begin{align}\label{ine-K2}
\mathcal K(q_1,w)-\mathcal K(q_2,w)\le \log\frac{\|w\|_{q_2}}{\|w\|_{q_1}}.
\end{align}

It follows from \eqref{AA02} that for almost every $t\ge 0$,
\begin{align}\label{AA42}
\frac{\ud}{\ud t}\log\|(u -v)(t,\cdot)\|_{r(t)}  \le& \frac{r'(t)}{r(t)}\mathcal K(r(t),(u-v)(t,\cdot))\notag\\
&-C_{r(t),p}\frac{\|(u -v)(t,\cdot)\|_{r(t)+p-2}^{r(t)+p-2}}{\|(u -v)(t,\cdot)\|_{r(t)}^{r(t)}} \frac{\theta p(r(t)+p-2)C}{N\varepsilon}\mathcal K(r(t)+p-2,(u-v)(t,\cdot)) \notag\\
&-C_{r(t),p}\frac{\|(u -v)(t,\cdot)\|_{r(t)+p-2}^{r(t)+p-2}}{\|(u -v)(t,\cdot)\|_{r(t)}^{r(t)}} \frac{\log\varepsilon}{C\varepsilon}.
\end{align}
Letting
\begin{align*}
\varepsilon:=C_{r(t),p} \frac{\|(u-v)(t,\cdot)\|_{r(t)+p-2}^{r(t)+p-2}}{\|(u-v)(t,\cdot)\|_{r(t)}^{r(t)}}\frac{\theta pr(t)(r(t)+p-2)C}{Nr'(t)}
\end{align*}
in \eqref{AA42} gives
\begin{align}\label{6212}
  \frac{\ud}{\ud t}\log\|(u -v)(t,\cdot)\|_{r(t)} \le & \frac{r'(t)}{r(t)}\Big[\mathcal K(r(t),(u-v)(t,\cdot)) - \mathcal K(r(t)+p-2,(u-v)(t,\cdot))\Big] \notag \\
  &-\frac{N}{\theta p}\frac{r'(t)}{r(t)}\frac{1}{r(t)+p-2} \log\left(\frac{\|(u -v)(t,\cdot)\|_{r(t)+p-2}^{r(r)+p-2}}{\|(u -v)(t,\cdot))\|_{r(t)}^{r(t)}}\right)\notag\\
  &-\frac{N}{\theta p} \frac{r'(t)}{r(t)} \frac{1}{r(t)+p-2} \log\left(\frac{r(t)}{r'(t)}\frac{C_{r(t),p}\theta pC(r(t)+p-2)}{N}\right).
\end{align}

Next, using \eqref{6212} and \eqref{ine-K2} with $q_1=r(t)$ and $q_2=r(t)+p-2$ we get
\begin{align*}
  \frac{\ud}{\ud t}\log\|(u -v)(t,\cdot)\|_{r(t)} \le & \frac{r'(t)}{r(t)} \log\frac{\|(u -v)(t,\cdot)\|_{r(t)+p-2}}{\|(u-v)(t,\cdot)\|_{r(t)}} \notag \\
  &-\frac{r'(t)}{r(t)} \frac{N}{\theta p(r(t)+p-2)} \log\left(\frac{\|(u -v)(t,\cdot)\|_{r(t)+p-2}^{r(r)+p-2}}{\|(u -v)(t,\cdot)\|_{r(t)}^{r(t)}}\right) \notag\\
  &-\frac{r'(t)}{r(t)} \frac{N}{\theta p(r(t)+p-2)} \log\left(\frac{r(t)}{r'(t)}\frac{C_{r(t),p}\theta pC(r(t)+p-2)}{N}\right).
\end{align*}
It follows from the preceding inequality that
\begin{align*}
  \frac{\ud}{\ud t}\log\|(u -v)(t,\cdot)\|_{r(t)} \le &  \frac{r'(t)}{r(t)} \Big[\log\|(u -v)(t,\cdot)\|_{r(t)+p-2} - \log\|(u -v)(t,\cdot)\|_{r(t)} -
  \frac{N}{\theta p} \log\|(u-v)(t,\cdot)\|_{r(t)+p-2}\\
 & + \frac{Nr(t)}{\theta p(r(t)-p+2)} \log\|(u-v)(t,\cdot)\|_{r(t)}\\
 & - \frac{N}{\theta p(r(t)+p-2)} \log\left(\frac{r(t)}{r'(t)} \frac{C_{r(t),p}\theta pC(r(t)+p-2)}{N}\right)
  \Big].
\end{align*}
This implies that
\begin{align}\label{BB2}
 \frac{\ud}{\ud t}\log \|(u -v)(t,\cdot)\|_{r(t)} \le &  \frac{r'(t)}{r(t)} \Big[\left(1-\frac{N}{\theta p}\right)\log\|(u -v)(t,\cdot)\|_{r(t)+p-2} 
 + \left(\frac{Nr(t)}{\theta p(r(t)-p+2)}-1\right)\log\|(u -v)(t,\cdot)\|_{r(t)}\notag\\
 &-\frac{N}{\theta pC(r(t)+p-2)} \log \left(\frac{r(t)}{r'(t)}\frac{C_{r(t),p}\theta pC(r(t)+p-2)}{N}\right)\Big].
\end{align}
Using H\"older's inequality, we have that
\begin{equation*}
\|(u -v)(t,\cdot)\|_{r(t)} \le |\Omega|^{\frac{p-2}{r(t)(r(t)+p-2)}}\|(u -v)(t,\cdot)\|_{r(t)+p-2}.
\end{equation*}
This implies that
\begin{align}\label{BB02}
\log\|(u-v)(t,\cdot)\|_{r(t)} \le \frac{p-2}{r(t)(r(t)+p-2)} \log|\Omega|+\log\|(u -v)(t,\cdot)\|_{r(t)+p-2}.
\end{align}
Since $1-\frac{N}{\theta p}<0$, it follows from \eqref{BB02} that 
\begin{align}\label{BB12}
\left(1-\frac{N}{\theta p}\right) \log\|(u -v)(t,\cdot)\|_{r(t)+p-2} \le & \left(1-\frac{N}{\theta p}\right) \log\|(u -v)(t,\cdot)\|_{r(t)}\notag\\
& - \left(1-\frac{N}{\theta p}\right) \frac{p-2}{r(t)(r(t)+p-2)}\log|\Omega|.
\end{align}
Using \eqref{BB12} we can deduce from \eqref{BB2} that 
\begin{align*}
  \frac{\ud}{\ud t}\log\|(u -v)(t,\cdot)\|_{r(t)} \le & -\frac{r'(t)}{r(t)} \frac{N}{\theta p} \frac{p-2}{r(t)+p-2} \log\|(u -v)(t,\cdot)\|_{r(t)}\notag\\
  & - \frac{N}{\theta p} \frac{r'(t)}{r(t)} \frac{1}{r(t)+p-2} \log\left(\frac{r(t)}{r'(t)}\frac{C_{r(t),p}\theta pC(r(t)+p-2)}{N}\right)\notag\\
  & + \frac{r'(t)}{r(t)} \left(\frac{N}{\theta p}-1\right) \frac{p-2}{r(t)+p-2} \log|\Omega|.
\end{align*}
We have shown \eqref{eq-diff-ine2}.
\end{proof}

Now, we are ready to give the proof of the main result of this section.

\begin{proof}[\bf Proof of Theorem \ref{theo-ultra}]
Let $q\in [1,\infty[$. 
Since $S^{D,\nu}$ is a submarkovian semigroup on $L^2(\Omega)$, it follows that $S^{D,\nu}$ extends consistently to a contraction strongly continuous semigroup on $L^q(\Omega)$.\\

{\bf Part I}: we assume first that $u_0$, $v_0\in L^\infty(\Omega)$. We proceed in two steps.\\

{\bf Step 1: the case $2\le p<\frac{N}{\theta}$.}
Let $u_0$, $v_0\in L^\infty(\Omega)$, $u:=S^{D,\nu}(\cdot)u_0$ and $v:=S^{D,\nu}(\cdot)v_0$. Let $r:[0,\infty[\to [2,\infty[$ be any increasing differentiable function. Let $Y(t):=\log\|(u -v)(t,\cdot)\|_{r(t)}$. It follows from \eqref{eq-diff-ine2} in Lemma \ref{lem-diff-ine-2} that $Y(t)$ satisfies the ordinary differential inequality
\begin{equation}\label{ODE-INE2}
    Y'(t)+\mathcal P(t)Y(t)+\mathcal Q(t)\le 0,
\end{equation}
where $\mathcal P(t)$ and  $\mathcal Q(t)$ are given by \eqref{eq-P} and \eqref{eq-Q}, respectively.
Consider the following associated first order ordinary differential equation for $\tau\ge 0$:
\begin{equation}\label{ODE-Eq2}
X'(\tau)+\mathcal P(\tau)X(\tau)+\mathcal Q(\tau)=0,\qquad X(0)\le Y(0).
\end{equation}
It is well-known that the unique solution $Y$ of \eqref{ODE-INE2} satisfies
 \[
 Y(t)\le X(t)\text{ for every } t\ge 0.
 \]
The unique solution of \eqref{ODE-Eq2} is given by
\begin{align*}
    X(\tau)=\exp{\left[-\int_0^\tau\mathcal P(\xi)\;\ud\xi\right]}\left(X(0)-\int_0^\tau\mathcal Q(\xi)\exp{\left[\int_0^\xi\mathcal P(z)\;dz\right]}\;\ud\xi\right).
\end{align*}
 Let $r(\tau):=\frac{qt}{t-\tau}$ where $0\le \tau<t$ so that $\frac{r'(\tau)}{r(\tau)}=\frac{1}{t-\tau}$. 
 Calculating we get 
 \begin{align*}
\mathcal P(\tau)=\frac{N(p-2)}{\theta p}\frac{1}{tq+(t-\tau)(p-2)} \;\hbox{ and }\; \mathcal Q(\tau)=\sum_{j=1}^5\mathcal Q_j(\tau),
 \end{align*}
 where
 \begin{align*}
\mathcal Q_1(\tau)=&\frac{N}{\theta p}\Big[tq+(t-\tau)(p-2)\Big]^{-1}\log\left(\frac{Cp^{p+1}}{N}\right) , \\
\mathcal Q_2(\tau)=&\frac{N}{\theta p}(p-1)\Big[tq+(t-\tau)(p-2)\Big]^{-1}\log(t-\tau) , \\
\mathcal Q_3(\tau)=&\frac{N}{\theta p}\Big[tq+(t-\tau)(p-2)\Big]^{-1}\log(tq-(t-\tau)) , \\
\mathcal Q_4(\tau)=&\frac{N}{\theta p}(p-1)\Big[tq+(t-\tau)(p-2)\Big]^{-1}\log\left(tq+(t-\tau)(p-2)\right) ,  \\
\mathcal Q_5(\tau)=&\frac{\theta p-N}{p}\frac{p-2}{tq+(t-\tau)(p-2)}\log|\Omega|.
 \end{align*}
A simple calculation gives
\begin{align}\label{Int-P2}
    \int_0^\tau\mathcal P(\xi)\;\ud\xi=\frac{N(p-2)}{\theta p}\int_0^\tau \frac{1}{tq+(t-\xi)(p-2)}\;\ud\xi=\frac{N}{\theta p}\log\left(\frac{t(q+p-2)}{tq+(t-\tau)(p-2)}\right).
\end{align}
We also have that
\begin{align}
\int_0^\tau\mathcal Q(\xi)\exp{\left[\int_0^\xi\mathcal P(z)\;dz\right]}\;\ud\xi
=\int_0^\tau \mathcal Q(\xi)\Big[\frac{\tau(q+p-2)}{\tau q+(\tau-\xi)(p-2}\Big]^{\frac{N}{\theta p}}\;\ud\xi=\sum_{i=1}^5K_i (\tau),
\end{align}
where
\begin{align*}
K_1 (\tau ) =&\frac{N}{\theta p}\Big[t(q+p-2)\Big]^{\frac{N}{\theta p}}\log\left(\frac{Cp^{p+1}}{N}\right)\int_0^\tau\Big[\tau q+(\tau-\xi)(p-2)\Big]^{-1-\frac{N}{\theta p}}\;\ud\xi , \\
K_2 (\tau ) =& \frac{N}{\theta p}\Big[t(q+p-2)\Big]^{\frac{N}{\theta p}}(p-1)\int_0^\tau\Big[\tau q+(\tau-\xi)(p-2)\Big]^{-1-\frac{N}{\theta p}}\log(\tau-\xi)\;\ud\xi , \\
K_3 (\tau ) =&\frac{N}{\theta p}\Big[t(q+p-2)\Big]^{\frac{N}{\theta p}}\int_0^\tau\Big[\tau q+(\tau-\xi)(p-2)\Big]^{-1-\frac{N}{\theta p}}\log[\tau q-(\tau-\xi)]\;\ud\xi , \\
K_4 (\tau ) =&-\frac{N}{\theta p}\Big[t(q+p-2)\Big]^{\frac{N}{\theta p}}(p-1)\int_0^\tau\Big[\tau q+(\tau-\xi)(p-2)\Big]^{-1-\frac{N}{\theta p}}\log[\tau q-(\tau-\xi)(p-2)]\;\ud\xi , \text{ and} \\
K_5 (\tau ) =&-\Big[(p-2)\left(\frac{N}{\theta p}-1\right)\log|\Omega|\Big]\Big[t(q+p-2)\Big]^{\frac{N}{\theta p}}\int_0^\tau\Big[\tau q+(\tau-\xi)(p-2)\Big]^{-1-\frac{N}{\theta p}}\;\ud\xi.
\end{align*}
An elementary calculation gives the following results:
\begin{align*}
    K_1 (\tau ) =&\frac{1}{p-2}\left(\Big[\frac{q+p-2}{q}\Big]^{\frac{N}{\theta p}}-1\right)\log\left(\frac{C\theta p^{p+1}}{N}\right) , \\
      K_2 (\tau ) =&\frac{p-1}{p-2} \left(\Big[\frac{q+p-2}{q}\Big]^{\frac{N}{\theta p}}-1\right)\log(\tau) +R_2 , \\
        K_3 (\tau ) =& \frac{1}{p-2} \left(\Big[\frac{q+p-2}{q}\Big]^{\frac{N}{\theta p}}-1\right)\log(\tau) +R_3 , \\
      K_4 (\tau ) =& \frac{p-1}{p-2} \left(\Big[\frac{q+p-2}{q}\Big]^{\frac{N}{\theta p}}-1\right)\log(\tau) +R_4 , \text{ and} \\  
      K_5 (\tau ) =&-\frac{N-\theta p}{N}\left(\Big[\frac{q+p-2}{q}]^{\frac{N}{\theta p}}-1\right)\log(|\Omega|),
\end{align*}
where $R_2$, $R_3$, $R_4$ depend only on $N$, $\theta$, $p$ and $q$ and are independent of $\tau$ and $|\Omega|$. Their explicit values are given by
\begin{align*}
R_2= & \frac{N}{\theta p}(p-1)\Big(q+p-2\Big)^{\frac{N}{\theta p}}\int_0^1\Big[q+\xi(p-2)\Big]^{-1-\frac{N}{\theta p}}\log(\xi)\;\ud\xi , \\
R_3=&\frac{N}{\theta p}\Big(q+p-2\Big)^{\frac{N}{\theta p}}\int_0^1\Big[q+\xi(p-2)\Big]^{-1-\frac{N}{\theta p}}\log(q-\xi)\;\ud\xi , \text{ and} \\
R_4=&-\frac{N}{\theta p}(p-1)\Big(q+p-2\Big)^{\frac{N}{\theta p}}\int_0^1\Big[q+\xi(p-2)\Big]^{-1-\frac{N}{\theta p}}\log(q+\xi(p-2))\;\ud\xi.
\end{align*}
Note that also $K_1$ is independent of $\tau$. We have shown that 
\begin{align*}
\lim_{\tau\to t}X(\tau)=&\Big(\frac{q}{q+p-2}\Big)^{\frac{N}{\theta p}}\Big [X(0)-\frac{1}{p-2}\Big(\Big(\frac{q+p-2}{q}\Big)^{\frac{N}{\theta p}}-1\Big]\log(t)\\
&-\Big[-\frac{N-\theta p}{N}\Big(\Big(\frac{q}{q+p-2}\Big)^{\frac{N}{\theta p}}-1\Big)\log|\Omega|+K_1+R_2+R_3+R_4\Big].
\end{align*}
It follows from the submarkovian property of the semigroup $S^{D,\nu}$ and $Y(\tau)\le X(\tau)$ for all $\tau\ge 0$ that,  for all $0\le \tau<t$, 
\begin{align}\label{CC02}
\|(u -v)(t,\cdot)\|_{r(\tau)} \le \|(u -v)(\tau,\cdot)\|_{r(\tau)} \le \exp\Big[\log\|(u-v)(\tau,\cdot )\|_{r(\tau)} \Big] = \exp\left(Y(\tau)\right)\le \exp\left(X(\tau)\right).
\end{align}
Since $\lim_{\tau\to t^-}r(\tau)=\infty$, taking the limit of \eqref{CC02} as $\tau\to t^-$ we get 
\begin{align*}
    \|(u -v)(t,\cdot)\|_{\infty} \le \lim_{\tau \to t^-} \exp\left(X(\tau)\right).
\end{align*}
Calculating and using the fact that $r(0)=q$ we obtain the estimate \eqref{est-65-3}.\\

{\bf Step 2: the case $2<\frac{N}{\theta}\le p<\infty$.}
Let $u_0$, $v_0\in L^\infty(\Omega)$, $u=S^{D,\nu}(\cdot)u_0$ and $v=S^{D,\nu}(\cdot)v_0$. We recall that in this case we have the continuous restriction operators $\mathring{W}^{\theta,\mathbf{S}}_{p,p}(\hat{\Omega}, \Omega)\to L^{\tau}(\Omega)$ for every $\tau\in [1,\infty]$ if $2<\frac{N}{\theta}< p<\infty$ and $\mathring{W}^{\theta,\mathbf{S}}_{p,p}(\hat{\Omega}, \Omega)\to L^{\tau}(\Omega)$ for every $\tau\in [1,\infty [$ if $\frac{N}{\theta}=p$. Therefore, a simple modification of the proofs of the previous lemmas and the results obtained in Step 1  will easily give the corresponding estimate \eqref{est-65-3} with $\alpha$, $\beta$ and $\gamma$ given by \eqref{abg2-b}.\\

{\bf Part II}:  we remove the requirement that $u_0$, $v_0\in L^\infty(\Omega)$. Indeed, let $u_0$, $v_0\in L^q(\Omega)$ and let $(u_{0,n})$, $(v_{0,n})$ be two sequences in $L^\infty(\Omega)$ that converge to $u_0$ and $v_0$ in $L^q(\Omega)$, respectively,  as $n\to\infty$. Let $u_n:=S^{D,\nu}(\cdot)u_{0,n}$ and $v_n:=S^{D,\nu}(\cdot)v_{0,n}$. It follows from the estimate \eqref{est-65-3} that for every $t>0$,
\begin{equation}\label{D12}
 \|(u_n -v_n)(t,\cdot)\|_\infty \le C\frac{|\Omega|^\alpha}{t^{\beta}}\|u_{0,n}-v_{0,n}\|_q^\gamma.
 \end{equation}
As $n\to\infty$, the right hand side of \eqref{D12} converges to
$\displaystyle C\frac{|\Omega|^\alpha}{t^{\beta}}\|u_0-v_0\|_q^\gamma$.
Thus, after taking subsequences if necessary,  $(u_n(t,\cdot))$  and $(v_n(t,\cdot))$ converge weakly$^\star$ in $L^\infty(\Omega)$ to some functions 
$\tilde u(t,\cdot)$ and $\tilde v(t,\cdot)$, respectively,  as $n\to\infty$. Using the weak$^\star$ lower semicontinuity of the $L^\infty$-norm we get the bound
\[ \|(\tilde u-\tilde v)(t,\cdot)\|_\infty \le C\frac{|\Omega|^\alpha}{t^{\beta}}\|u_{0}-v_0\|_q^\gamma.\]
The submarkovian property of the semigroup $S^{D,\nu}$ also implies that
\[\|(u_n -u)(t,\cdot)\|_q\le \|u_{0,n}-u_0\|_q \hbox{ and } \|(v_n-v)(t,\cdot)\|_q\le \|v_{0,n}-v_0\|_q.\]
This shows that $u_n(t,\cdot)\to u(t,\cdot)$ and $v_n(t,\cdot)\to v(t,\cdot)$ in $L^q(\Omega)$, as $n\to\infty$. Hence, $\tilde u(t,\cdot)=u(t,\cdot)$, $\tilde v(t,\cdot)=v(t,\cdot)$. 
We have shown \eqref{est-65-3} with $\alpha$, $\beta$ and $\gamma$ given by \eqref{abg2}.
\end{proof}

We have the following result as a direct consequence of Theorem \ref{theo-ultra}.

\begin{corollary}\label{cor-ultra}
 Take the assumptions of this section, and let $S^{D,\nu}$ be the submarkovian semigroup generated by the negative $j$-subgradient $-\partial_j\mathcal E_{B_D}$.
 Let $q\in[1,\infty]$ and $s\in [1,\infty[$. Then, there is a constant $C=C(N,p,q,s,\theta )>0$ such that for every $t>0$ and every $u_0$, $v_0\in L^2(\Omega)\cap L^{q}(\Omega)$ we have the estimate
 \begin{equation}\label{Cor-est-65-3}
 \|S^{D,\nu}(t)u_0-S^{D,\nu}(t)v_0\|_s \le C\frac{|\Omega|^{\alpha+\frac 1s}}{t^{\beta}}\|u_0-v_0\|_{q}^{\gamma},
 \end{equation}
where $(\alpha,\beta,\gamma)$ are given in \eqref{abg2} or \eqref{abg2-b} depending on the range of $p$.
\end{corollary}

\begin{proof}
Consider the estimate \eqref{est-65-3} with $q\in [1,\infty]$ and $(\alpha,\beta,\gamma)$ given by \eqref{abg2} or \eqref{abg2-b}
 depending on the range of $p$. Since $s\in [1,\infty[$, we have that
 \begin{equation}\label{ESS1}
 \|S^{D,\nu}(t)u_0-S^{D,\nu}(t)v_0\|_s\leq \|S^{D,\nu}(t)u_0-S^{D,\nu}(t)v_0\|_\infty|\Omega|^{\frac 1s}.
 \end{equation}
 Combining \eqref{est-65-3} and \eqref{ESS1} we get 
 \begin{align*}
 \|S^{D,\nu}(t)u_0-S^{D,\nu}(t)v_0\|_s \le C\frac{|\Omega|^\alpha|\Omega|^{\frac 1s}}{t^{\beta}}\|u_0-v_0\|_q^\gamma= C\frac{|\Omega|^{\alpha+\frac{1}{s}}}{t^{\beta}}\|u_0-v_0\|_q^\gamma,
 \end{align*}
 and we have shown \eqref{Cor-est-65-3}.
 \end{proof}

\begin{remark}
{\em 
We observe that for every $q\in [1,\infty[$, 
\begin{align*}
\lim_{p\to 2}\alpha(p,q)=0,\; \;\lim_{p\to 2}\beta(p,q)=\frac{N}{2\theta q}\;\mbox{ and } \lim_{p\to 2}\gamma(p,q)=1,
\end{align*}
so that the linear submarkovian semigroup $S^{D,\nu}$ ($p=2$) satisfies the estimate
 \begin{equation}\label{p=2}
 \|S^{D,\nu}(t)u_0\|_s \le C\frac{|\Omega|^{\frac 1s}}{t^{\frac{N}{2\theta q}}}\|u_0\|_q,
 \end{equation}
for every $s\in [1,\infty[$ and $u_0\in L^2(\Omega)\cap L^q(\Omega)$. This is consistent with the well-known results for the linear case contained in Davies \cite{Dav} and the references therein.
}
\end{remark}

\subsection{Neumann exterior conditions}\label{Sec-Ultra-Neum}

In this section we take the situation of Example \ref{ex.3.1}. We consider Neumann exterior conditions, that is, we consider the energy $\E : W^{\theta,{\mathbf S}}_{p,2} (\hat{\Omega} ,\Omega)\to [0,\infty[$ given by
\begin{equation}\label{func-N}
\E (\hat{u} ) = \frac{1}{2p}\iint_{{\mathbf S}} k(x,y)|\hat u(x)-\hat u(y)|^p \;\ud y\;\ud x .
\end{equation}
We let $j$ be the restriction operator from ${W}^{\theta,{\mathbf S}}_{p,2} (\hat{\Omega} , \Omega )$ into $L^2 (\Omega )$, and we denote by $S^N$ the semigroup on $L^2 (\Omega )$ generated by the negative $j$-subgradient $\partial_j\E$. 

It follows from the G\^ateaux differentiability of $\mathcal E$  on the space $W^{\theta,{\mathbf S}}_{p,2} (\hat{\Omega} ,\Omega )$ (Lemma \ref{lem.identification}) and the definition of the $j$-subgradient $\partial_j\E$ that for every $(u,f)\in\partial_j\E$ there exists an elliptic extension $\hat{u}\in W^{\theta,{\mathbf S}}_{p,2} (\hat{\Omega} ,\Omega )$ of $u$, such that for every $\hat{v}\in W^{\theta,{\mathbf S}}_{p,2}(\hat{\Omega},\Omega)$,
\begin{align}\label{INt-parts-N}
\int_{\Omega} f \hat{v} \;\ud x = \iint_{{\mathbf S}} k(x,y)|\hat{u}(x)-\hat{u}(y)|^{p-2}(\hat{u}(x)-\hat{u}(y))(\hat{v}(x)-\hat{v} (y)\; \ud y \; \ud x.
\end{align}

We assume that $\Omega\subseteq\R^N$ has the {\em $W_p^{\theta}$-extension property} in the sense that there is a bounded, linear extension operator $E: W_p^{\theta} (\Omega) \to W_p^{\theta}(\mathbb R^N)$ such that $Eu|_{\Omega}=u$ for every $u\in W_p^{\theta} (\Omega)$. In that case, there is a constant $C=C(N,\theta,p,\Omega)>0$ such that
\begin{equation}\label{S0}
\|Eu\|_{W_p^{\theta}(\mathbb R^N)}\le C\|u\|_{W_p^{\theta}(\Omega)}.
\end{equation}
Moreover, 
there is a constant $C=C(N,\theta,p,\Omega)>0$ such that for every $u\in W_p^{\theta}(\Omega)$, 
\begin{equation}\label{C-Eb}
\|u\|_{q}\le C\|u\|_{W_p^{\theta}(\Omega)},
\end{equation}
with $q$ given in \eqref{eq-q}, that is, we have the continuous 
embedding $W_p^{\theta}(\Omega)\hookrightarrow L^q(\Omega)$.

We have the following Sobolev type embedding. 

\begin{lemma}\label{lem-SE}
Let $1\le q\le \frac{Np}{N-\theta p}$ if $N>\theta p$, $q\in [1,\infty[$ if $N=\theta p$, and $q\in [1,\infty]$ if $N<\theta p$. 
Then, the following assertions hold.
\begin{enumerate}
\item There is a constant $C=C(N,\theta,p,\Omega)>0$ such that for every $\hat{u}\in W_{p,p}^{\theta,{\mathbf S}}(\hat{\Omega},\Omega)$,
\begin{equation}\label{Sob-Ne}
\|\hat{u}|_{\Omega}\|_{q}\le C\|\hat{u}\|_{W_{p,p}^{\theta,{\mathbf S}}(\hat{\Omega},\Omega)}.
\end{equation}
\item  There is a constant $C=C(N,\theta,p,\Omega)>0$ such that for every $\varepsilon>0$ and $\hat{u}\in W_{p,p}^{\theta ,{\mathbf S}}(\hat{\Omega},\Omega)$,
\begin{equation}\label{sob-St}
\|\hat{u}|_{\Omega}\|_{q}^p\le C\left(\varepsilon\int_{\Omega}|\hat{u}|^p\;\ud x+\iint_{\mathbf S} \frac{|\hat{u}(x)-\hat{u}(y)|^p}{|x-y|^{N+\theta p}}\;\ud y\;\ud x \right).
\end{equation}
\end{enumerate}
\end{lemma}

\begin{proof}
Since by assumption $\Omega$ has the extension property, it follows from \eqref{S0} and \eqref{C-Eb} that there is a constant $C=C(N,\theta,p,\Omega)>0$ such that for every $\hat{u}\in W_{p,p}^{\theta ,{\mathbf S}}(\hat{\Omega},\Omega)$,
\begin{equation*}
    \|\hat{u}|_{\Omega}\|_{q}  \le C \|\hat{u}|_{\Omega}\|_{W^{\theta}_{p}(\Omega)}  \le C \|\hat{u}\|_{W_{p,p}^{\theta , {\mathbf S}}(\hat{\Omega},\Omega)}
\end{equation*}
and we have shown \eqref{Sob-Ne}. Finally, \eqref{sob-St} is a direct consequence of \eqref{Sob-Ne}. 
\end{proof}

Here also we recall that $p\in [2,\infty[$.

\begin{lemma}[\bf Logarithmic Sobolev inequality for $W_{p,p}^{\theta,{\mathbf S}}(\hat{\Omega},\Omega)$] 
  Let $\hat{u}\in W_{p,p}^{\theta,{\mathbf S}}(\hat{\Omega},\Omega)$ be such that $\|\hat{u}\|_{p}=1$. Then, the following assertions hold.
  \begin{enumerate}
      \item For every $\varepsilon>0$, we have that
      \begin{equation}\label{LSI-1}
        \int_{\Omega} |\hat{u}|^p\log(|\hat{u}|)\;\ud x \le C(N,p,\theta) \left[-\log(\varepsilon)+\varepsilon C \iint_{\mathbf S} k(x,y)|\hat{u}(x)-\hat{u}(y)|^p \;\ud y\;\ud x +C\varepsilon\right],
      \end{equation}
where $C>0$ is the constant appearing in \eqref{Sob-Ne}.

\item For every $\varepsilon>0$ and $\varepsilon_1>0$, we have that
      \begin{equation}\label{LSI-2}
        \int_{\Omega} |\hat{u}|^p\log(|\hat{u}|)\;\ud x \le C(N,p,\theta) \left[-\log(\varepsilon_1)+\varepsilon_1 C\iint_{\mathbf S}k(x,y) |\hat{u}(x)-\hat{u}(y)|^p \;\ud y\;\ud x +\varepsilon_1C\varepsilon\right],
      \end{equation}
where $C>0$ is the constant appearing in \eqref{sob-St}.
  \end{enumerate}
In \eqref{LSI-1} and \eqref{LSI-2} the constant $C(N,p,\theta)$ is given by
\begin{equation*}
C(N,p,\theta)=
\begin{cases}
\frac{N}{\theta p^2}=\frac{N}{\theta p}\frac{1}{p}\;&\mbox{ if } 2\le  p<\frac{N}{\theta}\\
\frac{\theta}{N}=\frac{\theta p}{N}\frac{1}{p}&\mbox{ if } \frac{N}{\theta}\le  p<\infty.
\end{cases}
\end{equation*}
\end{lemma}

\begin{proof}
The proof follows the lines of the proof of Lemma \ref{LSI} by using Lemma \ref{lem-SE}. 
\end{proof}

The main concern of this subsection is to prove an $L^q-L^\infty$- H\"older-continuity property of the semigroup $S^N$ generated by the negative $j$-subgradient $-\partial_j\mathcal E$, where $\mathcal E$ is given by \eqref{func-N}. More precisely, we have the following result.

\begin{theorem}\label{theo-ultra-NEC}
 Under the assumptions of this subsection, for every $q\in [1,\infty[$, there are two constants $C_1=C(N,p,q,s,\theta)$ and $C_2=C(N,p,q,s,\theta)>0$ such that for every $t>0$, and $u_0$, $v_0\in L^2(\Omega)\cap L^q(\Omega)$,
 \begin{equation}\label{est-65-NE}
 \|S^N(t)u_0-S^N(t)v_0\|_\infty \le C_1\frac{|\Omega|^\alpha e^{C_2t}}{t^{\beta}}\|u_0-v_0\|_q^\gamma,
 \end{equation}
where $(\alpha,\beta,\gamma)$ are given by \eqref{abg2} or \eqref{abg2-b} depending on the range of $p$.

 If $p=2$, then the associated linear semigroup $S^N$ satisfies for all $t>0$ the estimate
 \begin{equation}\label{p=2-NE}
 \|S^N(t)u_0\|_\infty \le \frac{C_1e^{C_2t}}{t^{\frac{N}{2\theta q}}}\|u_0\|_q.
 \end{equation}
 \end{theorem}

 \begin{proof}
 The proof follows the lines of the proof of Theorem \ref{theo-ultra} with the necessary series of lemmas proved in Subsection \ref{sec-UlD}.  Let $q\in [1,\infty[$ and let $u_0, v_0\in L^\infty(\Omega)$. 
 We proceed in several steps by just giving the main differences with the proofs of the lemmas in Subsection \ref{sec-UlD} and Theorem \ref{theo-ultra}. Let $u(t):=S^N(t)u_0$ and $v(t):=S^N(t)v_0$.\\

 {\bf Step 1}:  for any $r\in [2,\infty[$ consider the function $f_r: ]0,\infty [\to ]0,\infty [$ defined by
\[
f_r(t):=\int_{\Omega}|u(t,x)-v(t,x)|^r\;\ud x = \| (u-v)(t,\cdot )\|_r^r .
\]
As in the proof of Lemma \ref{lem-deri-2} we have that $f_r$ is differentiable for a.e. $t\ge 0$, and
\begin{equation}\label{deri-est-2-NE}
\frac{\ud}{\ud t} f_r(t)\le -rC\iint_{\mathbf S} k(x,y)\left|(\hat{u}-\hat{v})(t,x)-(\hat{u}-\hat{v})(t,y)\right|^{p+r-2}\;\ud y\;\ud x,
\end{equation}
where $C\in ]0,1]$ is the constant appearing in \eqref{prod}, and  $\hat{u}(t,\cdot )$ and $\hat{v} (t,\cdot )$ are elliptic extensions of $u(t,\cdot)$ and $v(t,\cdot )$, respectively.\\

{\bf Step 2}:  let $r:[0,\infty[\to [2,\infty[$ be any increasing differentiable function. Proceeding exactly as in the proof of Lemma \ref{lem-deri-log-22} we obtain that for a.e. $t\ge 0$, 
\begin{align}\label{der-log-2-NE}
\frac{\ud}{\ud t}\log\|(u-v)(t,\cdot)\|_{r(t)}  \le& \frac{r'(t)}{r(t)}\int_{\Omega}\frac{|(u-v)(t,x)|^{r(t)}}{\|(u-v)(t,\cdot)\|_{r(t)}^{r(t)}}\log\left(\frac{|(u-v)(t,x)|}{\|(u-v)(t,\cdot)\|_{r(t)}}\right)\;\ud x\notag\\
&-\frac{C}{\|(u-v)(t,\cdot)\|_{r(t)}^{r(t)}} \iint_{\mathbf S} k(x,y)|(\hat{u}-\hat{v})(t,x)-(\hat{u}-\hat{v})(t,y)|^{p+r(t)-2}\;\ud y\;\ud x,
\end{align}
where $C\in ]0,1]$ is the constant appearing in \eqref{deri-est-2-NE}.\\

{\bf Step 3}:  let $r$ be as in Step 2. Proceeding as in the proof of Lemma \ref{lem-deri-log-23} by using the logarithmic Sobolev inequality \eqref{LSI-2}, we get that there is a constant $C>0$ such that for every $\varepsilon>0$, $\varepsilon_1>0$ and for a.e. $t\ge 0$, 
\begin{align}\label{AA02-NE}
&\frac{\ud}{\ud t}\log\|(u-v)(t,\cdot)\|_{r(t)}  \le \frac{r'(t)}{r(t)}\int_{\Omega}\frac{|(u-v)(t,x)|^{r(t)}}{\|(u-v)(t,\cdot)\|_{r(t)}^{r(t)}}\log\left(\frac{|(u-v)(t,x)|}{\|(u-v)(t,\cdot)\|_{r(t)}}\right)\;\ud x\notag\\
&-C_{r(t),p}\frac{\|(u-v)(t,\cdot)\|_{r(t)+p-2}^{r(t)+p-2}}{\|(u-v)(t,\cdot)\|_{r(t)}^{r(t)}}\frac{sp^2(r(t)+p-2)}{pNC\varepsilon}\int_{\Omega}\frac{|(u-v)(t,x)|^{r(t)+p-2}}{\|(u-v)(t,\cdot)\|_{r(t)}^{r(t)}}\log\left(\frac{|(u-v)(t,x)|}{\|(u-v)(t,\cdot)\|_{r(t)+p-2}}\right)\;\ud x \notag\\
&-C_{r(t),p}\frac{\|(u-v)(t,\cdot)\|_{r(t)+p-2}^{r(t)+p-2}}{\|(u-v)(t,\cdot)\|_{r(t)}^{r(t)}}\frac{\log\varepsilon_1}{C\varepsilon_1}+C_{r(t),p}\frac{\|(u-v)(t,\cdot)\|_{r(t)+p-2}^{r(t)+p-2}}{\|(u-v)(t,\cdot)\|_{r(t)}^{r(t)}}\varepsilon.
\end{align}

{\bf Step 4}:  let $r$ be as in Step 2. Using Step 3 and proceeding as in the proof of Lemma \ref{lem-diff-ine-2} by letting
\begin{align*}
\varepsilon_1:=C_{r(t),p}\frac{\|(u-v)(t,\cdot)\|_{r(t)+p-2}^{r(t)+p-2}}{\|(u-v)(t,\cdot)\|_{r(t)}^{r(t)}}\frac{\theta pr(t)(r(t)+p-2)C}{Nr'(t)}\; 
\mbox{ and  }
\varepsilon:=\frac{\|(u-v)(t,\cdot)\|_{r(t)+p-2}^{r(t)+p-2}}{\|(u-v)(t,\cdot)\|_{r(t)}^{r(t)}},
\end{align*}
we obtain that for a.e. $t\ge 0$, 
\begin{equation}\label{eq-diff-ine2-NE}
    \frac{\ud}{\ud t}\log\|(u-v)(t,\cdot)\|_{r(t)}\le -\mathcal P(t)\log\|(u-v)(t,\cdot)\|_{r(t)} -\mathcal Q(t),
\end{equation}
where here
\begin{align}\label{eqP}
    \mathcal P(t):=\frac{r'(t)}{r(t)}\frac{N}{\theta p}\frac{p-2}{r(t)+p-2}
\end{align}
and 
\begin{align}\label{eqQ}
    \mathcal Q(t):=&\frac{N}{\theta p}\frac{r'(t)}{r(t)}\frac{1}{r(t)+p-2}\log\left(\frac{r(t)}{r'(t)}\frac{C_{r(t),p}\theta pC(r(t)+p-2)}{N}\right)\notag\\
  &-\frac{r'(t)}{r(t)}\left(\frac{N}{\theta p}-1\right)\frac{p-2}{r(t)+p-2}\log|\Omega|+C_{r(t),p}.
\end{align}

{\bf Step 5}:  let $r$ be as in Step 2. Let $Y(t):=\log\|(u-v)(t,\cdot)\|_{r(t)}$. It follows from Step 4 that $Y(t)$ satisfies the ordinary differential inequality
\begin{equation}\label{ODE-INE2-NE}
    Y'(t)+\mathcal P(t)Y(t)+\mathcal Q(t)\le 0,
\end{equation}
where $\mathcal P$ and $\mathcal Q$ are given in \eqref{eqP} and \eqref{eqQ}, respectively.

The rest of the proof follows the proof of Theorem \ref{theo-ultra} Parts I and II by considering the ordinary differential inequality  \eqref{ODE-INE2-NE}.
We omit the full detail for the sake of brevity.
 \end{proof}

\section{Examples and Remarks}\label{section-ER}

In this section we give additional examples of nonlocal and fractional operators that enter in our framework. All the results obtained in the previous sections are valid for these operators.

\subsection{Fractional powers of linear elliptic operators in domains of $\mathbb R^N$}

In this subsection we let $\theta \in]0,1[$, $p\in]1,\infty[$ and $\Omega\subseteq\mathbb R^N$ an arbitrary open set.
We also let
\begin{equation}\label{space-ws0}
\ringring{W}^{\theta,{\mathbf S}}_{p,r}(\hat\Omega,\Omega) := \overline{C_c^\infty (\Omega)}^{{W_{p,p}^{\theta,{\mathbf S}}(\hat\Omega,\Omega)}}.
\end{equation}
Notice that the spaces $\mathring{W}^{\theta,{\mathbf S}}_{p,r}(\hat\Omega,\Omega)$ and $\ringring{W}^{\theta,{\mathbf S}}_{p,r} (\hat\Omega,\Omega)$ do not always coincide. We refer to Claus \& Warma \cite{ClWa20} and Grisvard \cite{Gris} for a complete discussion on this topic.

\begin{example}[\bf The fractional Laplace operator in $\R^N$]\label{flrn}
A prototype of a nonlocal operator and subgradient is a fractional power of the negative Laplace operator on $\R^N$. The Laplace operator on $L^2 (\R^N )$ is the operator $\Delta =\sum_{i=1}^N D_{ii}$ with domain $W^{2}_{2} (\R^N)$. The negative Laplace operator is a selfadjoint, nonnegative operator. 
Let $(e^{t\Delta})_{t\ge 0}$ be the semigroup on $L^2(\R^N)$ generated by the operator $\Delta$, that is, the Gaussian semigroup. For every $\theta\in ]0,1[$ and for every $u\in W^{2\theta}_{2} (\R^N )$ (see e.g. Di Nezza, Palatucci \& Valdinoci \cite{DNPaVa12} for a definition of $ W^{2\theta}_{2} (\R^N )$ using the Fourier transform), we define the fractional power $(-\Delta)^\theta$ as follows:
\begin{equation}\label{Lap-SG}
    (-\Delta)^{\theta}u=\frac{-\theta}{\Gamma(1-\theta)}\int_0^\infty\left(e^{t\Delta}u-u\right)\frac{\ud t}{t^{1+\theta}} .
\end{equation}
One can show that $(-\Delta)^{\theta}$ defined in \eqref{Lap-SG} is given exactly by
\begin{equation}\label{Lap-SG-2}
(-\Delta)^\theta u (x) = C_{N,\theta}\text{P.V.} \int_{\R^N} \frac{u(x)-u(y)}{|x-y|^{N+2\theta}} \; \ud y \qquad (x\in\R^N),
\end{equation}
where as before P.V. stands for the Cauchy principal value, and 
\begin{equation} \label{eq.cn}
C_{N,\theta}:=\frac{\theta 2^{2\theta}\Gamma\left(\frac{N+2\theta}{2}\right)}{\pi^{\frac N2}\Gamma(1-\theta)}=\left(\int_{\mathbb R^N}\frac{1-\cos(\xi_1)}{|\xi|^{N+2\theta}}\;\ud\xi\right)^{-1}.
\end{equation}
In fact, let $G$ be the kernel of the operator $e^{t\Delta}$ which is given by 
$$
G(t,x,y)=\frac{1}{(4\pi t)^{N/2}}e^{-\frac{|x-y|^2}{4t}},\; \;(t>0,\; x,y\in\R^N).
$$
Using \eqref{Lap-SG} we have that for every $u\in W_2^{2\theta}(\R^N )$ and $v\in W_2^{\theta}(\R^N )$,
\begin{align}\label{EP-R}
& \int_{\R^N} (-\Delta)^{\theta}u (x) v(x) \;\ud x \notag\\
=&\frac{-\theta}{\Gamma(1-\theta)}\int_0^\infty\int_{\R^N}\left[\int_{\R^N} G(t,x,y)u(y)v(x)\;\ud y-u(x)v(x)\right]\;\ud x\frac{\ud t}{t^{1+\theta}}\notag\\
=&\frac{-\theta}{\Gamma(1-\theta)}\int_0^\infty\int_{\R^N}\left[\int_{\R^N} G(t,x,y)(u(x)-u(y))v(x)\;\ud y-u(x)v(x )\left(\int_{\R^N}G(t,x,y)\;\ud y-1(x)\right)\right]\;\ud x\frac{\ud t}{t^{1+\theta}}\notag\\
=&\frac{-\theta}{\Gamma(1-\theta)}\int_0^\infty\int_{\R^N}\int_{\R^N} G(t,x,y)(u(x)-u(y))v(x)\;\ud y\;\ud x\frac{\ud t}{t^{1+\theta}}\notag\\
&+\int_0^{\infty}\int_{\R^N}u(x)v(x)\left(e^{t\Delta}1(x)-1(x)\right)\;\ud x\frac{\ud t}{t^{1+\theta}}.
\end{align}
Since $e^{t\Delta}1(x)=1$, it follows from \eqref{EP-R} that
\begin{align*}
\int_{\R^N} (-\Delta)^{\theta}u(x) v(x) \; \ud x = & \frac{-\theta}{\Gamma(1-\theta)}\int_0^\infty\int_{\R^N}\int_{\R^N} G(t,x,y)(u(x)-u(y))v(x)\;\ud y\;\ud x\frac{\ud t}{t^{1+\theta}}\\
=&\frac{-\theta}{2\Gamma(1-\theta)}\int_{\R^N}\int_{\R^N} \int_0^\infty G(t,x,y)\frac{\ud t}{t^{1+\theta}} (u(x)-u(y))v(x)\;\ud y\;\ud x\\ 
=&\frac{C_{N,\theta}}{2}\int_{\R^N}\int_{\R^N}\frac{(u(x)-u(y)(v(x)-v(y))}{|x-y|^{N+2\theta}}\;\ud y\;\ud x,\\
=&\int_{\R^N}\left(\frac{C_{N,\theta}}{2}\mbox{P.V.}\int_{\R^N}\frac{(u(x)-u(y)}{|x-y|^{N+2\theta}}\;\ud y\right)v(x)\;\ud x,
\end{align*}
where we have used that (by using a simple change of variables)
\begin{align*}
\frac{\theta}{\Gamma(1-\theta)}\int_0^\infty G(t,x,y)\;\frac{\ud t}{t^{1+\theta}}
=&\frac{\theta}{\Gamma(1-\theta)}\frac{1}{(4\pi)^{N/2}}\int_0^\infty e^{-\frac{|x-y|^2}{4t}}\frac{1}{t^{1+\theta+N/2}}\;\ud t\\
=&\frac{\theta}{\Gamma(1-\theta)}\frac{1}{(4\pi)^{N/2}}\frac{4^{1+\theta+N/2}}{4}\frac{1}{|x-y|^{N+2\theta}}\int_0^\infty e^{-\tau}\tau^{\theta+N/2-1}\;\ud\tau\\
=&\frac{\theta 2^{2\theta}\Gamma\left(\frac{N+2\theta}{2}\right)}{\pi^{\frac N2}\Gamma(1-\theta)} \frac{1}{|x-y|^{N+2\theta}}\\
=&\frac{C_{N,\theta}}{|x-y|^{N+2\theta}}.
\end{align*}

The fractional power of the negative Laplace operator is the subgradient on $L^2 (\R^N )$ of the quadratic form $\E_\theta : L^2 (\R^N ) \to [0,\infty ]$ given by 
\begin{equation}\label{Ener-E}
\E_\theta (u) = \begin{cases} 
 \displaystyle \frac{C_{N,\theta}}{4} \int_{\R^N} \int_{\R^N} \frac{(u(x) -u(y))^2}{|x-y|^{N+2\theta}} \; \ud y \, \ud x & \text{if } u\in W^{\theta}_2 (\R^N )  , \\[2mm]
 +\infty & \text{otherwise} ,
 \end{cases}
\end{equation}
so that it is a special case of Example \ref{ex.3.1} by taking $\Omega = \hat{\Omega} = \R^N$ and $\mathbf{S} = \R^N\times\R^N$. The kernel $k$ satisfies the condition \eqref{cond.k} for the given $\theta$ and for $p=2$.   
It is known (see e.g. Keyantuo, Seoanes \& Warma \cite{KSW} and the references therein) that the linear $C_0$-semigroup $S$ generated by $-(-\Delta)^\theta$ is submarkovian and analytic, and that it extrapolates to a positive, analytic contraction semigroup on $L^q(\Omega)$ for every $q\in [1, \infty[$. Due to the continuous embedding $W^{\theta}_2 (\R^N)\hookrightarrow L^{\frac{2N}{N-2\theta}}(\R^N)$, the semigroup $S$ is ultracontractive, that is, there are constants $C_1,C_2\ge 0$ such that for every $q\in[1,\infty[$, $u_0\in L^2(\Omega)\cap L^q(\Omega)$ and $t>0$, it satisfies the estimate 
\begin{equation*}
    \|S(t)u_0\|_\infty\le \frac{C_1e^{C_2t}}{t^{\frac{N}{2\theta q}}}\|u_0\|_q.
\end{equation*}
It has been shown in Blumenthal \& Getoor\cite{BlGe60}, Chen, Kim \& Song \cite{ChKiSo10}, and Chen \& Kumagai \cite{ChKu08} that the heat kernel $P_\theta (t,x,y)$ associated with this semigroup is given by
\begin{equation}\label{kernel-sg}
    P_\theta(t,x,y)=t^{-\frac{N}{2\theta}}\left(1+|x-y|t^{-\frac{1}{2\theta}}\right)^{-N-2\theta} \quad (t>0, \;x,y\in\R^N).
\end{equation}

For more information on the fractional Laplace operator on $\R^N$, we refer to Caffarelli \& Silvestre \cite{CaSi07},  Di Nezza, Palatucci \& Valdinoci \cite{DNPaVa12} and the references therein.
\end{example}

\begin{example}[\bf Fractional powers of linear second order elliptic operators on domains]
Fractional powers of more general linear second order elliptic operators equipped with Neumann, Dirichlet or Robin boundary conditions also fit into the setting in this article or, more precisely, into the setting of Examples \ref{ex.3.1} and \ref{ex.6.2}. Let $\Omega\subseteq\R^N$ be an arbitrary open set. Take coefficients $a_{ij}\in L^\infty(\Omega)$ ($1\le i,j\le N$), and assume that $a_{ij}=a_{ji}$ for all $1\le i,j\le N$. Assume further that the $a_{ij}$ are {\em uniformly elliptic} in the sense that there exists $\alpha >0$ such that 
\[
\sum_{i,j=1}^N a_{ij}(x)\xi_i\xi_j\ge \alpha |\xi|^2\qquad (x,\in\Omega, \, \xi\in\R^N) .
\]
\begin{enumerate}
\item[(1)] {\bf Neumann boundary conditions}: we define the quadratic form $\E : L^2 (\Omega ) \to [0,\infty]$ by
\[
\E_N (u) = \begin{cases} 
\displaystyle \frac12 \int_\Omega \sum_{i,j=1}^N a_{ij} D_ju D_iu & \text{if } u\in W^{1}_2 (\Omega ) , \\[2mm]
 +\infty & \text{otherwise} ,
 \end{cases}
\]
where $W^{1}_2 (\Omega ):=\{u\in L^2(\Omega):\; |\nabla u|\in L^2(\Omega)\}$ is the classical first order Sobolev space. Due to the boundedness and ellipticity of the coefficients, the quadratic form is proper, lower semicontinuous and convex. The subgradient $A_N = \partial\E_N$ is a linear, selfadjoint, positive semidefinite operator on $L^2 (\Omega )$. This operator is a realization on $L^2 (\Omega )$ of the second order elliptic operator formally given by
\begin{equation}\label{SOEO}
A u= - \sum_{i,j=1}^N D_i\left(a_{ij}(x)D_ju\right) \text{ in } \Omega ,
\end{equation}
and complemented by Neumann boundary conditions on $\partial\Omega$. Let $(e^{-tA_N})_{t\ge 0}$ be the semigroup on $L^2(\Omega)$ generated by $-A_N$.
One may define fractional powers of $A_N$ as in \eqref{Lap-SG}, that is, for $u\in \dom{A_N}$ and $\theta\in ]0,1[$, we let
\begin{equation}\label{SG-AN}
    (A_N)^{\theta}u=\frac{-\theta}{\Gamma(1-\theta)}\int_0^\infty\left(e^{-tA_N}u-u\right)\frac{\ud t}{t^{1+\theta}}.
\end{equation}
One may also define the fractional powers of $A_N$ for example, by using the holomorphic functional calculus, the Hille-Phillips functional calculus, or the spectral theorem. By the spectral theorem, $A_N$ is unitarily equivalent to a multiplication operator $M$ on $L^2 (X)$, where $X$ is an abstract measure space and $Mg := mg$ for some measurable function $m:X\to [0,\infty [$ (see e.g. Reed \& Simon \cite[Theorem VIII.4, p. 260]{ReSi80I}). This means that there exists a unitary operator $U : L^2 (\Omega ) \to L^2 (X)$ such that $U$ is a bijection from $\dom{A_N}$ onto $\dom{M}$ and $A_N = U^*MU$. For a given $\theta\in ]0,1[$ one simply defines the fractional power $m^\theta$ of the measurable function $m$, the associated multiplication operator $M^\theta g := m^\theta g$ with maximal domain, and then the fractional power $(A_N)^\theta := U^* M^\theta U$ with maximal domain. 

In the special case when $A_N$ has compact resolvent or, equivalently, when the embedding of $W^{1}_2 (\Omega )$ into $L^2 (\Omega )$ is compact (for example, when $\Omega$ is bounded and has Lipschitz continuous boundary), then there exists an orthonormal basis $(e_n)$ of eigenvectors of $A_N$ associated with nonnegative eigenvalues $(\lambda_n)$. In this case, $L^2 (X) = \ell^2 (\N )$ and $U : L^2 (\Omega ) \to \ell^2 (\N )$ is the unitary operator associated with the orthonormal basis $(e_n)$, and
\[
A_N u = \sum_{n\in\N} \lambda_n \langle u,e_n\rangle_{L^2(\Omega )} e_n \text{ and } (A_N)^\theta u = \sum_{n\in\N} \lambda_n^\theta \langle u,e_n\rangle_{L^2(\Omega )} e_n .
\]

The operator $(A_N)^\theta$ is again a linear, selfadjoint, positive semidefinite operator. It is the subgradient of the proper, lower semicontinuous, convex energy $\E_\theta^N : L^2 (\Omega ) \to [0,\infty ]$ given by
\[
\E_\theta^N (u) = \begin{cases} 
\displaystyle \frac12 \| (A_N)^{\frac{\theta}{2}} u \|_{L^2 (\Omega )}^2 & \text{if } u\in \dom{(A_N)^{\frac{\theta}{2}}} , \\[2mm]
 +\infty & \text{otherwise} .
 \end{cases}
\]

It has been shown by Caffarelli \& Stinga \cite[Sections 1 and 2]{CaSt} that $\dom{(A_N)^{\frac{\theta}{2}}} = W^{\theta}_2 (\Omega )$. In addition, 
using \eqref{SG-AN} we get that for every  $u\in \dom{(A_N)^{\theta}}$ and $v\in W^{\theta}_2 (\Omega )$,
\begin{equation*}
\langle (A_N)^{\theta}u,v\rangle_{L^2(\Omega )}=\int_{\Omega}\int_{\Omega} k_\theta^N(x,y)(u(x)-u(y))(v(x)-v(y))\;\ud y\;\ud x,
\end{equation*}
so that for every $u\in W^{\theta}_2 (\Omega )$ one has 
\begin{equation} \label{int-rep-N}
\E_\theta^N (u) = \frac14 \int_{\Omega}\int_{\Omega} k_\theta^N(x,y) (u(x)-u(y))^2 \; \ud y \; \ud x ,
\end{equation}
where the kernel $k_\theta^N$ is given by 
\begin{align}\label{Ks}
0\le k_\theta^N (x,y) := \frac{\theta}{\Gamma(1-\theta)} \int_0^\infty\frac{G_\Omega^N(t,x,y)}{t^{1+\theta}}\;\ud t , \quad x,y\in\Omega,
\end{align}
and $G_{\Omega}^N$ is the associated heat kernel. For example, if $A_N$ has a compact resolvent, then 
\begin{align*}
G_\Omega^N (t,x,y)=\sum_{n=1}^\infty e^{-t\lambda_{n}} e_n(x) e_n (y),\;t>0,\; x,y\in\Omega ,
\end{align*}
and by \cite[Section 1]{CaSt},
$$0\le k_\theta^N (x,y)  \le \frac{C_{N,\theta}}{|x-y|^{N+2\theta}}, \quad x,y\in\Omega.$$
As a consequence, the fractional powers $(A_N)^\theta$ of the linear second order elliptic operators with Neumann boundary conditions are nonlocal operators which fit into the setting of this article; see especially Example \ref{ex.3.1} with $\hat{\Omega} = \Omega$ and $p=2$. It follows that the operators $-(A_N)^\theta$ generate linear, submarkovian, analytic $C_0$-semigroups $S_\theta^N$ which extrapolate to positive, analytic contraction semigroups on $L^q (\Omega )$ for every $q\in [1,\infty[$. Due to the continuous embedding 
\begin{equation}\label{inWe conj1}
W^{\theta}_2 (\Omega) \hookrightarrow
\begin{cases}
L^{\frac{2N}{N-2\theta}}(\Omega)\;&\mbox{ if }\; N>2\theta,\\[2mm]
L^q(\Omega),\;q\in[1,\infty [\;&\mbox{ if }\; N=2\theta,\\[2mm]
C^{0,\theta -\frac{N}{2}}(\bOm)\;&\mbox{ if }\; N<2\theta ,
\end{cases}
\end{equation}
the resulting linear semigroups are ultracontractive. More precisely, there are constants $C_1,C_2\ge 0$ such that for every $t>0$ and $u_0\in L^2(\Omega)\cap L^q(\Omega)$ ($q\in [1,\infty[$), 
 \begin{equation*}
 \|S_{\theta}^N(t)u_0\|_\infty \le \frac{C_1e^{C_2t}}{t^{\frac{N}{2\theta q}}}\|u_0\|_q.
 \end{equation*}

\item[(2)] {\bf Dirichlet boundary conditions}:
under the same conditions on the coefficients, one could also start with the quadratic form $\E_D : L^2 (\Omega ) \to [0,\infty]$ given by
\[
\E_D (u) = \begin{cases} 
\displaystyle \frac12 \int_\Omega \sum_{i,j=1}^N a_{ij} D_ju D_iu & \text{if } u\in \ringring{W}^{1}_2 (\Omega ) , \\[2mm]
 +\infty & \text{otherwise} .
 \end{cases}
\]
In this case, the subgradient $A_D = \partial\E_D$ is a realization of the second order elliptic operator \eqref{SOEO} equipped with Dirichlet boundary conditions. One may continue as above and define the fractional powers $(A_D)^\theta$ and the corresponding quadratic forms $\E_\theta^D$.  When $\Omega$ is bounded, one has 
\begin{align*}
\dom{(A_D)^\frac{\theta}{2}}  
=& \begin{cases}
\ringring{W}^{\theta}_2 (\Omega) & \mbox{ if } \theta\ne 1/2 , \\[2mm]
\displaystyle \dot{W}^{1/2}_2 (\Omega):=\left\{ u\in W^{1/2}_2 (\Omega) \st \int_{\Omega}\frac{|u(x)|^2}{\mbox{dist}(x,\partial\Omega)}\;\ud x<\infty\right\} &\mbox{ if } \theta =1/2. 
\end{cases} 
\end{align*}
Moreover, it has been shown in \cite[Theorem 2.3]{CaSt} that for every $u\in \dom{(A_D)^\frac{\theta}{2}}$, one has
\begin{align}\label{int-rep}
\E_\theta^D (u) = \int_{\Omega}\int_{\Omega} k_\theta^D(x,y) (u(x)-u(y))^2 \; \ud y \; \ud x +\int_{\Omega}\kappa_\theta (x)u(x)^2 \;\ud x,
\end{align}
where the kernel $k_\theta^D$ is defined similarly as in Part (1), with the heat kernel $G_\Omega^N$ replaced by the heat kernel $G_\Omega^D$ which is defined in exactly the same manner, replacing the orthonormal basis of $A_N$ (Neumann boundary conditions) by the orthonormal basis of $A_D$ (Dirichlet boundary conditions). If the coefficients $a_{ij}$ are H\"older continuous with exponent $\alpha\in ]0,1[$, then a two sided and precise estimate of $k_\theta^D(x,y)$ has been given in \cite[Theorem 2.4]{CaSt}. More precisely, there is a constant $c=c(N,\Omega)>0$ and $\eta\le 1\le \rho$ such that
\[
c^{-1}\min\left(1,\frac{\phi_0(x)\phi_0(y)}{|x-y|^{2\rho}} \right)\frac{1}{|x-y|^{N+2\theta}}\le k_\theta^D(x,y)\le c\min\left(1,\frac{\phi_0(x)\phi_0(y)}{|x-y|^{2\eta}}\right)\frac{1}{|x-y|^{N+2\theta}},
\]
where $\phi_0$ is the first nonnegative eigenfunction of the operator $A_D$. Here, for some constant $C=C(N,\alpha,\Omega)$ we have that
\[
C^{-1}\mbox{dist}(x,\partial\Omega)^{\rho}\le \phi_0(x)\le C\mbox{dist}(x,\partial\Omega)^{\eta}.
\]
In addition, 
\begin{align*}
0\le \kappa_\theta (x)=\frac{\theta}{\Gamma(1-\theta )} \int_0^\infty \Big(1-e^{-tA_D}1(x)\Big)\;\frac{\ud t}{t^{1+\theta}},\quad x\in\Omega .
\end{align*}
As a consequence, the fractional powers $(A_D)^\theta$ of the linear second order elliptic operators with Dirichlet boundary conditions are nonlocal operators which fit into the setting of this article with $\hat{\Omega} = \Omega$ and $p=2$. It follows that the operators $-(A_D)^\theta$ generate linear, submarkovian, analytic $C_0$-semigroups $S_\theta^D$ which extrapolate to positive, analytic contraction semigroups on $L^q (\Omega )$ for every $q\in [1,\infty[$.

Since the embedding \eqref{inWe conj1} holds with $W^{\theta}_2 (\Omega)$ replaced with
$\ringring{W}^{\theta}_2(\Omega)$, we also have that the resulting semigroups are ultracontractive. More precisely, there is a constant $C\ge 0$ such that for every $u_0\in L^2(\Omega)\cap L^q(\Omega)$ ($q\in [1,\infty[$) and $t>0$, one has
\begin{equation*}
    \|S_\theta^D(t)u_0\|_\infty\le \frac{C}{t^{\frac{N}{2\theta q}}}\|u_0\|_q.
\end{equation*}
\end{enumerate}
For more details on these two examples, we refer the interested reader to the paper \cite{CaSt}. We also refer to Haase \cite{Hs06}, Lunardi \cite{Lu95}, Arendt, ter Elst \& Warma \cite{ATW} for a theory of fractional powers of general sectorial operators.
\end{example}

\subsection{Other nonlocal operators in domains of $\mathbb R^N$}\label{sec.lnnod}

In the following series of examples we take up Examples \ref{ex.3.1}, \ref{ex.6.1} and \ref{ex.6.2}. We let $\Omega\subseteq \hat{\Omega} \subseteq\R^N$ be open subsets, ${\mathbf S}\subseteq\hat{\Omega}\times\hat{\Omega}$ satisfies the thickness condition, $k:\hat{\Omega}\times\hat{\Omega}\to [0,\infty]$ a symmetric  kernel satisfying condition \eqref{cond.k} for some $\theta\in ]0,1[$, and $p\in ]1,\infty[$. We let $\Phi : \hat{\Omega}\times\hat{\Omega}\to [0,\infty ]$ be given by
\[
\Phi (x,y,s) = \frac{1}{p} k(x,y) |s|^p \, 1_{\mathbf S} (x,y) .
\]
Then $\Phi$ satisfies the standard conditions. The associated space $W^{\Phi ,2} (\hat{\Omega},\Omega )$ is equal to the Banach space $W^{\theta,{\mathbf S}}_{p,2} (\hat{\Omega} , \Omega )$.
We further let $B:\hat{\Omega}\times\R\to [0,\infty]$ be a function satisfying the condition \eqref{cond.b}, and we let $\nu$ be a Borel measure on $\tilde{\Omega}$ which is absolutely continuous with respect to the capacity $\Cap_\E$, where $\E$ is the "free" energy $\E$ from \eqref{energy.e}. 

We consider the energy $\E_{B,\nu} : W^{\theta,{\mathbf S}}_{p,2}(\hat{\Omega},\Omega) \to [0,\infty ]$ given by 
\[
\E_{B,\nu} (\hat{u} ) = \iint_{{\mathbf S}} k(x,y) |\hat{u}(x)-\hat{u}(y)|^p \;\ud y\;\ud x + \int_{\tilde{\Omega}} B(x,\hat{u}(x))\; \ud\nu (x).
\]
Its effective domain is $W^{\theta,{\mathbf S}}_{p,2}(\hat{\Omega},\Omega) \cap L^B (\tilde{\Omega},\nu )$.

\begin{example}[\bf Regional "fractional" $p$-Laplace operators] \label{ex.9.3}
This is the case where $\Omega=\hat{\Omega}\subseteq\R^N$ is an arbitrary open set with boundary $\partial\Omega$, and ${\mathbf S} = \Omega\times\Omega$, $k(x,y) = |x-y|^{-(N+\theta p)}$ (but we can also deal with a more general kernel satisfying \eqref{cond.k}). As we want to concentrate on the boundary conditions, we assume that $B$ is only defined on $\partial\Omega$ (instead of all $\bar{\Omega}$). We recall from \eqref{eq.cpv} that $u\in W^\theta_p (\Omega)$ is a weak solution of the stationary problem
\begin{equation} \label{eq.RFPL}
\text{P.V.}\int_{\Omega}\frac{|u(x)-u(y)|^{p-2}(u(x)-u(y))}{|x-y|^{N+\theta p}}\;\ud y = f(x) 
\end{equation}
with $f\in L^1_{loc} (\Omega )$,  if for every $v\in {\mathcal D} := C^\infty_c (\Omega )$
\[
\frac12 \int_\Omega \int_\Omega \frac{1}{|x-y|^{N+\theta p}} (u(x)-u(y)) (v(x)-v(y)) \; \ud y \; \ud x = \int_\Omega f(x) v(x) \; \ud x ,
\]
and if $(u,f)\in\partial\E_{B,\nu}$, then $u$ is a weak solution of \eqref{eq.RFPL}. For more details on the operator $\partial\E_{B,\nu}$ we refer to Gal \& Warma \cite{Ga-Wa-CPDE}, Guan \& Ma \cite{GuMa1,GuMa2}, Warma \cite{War-NODEA,War18} and their references. This operator is known in the literature as the regional fractional $p$-Laplace operator. As we have observed in our general setting, for such operators we do not have exterior conditions. We just have boundary conditions as in the classical local case. In addition we shall see that the boundary operator is a local operator.
Here, we discuss Neumann and Robin boundary conditions, since Dirichlet boundary conditions are very classical. Before doing that we introduce an integration by parts formula associated with this nonlocal operator.

\begin{enumerate}
\item[(i)] {\bf The strong $(\alpha,p)$-normal derivative}: assume that $\Omega\subseteq\R^N$ has a $C^{1,1}$ boundary $\partial\Omega$. 
For $p\in ]1,\infty[$, $\alpha\in [0,2[$,  $u\in C^1(\Om)$ and $z\in\pOm$, we define the operator 
\begin{align}\label{B-OP}
\mathcal N_{\partial\Omega}^{2-\alpha,p} u(z):=&-\lim_{t\downarrow 0} \left|\frac{\ud u(z+t\vec{\nu}(z))}{\ud t}\right|^{p-2}\frac{\ud u(z+t\vec{\nu}(z))}{\ud t}t^{\alpha(p-1)}\notag\\
=&-\lim_{t\downarrow 0} \left|\frac{\ud u(z+t\vec{\nu}(z))}{\ud t}t^\alpha\right|^{p-2}\frac{\ud u(z+t\vec{\nu}(z))}{\ud t}t^{\alpha},
\end{align}
provided that the limit exists, where $\vec{\nu}$ denotes the unit outer normal vector. Observe that $\mathcal N_{\partial\Omega}^{2-\alpha,p}$ is a local operator. It is easy to see that 
if $u\in C^1(\overline{\Omega})$, then $\mathcal N_{\partial\Omega}^{2-\alpha,p} u(z)=0$ for every $z\in\partial\Omega$. Next, we introduce the class of functions for which this normal derivative is not always zero on $\partial\Omega$.

\begin{itemize}
\item Let $\alpha\in ]1,2[$ and 
\[
C_\alpha^1(\bOm):=\Big\{u \st u(x)=f(x)\text{\rm dist}(x,\partial\Omega)^{\alpha-1} +g(x) \text{ for every } x\in\Om,\;\text{ for some }\, f,g\in C^1(\bOm)\Big\}.
\]
Then, $\mathcal N_{\partial\Omega}^{2-\alpha,p} u(z)$ exists for every $z\in\pOm$ and is given by
\begin{equation*}
\mathcal N_{\partial\Omega}^{2-\alpha,p}u(z)=-(\alpha-1)|\alpha-1|^{p-2}|f(z)|^{p-2}f(z).
\end{equation*}
We refer to Warma \cite{War-NODEA} for the proof. The linear case $p=2$ is contained in Guan \cite{Gua}.

\item Let $\max\{\frac{1}{p}, \frac{p-1}{p}\}<\theta<1$ and $\alpha:=\frac{p\theta-1}{p-1}+1$. Consider the nonlocal problem given in \eqref{eq.RFPL}. Let also
\[
C_{\alpha}^2(\bOm):=\Big\{u \st u(x)=f(x)\text{\rm dist}(x,\partial\Omega)^{\alpha-1}+g(x) \text{ for every } x\in\Om, \text{ for some } f,g\in C^2(\bOm)\Big\}.
\]
Then, for every $u\in C_{\alpha}^2(\bOm)$, $f\in L^1 (\Omega )$
and $v\in W^{\theta}_p (\Om)\cap C(\overline{\Omega})$, we have that
\begin{align}\label{IPF-R}
\int_\Om fv \;\ud x=&\frac{C_{N,p,\theta}}{2}\int_\Om\int_\Om|u(x)-u(y)|^{p-2} \frac{(u(x)-u(y))(v(x)-v(y))}{|x-y|^{N+p\theta}}\;\ud y\;\ud x\notag\\
&-B_{N,p,\theta}\int_{\pOm}v\mathcal N_{\partial\Omega}^{2-\alpha,p}u\;d\sigma,
\end{align}
where 
$B_{N,p,\theta}$ is a normalizing constant depending only on $N,\theta,p$ and $\sigma$ is the usual Lebesgue surface measure on $\partial\Omega$.
\end{itemize}

\item[(ii)] {\bf The weak $(\alpha,p)$-normal derivative}: now assume that $\Omega\subseteq\R^N$ is a bounded domain. 
For general $p\in ]1,\infty[$ a weak version of the integration by parts formula \eqref{IPF-R} can be obtained as follows: let $\max\{\frac{1}{p}, \frac{p-1}{p}\}<\theta<1$ and $\alpha:=\frac{p\theta-1}{p-1}+1$. Observe that $1<\alpha<2$.
Let $u\in W^{\theta}_p (\Omega)$ be a weak solution of \eqref{eq.RFPL} with $f\in L^2(\Omega)$. We say that $u$ has a weak {\em $(\alpha,p)$-normal derivative} if there exists a function $g\in L^{\frac{p}{p-1}}(\partial\Omega)$ such that for every $v\in W^{\theta}_p (\Omega)\cap C(\overline{\Omega})$,
\begin{equation}\label{int-p-re}
\int_{\Omega} fv\; \ud x=\frac{C_{N,\theta,p}}{2} \int_{\Omega}\int_{\Omega}|u(x)-u(y)|^{p-2}\frac{(u(x)-u(y))(v(x)-v(y))}{|x-y|^{N+\theta p}}\;\ud y\;\ud x - B_{N,p,\theta} \int_{\partial\Omega}g v\ud\sigma.
\end{equation}
In that case the function $g$ is uniquely determined by  \eqref{int-p-re}. We write $B_{N,p,\theta}\mathcal N^{2-\alpha,p}_{\partial\Omega}u=g$ and call $g$ {\em the weak $(\alpha,p)$-normal derivative of $u$}. This definition makes sense as soon that  $\sigma(\partial\Omega)<\infty$. We refer to \cite{War-NODEA} for more details.
\end{enumerate}

\begin{enumerate}
\item[(1)] {\bf The Neumann boundary conditions}:\label{ex-Rpl}
we assume that the bounded open set $\Omega$ has a Lipschitz continuous boundary. 

Let  $\max\{\frac{1}{p}, \frac{p-1}{p}\}<\theta<1$ and set $\tilde{\alpha}:=\frac{p\theta-1}{p-1}+1$. Define the energy
\begin{equation*}
\E_N(u):=\frac{C_{N,\theta,p}}{{2p}}\int_{\Omega}\int_{\Omega}\frac{|u(x)-u(y)|^p}{|x-y|^{N+\theta p}}\;\ud y\;\ud x,\qquad \dom\E_N=W^{\theta}_p (\Omega)\cap L^2(\Omega).
\end{equation*}
The subgradient $\partial\E_N$ of $\E_N$ is the realization of the nonlocal operator given in \eqref{eq.RFPL} with the Neumann boundary conditions $\mathcal N_{\partial\Omega}^{2-\tilde{\alpha},p}u=0$ on $\partial\Omega$. The energy $\E_N$ is a (nonlinear) Dirichlet form in the sense given at the beginning of Section \ref{NDF}.
If $p\in [2,\infty[$, then the semigroup $S^N$ generated by the negative subgradient $-\partial\E_N$ satisfies the $L^q-L^\infty$-H\"older estimates. That is, there are constants $C_1,C_2\ge 0$ such that for every $u_0,v_0\in L^2(\Omega)\cap L^q(\Omega)$ ($q\in [1,\infty[$) and $t>0$,
\begin{equation}\label{NB1}
    \|S^N(t)u_0-S^N(t)v_0\|_\infty\le C_1\frac{|\Omega|^{\alpha}e^{C_2t}}{t^{\beta}}\|u_0-v_0\|_q^\gamma,
\end{equation}
where $(\alpha,\beta,\gamma)$ is given in \eqref{abg2} or in \eqref{abg2-b} depending the range of $\theta$. The proof of \eqref{NB1} follows the lines of the proof of Theorem \ref{theo-ultra-NEC}.

\item[(2)] {\bf The Robin boundary conditions}:
we assume that the bounded open set $\Omega$ has a Lipschitz continuous boundary.
Let $\kappa\in L^\infty(\partial\Omega)$ be such that $\kappa(x)\ge \kappa_0>0$ $\sigma$-a.e. on $\partial\Omega$ for some constant $\kappa_0$. For $p\in ]1,\infty[$, $\max\{\frac{1}{p}, \frac{p-1}{p}\}<\theta<1$ and $\tilde{\alpha}:=\frac{p\theta-1}{p-1}+1$ we consider the energy
\begin{equation*}
\E_{B,\sigma}(u) =\frac{C_{N,\theta,p}}{{2p}}\int_{\Omega}\int_{\Omega}\frac{|u(x)-u(y)|^p}{|x-y|^{N+\theta p}}\;\ud x\;\ud y+\frac 1p\int_{\partial\Omega}\kappa |u|^p\;d\sigma,\qquad \dom\E_{B}=W^{\theta}_p (\Omega)\cap L^2(\Omega).
\end{equation*}
Here, $B(x,s) := \frac{1}{p} \kappa(x) |s|^p$.
The subgradient $\partial\E_{B,\sigma}$ of $\E_{B,\sigma}$ is the realization of the nonlocal operator given in \eqref{eq.RFPL} with the Robin boundary conditions $\mathcal N_{\partial\Omega}^{2-\tilde{\alpha},p}u+\kappa|u|^{p-2}u=0$ on $\partial\Omega$. Here also the energy $\E_{B,\sigma}$ is a (nonlinear) Dirichlet form. Moreover, if $p\in [2,\infty[$, then the semigroup $S^{B,\sigma}$ generated by $-\partial\E_{B,\sigma}$ satisfies the $L^q-L^\infty$-H\"older estimates in the sense that there is a constant $C\ge 0$ such that for every $u_0,v_0\in L^2(\Omega)\cap L^q(\Omega)$ ($q\in [1,\infty[$) and $t>0$,
\begin{equation}\label{RB1}
    \|S^{B,\sigma}(t)u_0-S^{B,\sigma}(t)v_0\|_\infty\le C\frac{|\Omega|^\alpha}{t^{\beta}}\|u_0-v_0\|_q^\gamma,
    \end{equation}
    where $(\alpha,\beta,\gamma)$ is given in \eqref{abg2} or in \eqref{abg2-b} depending the range of $\theta$.

Let $\tilde{\kappa}(x,s):=\kappa(x)|s|^{p-2}s$. It has been shown in Biegert \& Warma \cite[Example 4.17]{BW-E} that $\tilde{\kappa}$ satisfies the following growth condition: there is a constant $C\in]0,1]$ such that
\begin{equation}\label{BW-E}
 C\left|\tilde{\kappa}(x,\xi-\eta) \right|\le \left|\tilde{\kappa}(x,\xi)-\tilde{\kappa}(x,\eta)\right| 
\end{equation}
for $\sigma$-a.e. $x\in\partial\Omega$ and for all $\xi,\eta\in\R$.

The proof of \eqref{RB1} follows the lines of the proof of Theorem \ref{theo-ultra} and uses heavily the growth condition \eqref{BW-E} satisfied by $\tilde{\kappa}$.
\end{enumerate}

For more details on the regional fractional Laplace operator we refer to  \cite{Ga-Wa-CPDE,GaWa2021,Gua,GuMa1,GuMa2,Wa15,War-NODEA} and their references.
\end{example}

\begin{remark}
Regarding Example \ref{ex.9.3} we make the following comments.
\begin{enumerate}
\item  As we have observed in Example \ref{ex.9.3}(i) in the integration by parts formulas \eqref{IPF-R} and \eqref{int-p-re}, if the function $u$ is smooth, say, $u\in C^1(\overline{\Omega})$, then there is no boundary boundary term, that is, the $(\alpha,p)$-normal derivative $\mathcal N^{2-\alpha,p}_{\partial\Omega}u=0$ on $\partial\Omega$. This is not the same as in the classical local case, but there is an explanation. Observe that for the regional fractional $p$-Laplace operator defined in \eqref{eq.RFPL} interactions are only allowed inside $\overline{\Omega}$. This is not the same as the fractional $p$-Laplace operator where we have interactions between $\Omega$ and $\mathbb R^N\setminus\Omega$, that allows boundary and exterior conditions. For that reason, in the integration by parts formulas for the regional fractional $p$-Laplace operator, if the function $u$ is smooth, the energy has already accumulated enough information in $\Omega\times\Omega$ so that there is no needed information left for the boundary. We refer to \cite{Gua,GuMa1,GuMa2} and their references for some probabilistic interpretations on this topic.

\item In Example \ref{ex.9.3} we have always assumed that $\max\{\frac{1}{p}, \frac{p-1}{p}\}<\theta<1$. This is not a restriction. In fact it is known that if $\Omega$ has a Lipschitz continuous boundary and $0<\theta\le \max\{\frac{1}{p}, \frac{p-1}{p}\}$, then the spaces $W_p^{\theta}(\Omega)$ and $\ringring{W}^{\theta}_p (\Omega )$ coincide with equivalent norms so that Dirichlet, Neumann and Robin boundary conditions are the same (see e.g. Grisvard \cite[Chapter 1]{Gris} for the case $p=2$ and Warma \cite{Wa15} for general $p\in ]1,\infty[$).

\item For domains $\Omega$ with a fractal geometry like the Cantor set, the von Koch curve, in the integration by parts formula \eqref{int-p-re}, one should replace the surface measure $\sigma$ with the restriction to $\partial\Omega$ of the $d$-dimensional Hausdorff measure $\mathcal H_d$ on $\partial\Omega$ where $d$ denotes the Hausdorff dimension of the boundary $\partial\Omega$. In that case the range of $\theta$ will depend on $d$. For more details on this topic we refer to \cite{Wa15}. Since this is not the main goal of this paper we will not go into details.
\end{enumerate}
\end{remark}

\begin{example}[\bf The "fractional" $p$-Laplace operator with nonlocal exterior conditions]\label{Ex-pN}
This is the case where $\Omega \subsetneq\hat{\Omega}\subseteq\R^N$ are open subsets of $\R^N$ equipped with the Lebesgue measure, that is, the situation where one explicitly has exterior conditions. 

When $\Omega\subsetneq\hat{\Omega} = \R^N$ and $p=2$, wellposedness was studied by Dipierro, Ros-Oton \& Valdinoci \cite{DiROVa17}, and later Claus \& Warma  \cite{ClWa20} proved that the corresponding semigroup on $L^2 (\Omega )$ is submarkovian. The case $p\not= 2$ (but still $\Omega\subsetneq\hat{\Omega} =\R^N$) was considered by Foghem \cite{Fog2025}. 

There is no reason to restrict to $\hat{\Omega} =\R^N$: why is the exterior of $\Omega$ only $\R^N\setminus\Omega$? Why cannot $\hat{\Omega}$ be an other set containing $\Omega$, for example a larger Riemannian manifold? Another question is: why should ${\mathbf S}$ be chosen to be equal to $(\Omega \times\hat{\Omega}) \cup (\hat{\Omega}\times \Omega)$? For example, ${\mathbf S}$ could be equal to the set $\{ (x,y)\in\hat{\Omega}\times\hat{\Omega} \st |x-y| \leq 1\}$, that is, nonlocal behaviour occurs, but there is no interaction between points at a distance larger than $1$. These were actually the motivating questions for this article. In fact, for every choice of $\hat{\Omega}$ (let $\hat{\Omega}$ be an open subset of $\R^N$, for simplicity) and for every ${\mathbf S}$ satisfying the thickness condition one obtains a semigroup $S^N$ associated with Neumann exterior conditions, even if $p$ and $k$ are fixed. Are some of these semigroups for various choices of $\hat{\Omega}$ and $\mathbf{S}$ equal? 

Of course, apart from Neumann exterior conditions, one also has the possibility of considering Robin exterior conditions and Robin boundary conditions (see Section \ref{sec.robin}).
\end{example}

\begin{example}[\bf Extensions of the preceding examples]
As far as wellposedness and the question of Dirichlet forms are concerned, one may consider more general cases.
\begin{itemize}
 \item[(1)] The function need not satisfy the condition \eqref{cond.k}, but one may consider more general functions $k$, for example functions satisfying the condition
 \begin{equation}
     C_1 \min \{ \frac{1}{|x-y|^{N+p\theta_1}}, \frac{1}{|x-y|^{N+p\theta_2}} \} \leq k(x,y) \leq C_2 \max \{ \frac{1}{|x-y|^{N+p\theta_1}}, \frac{1}{|x-y|^{N+p\theta_2}} \} .
 \end{equation}
 This includes for example,
 \[
 k(x,y) = \frac{1}{|x-y|^{N+p \theta(x,y)}}
 \]
 for some function $\theta : \hat{\Omega}\times\hat{\Omega} \to [\theta_1 , \theta_2]$ and some constants $\theta_1$, $\theta_2\in ]0,1[$, $\theta_1 \leq \theta_2$, or
 \[
 k(x,y) = \frac{1}{|x-y|^{N+p\theta_1}} + \frac{1}{|x-y|^{N+p\theta_2}} .
 \]

 \item[(2)] Also the exponent $p$ in Examples \ref{ex.3.1}, \ref{ex.6.1}, \ref{ex.6.2}, 
 could depend on $(x,y)$, that is, one could consider the kernel
 \[
 \Phi(x,y,s) = k(x,y) |s|^{p(x,y)}
 \]
 for some symmetric function $p:\hat{\Omega}\times\hat{\Omega} \to ]1,\infty [$. This kernel and its Fenchel conjugate satisfy the $\Delta_2$-condition if there exist $p_{min}$, $p_{max}\in ]1,\infty [$ such that 
 \[
 p_{min} \leq p(x,y) \leq p_{max} \text{ for almost every } (x,y)\in\hat{\Omega}\times\hat{\Omega} .
 \]

 \item[(3)] The analysis does not change when one considers Riemannian manifolds $\Omega\subseteq \hat{\Omega}$. If 
 $$\Phi (x,y,s) = \frac{1}{p} k(x,y) |s|^p 1_{\mathbf S} (x,y),$$ then $k(x,y)$ may or may not depend on the distance $d(x,y)$ determined by the Riemannian metric. The Riemannian manifolds could be boundaries of open sets in $\R^N$. The Dirichlet-to-Neumann operator on the boundary of a half space in $\R^N$ is a fractional Laplace operator on $\R^{N-1}$, and the Dirichlet-to-Neumann operator on the boundary of a general open set certainly has nonlocal behaviour. It is an open question whether these Dirichlet-to-Neumann operators are operators of the form considered here. Note that Dirichlet-to-Neumann operators and $p$-Dirichlet-to-Neumann operators are subgradients on $L^2 (\partial\Omega )$. 

 \item[(4)] In general metric spaces $\Omega\subseteq \hat{\Omega}$, the kernel $k(x,y)$ may depend on the distance $d(x,y)$ in $\hat{\Omega}$, but if $\Omega\subseteq\hat{\Omega}$ are only measured spaces, then $k$ (together with the measure $\mu$) seems to determine a kind of geometric structure.
 \end{itemize}
\end{example}

\subsection{Nonlocal operators on graphs}

\begin{example}[\bf Graphs]
We have considered graphs in Examples \ref{ex.3.2} and \ref{ex.6.3}. Graph-Laplace operators or Graph-$p$-Laplace operators are genuinely nonlocal operators. The literature is huge. We are not aware of examples when $\Omega \subsetneq\hat{\Omega}$. 
\end{example}

\subsection{Nonlocal operators on metric random walk spaces}

The example in this section only serves as information for the interested reader.

\begin{example}[\bf Metric random walk spaces]
 Let $(\hat{\Omega} ,d)$ be a Polish space (a separable, complete, metric space) equipped with its Borel $\sigma$-algebra $\mathcal{B} (\hat{\Omega})$. A {\em random walk} on $\hat{\Omega}$ is a family $m=(m_x)_{x\in\hat{\Omega}}$ of Borel probability measures $m_x$ on $\hat{\Omega}$ such that
 \begin{itemize}
     \item[(a)] the measures $m_x$ depend measurably on $x$ in the sense that for every Borel set $A\subseteq\hat{\Omega}$ and every Borel set $B\subseteq\R$ the set
     \[
     \{ x\in \hat{\Omega} \st m_x (A) \in B \} 
     \]
     is a Borel subset of $\hat{\Omega}$, and
     \item[(b)] each measure $m_x$ has finite first moment in the sense that there exists some $x_0\in\hat{\Omega}$ such that for every $x\in\hat{\Omega}$
     \[
     \int_{\hat{\Omega}} d(x_0 ,y) \; dm_x (y) < \infty ;
     \]
 \end{itemize}
 see Ollivier \cite[Definition 1]{Ol09}. A {\em metric random walk space} is a triple $(\hat{\Omega},d,m)$ where $(\hat{\Omega} ,d)$ is a Polish space and $m$ is a random walk on $\hat{\Omega}$. 

 Mazon, Solera \& Toledo \cite{MaSoTo20} and Solera \& Toledo \cite{SoTo23} consider evolution equations on measurable subsets $\Omega\subseteq\hat{\Omega}$ which are coupled with stationary equations on $\hat{\Omega}\setminus\Omega$, much as is done in this article. Existence of solutions of these evolution equations is obtained using the theory of $m$-accretive operators. Examples of metric random walk spaces include graphs and subsets of $\R^N$ (see \cite[Example 1.2]{MaSoTo20}, \cite[Examples 1.1, 1.2]{SoTo23}). 
 
 For example, if $\hat{\Omega}$ is a countable set and $(b,0)$ is a locally finite graph on $\hat{\Omega}$, then one may define
 \[
 m_x (A) := \frac{\sum_{y\in A} b(x,y)}{\sum_{y\in\hat{\Omega}} b(x,y)}  \quad (A\subseteq\hat{\Omega}) ;
 \]
 the assumption that $b$ is locally finite ensures that the two sums are finite (really, one would only need that the sums converge). Then $m=(m_x)_{x\in\hat{\Omega}}$ is a random walk on $\hat{\Omega}$. 

 Another example is the case $\hat{\Omega} = \R^N$ and $k:\R^N\times\R^N\to [0,\infty[$ is a measurable kernel such that $\int_{\R^N} k(x,y) \;\ud y =1$ for every $x\in\R^N$. Then one may define
 \[
 m_x (A) := \int_A k(x,y) \;\ud y \quad (A\subseteq\R^N \text{ measurable}) 
 \]
 and $m=(m_x)_{x\in\hat{\Omega}}$ is a random walk on $\R^N$. In this example, the assumption on the kernel excludes singular kernels as they appear for example in the fractional Laplace operator or fractional $p$-Laplace operator considered above. 

 In \cite[Example 1.2]{MaSoTo20} one can find more examples of metric random walk spaces. The abstract problems considered in \cite{SoTo23} have probably the shortcoming that the assumptions on the kernel are more restrictive than the assumptions here, but the involved operators are not necessarily subgradients ("symmetric coefficients") but also general Leray-Lions operators ("nonsymmetric coefficients"). We think that one should be aware of this parallel theory, but we do not go into more details here.
\end{example}




\end{document}